\documentclass[11pt,reqno]{amsart}
\usepackage{amssymb,mathrsfs,color,mathtools,empheq, verbatim, epstopdf}
\usepackage{pinlabel}
\mathtoolsset{showonlyrefs}
\usepackage{hyperref} 
\usepackage{enumerate, bbm}
\usepackage{tensor}
\usepackage{graphicx}
\usepackage{xcolor} 
\usepackage{slashed}
\usepackage{cite}
\usepackage[shortlabels]{enumitem}
\usepackage{geometry}
\numberwithin{equation}{section}

\newtheorem{thm}{Theorem}[section]

\newtheorem{lem}[thm]{Lemma}
\newtheorem{prop}[thm]{Proposition}

\theoremstyle{definition} 
\newtheorem{rem}[thm]{Remark}

\theoremstyle{remark}
\newtheorem{claim}[thm]{Claim}

\newcommand{\La}{\big\langle}
\newcommand{\Ra}{\big\rangle}

\newcommand{\bs}[1]{\boldsymbol{#1}}

\newcommand{\supp}{\operatorname{supp}}

\renewcommand{\ker}{\operatorname{ker}}

\newcommand{\Id}{\operatorname{Id}}

\definecolor{deepgreen}{cmyk}{1,0,1,0.5}

\newcommand{\C}{\mathbb{C}}

\newcommand{\R}{\mathbb{R}}

\newcommand{\Z}{\mathbb{Z}}

\newcommand{\al}{\alpha}
\newcommand{\be}{\beta}
\newcommand{\ga}{\gamma}
\newcommand{\de}{\delta}

\newcommand{\om}{\omega}
\newcommand{\lam}{\lambda}
\newcommand{\te}{\theta}

\newcommand{\ta}{\tau}

\newcommand{\De}{\Delta}

\newcommand{\p}{\partial}
\newcommand{\na}{\nabla}

\makeatletter

\newcommand{\Rmnum}[1]{\expandafter\@slowromancap\romannumeral #1@}
\makeatother

\newcommand{\ti}{\widetilde}

\newcommand{\Del}[1]{}

\newcommand{\mand}{{\ \ \text{and} \ \  }}

\newcommand{\mas}{{\ \ \text{as} \ \ }}

\renewcommand\Re{\mathrm{Re}\,}

\definecolor{green}{rgb}{0,0.8,0}

\renewcommand{\Re}{\mathrm{Re}}

\newcommand{\ud}{\mathrm{d}}

\newcommand{\calA}{\mathcal A}

\newcommand{\calE}{\mathcal E}
\newcommand{\calF}{\mathcal F}
\newcommand{\calG}{\mathcal G}
\newcommand{\calH}{\mathcal H}

\newcommand{\calL}{\mathcal L}

\newcommand{\calP}{\mathcal P}

\newcommand{\calR}{\mathcal R}

\newcommand{\calT}{\mathcal T}

\newcommand{\calV}{\mathcal V}

\allowdisplaybreaks[4]

\begin{document}
	\title[ ]{Smooth asymptotically Self-similar blow up for the scalar Zakharov system in dimension $3$ }
	\author{Yuchen Yin}
    \address{School of Mathematical Sciences,
University of Science and Technology of China, Hefei 230026, Anhui, China.}
\email{yuchenyin@mail.ustc.edu.cn}

\begin{abstract}
    We construct blow up solutions to the scalar Zakharov system in dimension $3$. More precisely, we construct the solutions $(\psi,n)$ defined for $t\in[0,T)$, blowing up at $t=T$ and taking the form
    \begin{align}
        \psi(t,x)&\sim\frac{1}{T-t}P\left(\frac{|x|}{\sqrt{3}(T-t)^{\frac{2}{3}}}\right)e^{i(T-t)^{-1/3}},\\
        n(t,x)&\sim\frac{1}{3(T-t)^{4/3}}N\left(\frac{|x|}{\sqrt{3}(T-t)^{\frac{2}{3}}}\right)
    \end{align}
    in the region $|x|\ll 1$. Here $(P,N)$ is the solution to a specified ODE system related to Zakharov system. 
    The solutions we construct are smooth for $0\leq t<T$. Moreover, there exist $n_0,n_1\in C^{\infty}(\R^3\backslash\{0\})$ with bounded support such that  for any $r>0$, $(\psi,n,\partial_t n)$ converge to $(0,n_0,n_1)$  in $C^{\infty}(|x|>r)$.

    We correct the profile $(P,N)$ and obtain a sequence of approximate solutions that coincide in the region $|x|\gtrsim 1$ with higher and higher accuracy.  Then we apply perturbative analysis to construct solutions with  arbitrarily high  finite  regularity. The fixed point argument in this procedure actually yields that these solutions are the same for $t$ sufficiently close to $T$, which combined with the persistence of high regularity gives the existence of smooth asymptotically self-similar blow up solutions. 
\end{abstract}
\maketitle
\section{Introduction}\label{sec:intro}
\subsection{Setting of the problem}
Consider the Zakharov system
\begin{align}\label{eq:Za}
    \begin{cases}
        i\p_t \psi+\De\psi=n\psi,\\
        \p_{tt}n-\De n=\De |\psi|^2,
    \end{cases}
\end{align}
where $\psi(t,x),n(t,x)$ are functions defined in $\R\times\R^d$. The system ~\eqref{eq:Za} is a scalar version of the effective model for Langmuir waves, also referred to as electron plasma waves, which take place in a non-magnetized or weakly magnetized plasma.  It describes the strong plasma turbulence driven by the nonlinear coupling between electron plasma waves and ion density fluctuations. See \cite{Phy1972} for details. 

The Zakharov system is Hamiltonian,  where the energy 
\begin{equation}
    \begin{aligned}
        E(\psi,n,\p_t n)_{|_t}=\frac{1}{2}\int_{\R^d}\left(|\na\psi(t)|^2+\frac{1}{2}||\na|^{-1}\p_t n(t)|^2+\frac{1}{2}|n(t)|^2+ n(t)|\psi(t)|^2\right)\ud x
    \end{aligned}
\end{equation}
and mass $\|\psi\|_{L^2}$ are conserved along the flow.  The energy space is given by $ \dot H^1(\R^d)\times L^2(\R^d)\times\dot H^{-1}(\R^d)$.

\medskip
 
The Cauchy theory of ~\eqref{eq:Za} has been studied extensively over the past decades.  For $1\leq d\leq 3$, Ozawa and Tsutsumi \cite{OT-92} established the local wellposedness for solutions with initial data in $H^2\times H^1\times L^2$ and proved the persistence of high regularity. For solutions with initial data in the energy space, Ginibre, Tsutsumi and Velo \cite{GTV-97} proved the global wellposedness for $d=1$.  The local wellposedness for $d=2,3$ was proved by  Bourgain and Colliander \cite{BC-96}.    We also refer to \cite{AA-88},  \cite{SW-86}, \cite{KPV-95}, \cite{OT-9394}, \cite{S-22}, \cite{BGHN-15}, \cite{CHN-23}.

Concerning the large time dynamic theory, we refer to  \cite{S-04}, \cite{GV-06}, \cite{OT-9394} for the study of final data problems (namely wave operators).  For $d=4$, Bejenaru, Guo, Herr and Nakanishi \cite{BGHN-15} established the small energy scattering, Guo and Nakanishi \cite{GN-21} extended the smallness to the threshold given by the ground state in the radial case. Global wellposedness and scattering for small energy solutions was established by Guo and Nakanishi \cite{GN-14} in radial case for dimension $d=3$. For large data below the ground state, scattering was showed by Guo, Nakanishi and Wang \cite{GNW-13} in radial case. Very recently,  Guo and Nakanishi \cite{GN-26} studied the global dynamics for the 3D Zakharov system with radial initial data of energy slightly above the ground state energy, proving that the initial data set splits into nine nonempty, pairwise disjoint regions in which the solutions have distinct behaviors: growup, scattering, or trapped by the ground state, as time goes in the forward or backward directions. However, there are  few results on the existence of the  finite time blow up (referred to as \textit{wave collapse}) solutions. To the author's knowledge, the only two examples are the self-similar blow up solutions in $d=2$ constructed by Glangetas and Merle in \cite{GM-94-1}, \cite{GM-94-2}, and the blow up solutions with finite energy norm for $d=4$ constructed by Krieger and Schmid \cite{KS-24-1}, \cite{KS-24-2}.

\medskip

Wave collapse phenomenon of Zakharov system in dimensions three is of fundamental importance in plasma physics. It  describes Langmuir collapse, in which the electric field becomes increasingly intense and concentrated in a shrinking region, accompanied by a deepening plasma density cavity.  Wave collapse provides a key mechanism for localizing energy in ~\eqref{eq:Za} as a conservative system.  Furthermore, the existence of a singularity motivates us to determine the conditions that leads to blow up and to reconstruct the solution before its ultimate divergence. An accurate description of the solution near the blow up time appears to be a prerequisite for investigating additional physical processes that limit the collapse. We refer to \cite{B-98} for further discussions. However, the rigorous construction of finite time blow up solutions to ~\eqref{eq:Za} in dimension three  is a long-standing open problem. It was conjectured in the 1990s by Merle \cite{M-96}, and is still listed as an open problem in  the very recent work of  Guo and Nakanishi \cite{GN-26}. See also the  survey by Cher, Czubak and Sulem \cite{CCS-16}.

We now and henceforth restrict to dimension $d=3$. We are concerned with the asymptotically self-similar blow up solutions to ~\eqref{eq:Za},  which have been observed numerically, see \cite{PSSW-91}, \cite{LPSSW-92}.  As studied in \cite{Phy1975},\cite{Phy1981}, such solutions exist only asymptotically close to collapse and have the universal form 
\begin{equation}
    \label{eq:selfsimilar-formal}
    \begin{aligned}
        \psi(t,x)&\sim\frac{1}{T-t}P\left(\frac{|x|}{\sqrt{3}(T-t)^{\frac{2}{3}}}\right)e^{i(T-t)^{-1/3}},\\
        n(t,x)&\sim\frac{1}{3(T-t)^{4/3}}N\left(\frac{|x|}{\sqrt{3}(T-t)^{\frac{2}{3}}}\right).
    \end{aligned}
\end{equation}
Define $\eta=\frac{|x|}{\sqrt{3}(T-t)^{\frac{2}{3}}}\geq 0$,  the functions $P(\eta)$ and $N(\eta)$ are radially symmetric scalar functions satisfying the coupled system of ODEs
\begin{equation}
    \label{eq:PN-dominate}
    \begin{aligned}
        \begin{cases}
            \De_{\eta}P-P-NP =0,\\
            \frac{2}{9}\left(2\eta^2N''+13\eta N'+14N\right) =\De_{\eta}\left(|P|^2     \right), 
        \end{cases}
    \end{aligned}
\end{equation}
where $\De_{\eta}=\p_{\eta}^2+2\eta^{-1}\p_{\eta}$.

 This type of blow up is referred as supersonic collapse, as the pressure term $\De n$ in ~\eqref{eq:Za} is of lower order than $\p_{tt}n$ when inserting the expressions ~\eqref{eq:selfsimilar-formal}.

As discussed in \cite{Phy1981} and proved in \cite{Ma-01}, there exits an infinite number of solutions $(\ti P_k,\ti N_k)$ to  ~\eqref{eq:PN-dominate} such that, for all integers $k\geq 1$,
\begin{align}
    0<\ti P_k(\eta)<\ti P_{k+1}(\eta)\ \ \mand \ \ti N_k(\eta)<0.
\end{align}
The profiles $\ti P_k$ decay exponentially $|\ti P_k(\eta)|\leq C_ke^{-\de\eta}$ for some $\de>0$ and for $\eta$ large, while the profiles $|\ti N_k|\leq \frac{C_k}{1+\eta^2}$ decay algebraically. For more properties of $(\ti P_k,\ti N_k)$, see Lemma ~\ref{lem:profile}.

In this paper we establish the following result.

\begin{thm}
    \label{thm:main}
    For any integer $k\geq 1$,  there exists a time $T>0$ and radial solution  $(\psi,n)$ to ~\eqref{eq:Za} such that
    \begin{equation}\label{eq:main-0}
        \begin{aligned}
             \psi\in& C_{\rm{loc}}([0,T);C^{\infty}(\R^3,\C))\cap C_{\rm{loc}}^1([0,T);C^{\infty}(\R^3,\C)),\\
            n\in&C_{\rm{loc}}([0,T);C^{\infty}(\R^3,\R))\cap C_{\rm{loc}}^1([0,T);C^{\infty}(\R^3,\R)).
        \end{aligned}
    \end{equation}
    Let $r=|x|$, $\eta=\frac{|x|}{\sqrt{3}(T-t)^{\frac{2}{3}}}$, the solution satisfies the following properties:
    \begin{enumerate}
        \item For all $0\leq t<T$, there hold
        \begin{equation}\label{eq:main-rate}
            \begin{aligned}
                \|\na_x\psi(t)\|_{L^2}\gtrsim_k (T-t)^{-\frac{2}{3}},\ \|n(t)\|_{L^2}\gtrsim_k (T-t)^{-\frac{1}{3}}.
            \end{aligned}
        \end{equation}
        In particular, the solution blows up at time $T$.
        \item  For any integer $s\geq 0$ and $\de>0$,   there holds
        \begin{equation}\label{eq:main-core}
            \begin{aligned}
                \left\|(T-t)e^{-i(T-t)^{-\frac{1}{3}}}\psi\left(t,\sqrt{3}(T-t)^{\frac{2}{3}}\eta\right)-\ti P_k(\eta) \right\|_{H_{\eta}^s} &\lesssim_{s,\de}  (T-t)^{\frac{1}{3}-\de} ,\\
                \left\|3(T-t)^{\frac{4}{3}}n\left(t,\sqrt{3}(T-t)^{\frac{2}{3}}\eta\right)-\ti N_k(\eta) \right\|_{H_{\eta}^s} &\lesssim_s  (T-t)^{\frac{1}{3}},
            \end{aligned}
        \end{equation}
        where $H_{\eta}^s$ denotes the Sobolev space $H^s$ with variable $\eta$.
        \item There exist real-valued radial functions $n_0,n_1\in C^{\infty}(\R^3\backslash\{0\})$ supported in $\{|x|\leq 2\}$ such that for any integers $s\geq 0$, $K\geq 1$, any $r_0>0$ and every $0\leq t<T$, there hold
        \begin{equation}\label{eq:main-remote}
            \begin{aligned}
                \|\psi(t)\|_{H^{s+2}(r>r_0)}+\|\p_t\psi(t)\|_{H^s(r>r_0)}\lesssim_{r_0,s,K}  & (T-t)^K,\\
                \|n(t)-n_0 \|_{H^{s+1}(r>r_0)}\lesssim_{r_0,s,K} &T-t,\\
                \|\p_tn(t)-n_1 \|_{H^{s}(r>r_0)}\lesssim_{r_0,s,K}  &T-t.
            \end{aligned}
        \end{equation}
    \end{enumerate}
\end{thm} 

\begin{rem}
    The solutions we construct are $\bf{asymptotically}$ self-similar, as the nontrivial solutions of the form ~\eqref{eq:selfsimilar-formal} with $\sim$ replaced by $=$ do not exist. This is a very different feature, compared to  the works on exact self-similar blow up solutions to focusing cubic Schr\"odinger equations,  whose ground state solution $e^{i\om t}Q_{\om}(x)$ yields a  standing wave solution $(e^{i\om t}Q_{\om},-|Q_{\om}|^2 )$ to ~\eqref{eq:Za}. We refer to \cite{KL-95}, \cite{RK-02}, \cite{DS-26}.
\end{rem}

\begin{rem}
     The stability of the solutions we constructed is an interesting topic for future study. From the numerical perspective, it was observed in \cite{LPSSW-92} that for a large class of initial data, the related solutions asymptotically display a self-similar collapse described above, with the profile being $(\ti P_1,\ti N_1)$. It would  be worthwhile to describe this phenomenon rigorously.
\end{rem}

\subsection{Strategy of the proof}

In general, we first construct solutions $(\psi,n)$ of finite regularity $H^{s+2}\times H^{s+1}$ with any prescribed $s$, then we prove that these solutions are the same for $t$ sufficiently close to $T$, which combined with the Cauchy theory for ~\eqref{eq:Za} stated in Lemma~\ref{lem:Cauchy} finishes the proof.

\bigskip

\noindent
\textit{Step 1.} Construction of solutions of finite regularity.
We  construct solutions $(\psi,n)$ of finite regularity $H^{s+2}\times H^{s+1}$ with any prescribed $s$,  for which we take a classical approach: we first construct approximate solutions at different region according to the reduced equation and glue them together to obtain a global approximate solution. Then the actual solution is obtained by a fixed point argument. We refer to \cite{KST-08}, \cite{Ga-14}, \cite{KS-24-1}, \cite{KS-24-2}, \cite{Li-26-NLS} for the applications of such strategy. 

\bigskip

\noindent
\textit{Step 1.1.} Construction of approximate solutions in the core region $r\lesssim(T-t)^{\frac{1}{2}}$.\footnote{The scale  $(T-t)^{\frac{1}{2}}$ is not essential  and can be replaced by $(T-t)^\ga$ for any   $0<\ga<\frac{2}{3}$.} 
Denote by $\tau=(T-t)^{\frac{1}{3}}$. We consider the ansatz 
\begin{align}
    \psi(t,x)&=\frac{1}{\tau^3}P\left(\tau,\eta\right)e^{i\tau^{-1}},\\
    n(t,x)&=\frac{1}{3\tau^4} N\left(\tau,\eta\right),
\end{align} 
then ~\eqref{eq:Za} turns into
\begin{equation}
    \label{eq:PN}
    \begin{aligned}
    \left\{
    \begin{array}{l}
        \De_{\eta}P-P-NP=-3i\tau P-2i\tau\eta\p_{\eta}P+i\tau^2\p_{\tau}P,\\
        \frac{2}{9}\left(2\eta^2\p_{\eta}^2N+13\eta\p_{\eta}N+14N\right) =\De_{\eta}(|P|^2)\\
        \quad +\frac{1}{3}\tau^2\left(\p_{\eta}^2N+2\eta^{-1}\p_{\eta}N\right) 
         -\frac{1}{9}\left(\tau^2\p_{\tau}^2N-4\tau\eta\p_{\tau}\p_{\eta}N-10\tau\p_{\tau}N\right).
    \end{array}
    \right.
    \end{aligned}
\end{equation}
We  take $Q(\tau,\eta)=e^{i\frac{1}{2}\tau\eta^2}P(\tau,\eta)$ to cancel the term $-3i\tau P-2i\tau\eta\p_{\eta}P$ and obtain a  system 
\begin{equation}
    \label{eq:QN}
    \begin{aligned}
    \left\{
        \begin{array}{l}
            \De_{\eta}Q-Q-NQ= -\frac{1}{2}\tau^2\eta^2 Q +i\tau^2\p_{\tau}Q ,\\
        \frac{2}{9}\left(2\eta^2\p_{\eta}^2N+13\eta\p_{\eta}N+14N\right) =\De_{\eta}(|Q|^2)\\
        \quad +\frac{1}{3}\tau^2\left(\p_{\eta}^2N+2\eta^{-1}\p_{\eta}N\right) 
         -\frac{1}{9}\left(\tau^2\p_{\tau}^2N-4\tau\eta\p_{\tau}\p_{\eta}N-10\tau\p_{\tau}N\right).
        \end{array}
        \right.
    \end{aligned}
\end{equation}
In the region $\ta\eta^2\lesssim 1$, i,e, $r\lesssim (T-t)^{\frac{1}{2}}$, we consider the approximate solution of the form
\begin{equation}
    \label{eq:QN-expansion}
    \begin{aligned}
        \left\{
        \begin{array}{rl}
            Q(\tau,\eta)&=\mathop{\sum}\limits_{m=0}^{\infty} \tau^m Q_{m}(\log\tau. \eta),\\
            N(\tau,\eta)&=\mathop{\sum}\limits_{m=0}^{\infty} \tau^m N_{m }(\log\tau. \eta),
        \end{array}
        \right.
    \end{aligned}
\end{equation}
where $Q_m(\log\tau. \eta)$ and $N_{m }(\log\tau. \eta)$  are polynomials of $\log \tau$. It is natural to take $Q_0(\log\tau,\eta)=\ti P_k(\eta)$, $N_0(\log\tau,\eta)=\ti N_k(\eta)$, where $(\ti P_k,\ti N_k)$ are solutions to ~\eqref{eq:PN-dominate}. However, the principal operator appearing in the equations for $Q_{m}$ is $\De_{\eta}-1-N_0$, which is self-adjoint and has kernels $\ti P_k,i\ti P_k$. Then we introduce a auxiliary parameter $\te$ to make the equations for $Q_m$ satisfying Fredholm conditions. More precisely, let $R(\tau,\eta)=e^{i\te(\tau)}Q(\tau,\eta)$ and take the ansatz 
\begin{equation}\label{eq:intro-RNte}
    \begin{aligned}
        \left\{
        \begin{array}{rl}
            R(\tau,\eta)&=\mathop{\sum}\limits_{m=0}^{\infty} \tau^m R_{m}(\log\tau. \eta),\\
            N(\tau,\eta)&=\mathop{\sum}\limits_{m=0}^{\infty} \tau^m N_{m }(\log\tau. \eta),\\
            \te(\tau)&=\mathop{\sum}\limits_{m=0}^{\infty} \tau^m \te_{m }(\log\tau),
        \end{array}
        \right.
    \end{aligned}
\end{equation}
we obtain the system for $(R_m.N_m,\te_{m-1})$ as in ~\eqref{eq:RNtheta-iteration}. We will construct the solutions to ~\eqref{eq:RNtheta-iteration} iteratively. We take $(R_0,N_0,\te_0)=(\ti P_k,\ti N_k,0)$. To solve the first equation of $R_m$, we need the source term being orthogonal with $\ti P_k$ and $i\ti P_k$. The former condition can be satisfied once we choose suitable $\te_{m-1}$, and the latter condition can be satisfied by taking suitable $\La R_{m-1}, \ti P_k\Ra$. Then ~\eqref{eq:RNtheta-iteration} are reduced to the nonlocal equations for $N_m$ as in ~\eqref{eq:Nn},  
\begin{align}
    \calL_m^NN_m+\calA  (N_m)= \ti F_m^{N},
\end{align}
where $\calL_{m}^N$ is a $\C$-linear operator, $\calA $ is an nonlocal $\R$-linear operator, and $\ti F_m^N$ only depends on $\{R_\ell,N_\ell,\te_{\ell-1}\}_{0\leq\ell\leq m-1}$. Inserting the ansatz that $N_m$ is a polynomial of $\log \tau$, we obtain the system, see ~\eqref{eq:Nn-1},
\begin{align}\label{eq:intro-Nmj}
    (\calL_{m,0}^N+\calA )N_{m,j}+\frac{j+1}{9}(-4\eta\p_{\eta}+2m-11)N_{m,j+1}+\frac{(j+1)(j+2)}{9}N_{m,j+2}=\ti F_{m,j},
\end{align}
where
\begin{align}
     \calL_{m,0}^N= \frac{1}{9}\left(4\eta^2\p_{\eta}^2+(26-4m)\eta\p_{\eta}+m^2-11m +28\right).
\end{align}

We have to solve it in suitable space. For $m=1$ it suffices to take $R_1=N_1=0$. For $m\geq 2$, we  impose the suitable bounds on $R_m,N_m$ and construct them inductively. The bounds we choose are ~\eqref{eq:RN-decay}, where the exponential decay of $R_m$ is inspired from the decay properties of $\ti P_k$, and the bound for $N_m$ is inspired by the identity $\calL_{m,0}^N\eta^{\frac{m-4}{2}}=\calL_{m,0}^N\eta^{\frac{m-7}{2}}=0$.  The choice of the space $X_m$  in ~\eqref{eq:Xn-def} is inspired by two observations. The first comes from the expression of $\calL_{m,0}^N$, which makes it natural to use Frobenius theorem to invert $\calL_{m,0}^N$ near the origin.  The second comes from the behavior of $R_m, N_m$ at infinity ~\eqref{eq:RN-decay}. The main difficulty here is that $\calA$ is an nonlocal operator, which prevents us to solve $N_{m,j}$ directly using the ODE techniques. The key observation is that $\calA$ is actually a compact operator from a subspace of $X_m$ to $X_m$, which allows us to use some abstract results on Fredholm theory to give possible $N_{m,j}$  without solving the equations explicitly . We will explain this in details.

We study the operator $\calL_{m,0}^N$ and $\calA$ separately.  We take the domain of $\calL_{m,0}^N$ being those functions living in $X_m$ when applying $(\eta\p_{\eta})^a$ for $a=0,1,2$.  Using Frobenius method, we find that  $\{\calL_{m',0}^N\}_{2\leq m'\leq m}$ is a family of Fredholm operators from $\ti X_m$ to $X_m$ with the parameter $m'$. We extend $\calL_{m,0}^N$ to a family of operators $\{\calL_{\lam,0}^N\}$ with the same domain for $\lam\in \C$ with  $\Re\lam<2m+6$.  This is also a family of Fredholm operator and is analytic in $\lam$.   As for $\calA$, it follows from the expression of $\calA $ with Green function that $\calA $ is a compact operator from $\ti X_m$ to $X_m$, see Claim ~\eqref{claim:An}. Then $\{\calL_{\lam,0}^N+\calA \}$ is also a family of Fredholm operator from $\ti X_m$ to $X_m$  and is analytic in $\lam$. A direct calculation yields that there exists $\lam_0$ with $\Re\lam_0\ll-1$ such that $\calL_{\lam_0,0}^N+\calA$ is invertible, which allows us to apply analytic Fredholm theory to $\{\calL_{\lam,0}^N+\calA \}$. This yields an explicit structure for  $\calL_{m,0}^N+\calA $ even when it is not invertible. More precisely, for $m\geq 2$, there exists an  integer $q_m\geq 0$ such that  
\begin{align}\label{eq:intro-L-1}
    (\calL_{\lam,0}^N+\calA)^{-1}=\sum_{\ell=0}^{\infty}(\lam-{m})^{\ell}\ti\calH_{{m},\ell}+\sum_{\ell=1}^{q_{m}}\frac{1}{(\lam-{m})^{\ell}}\calH_{{m},\ell}   
\end{align}
for some $\ti \calH_{m,\ell},\calH_{m,\ell}$ bounded from $X_m$ to $\ti X_m$. Since $\calL_{\lam,0}^N+\calA$ is a quadratic polynomial of $\lam$, ignoring the terms with positive power of $\lam-m$, we have the following formal identity for $|\lam-m|\ll1$:
\begin{align}\label{eq:intro-L-sum}
    &(\calL_{\lam,0}^N+\calA)\left(\sum_{j=0}^{\infty} \frac{j!}{(\lam-m)^j}N_{m,j} \right)\\
    \sim&\sum_{j=0}^{\infty}\frac{j!}{(\lam-m)^j}\left((\calL_{m,0}^N+\calA )N_{m,j}+\frac{j+1}{9}(-4\eta\p_{\eta}+2m-11)N_{m,j+1}+\frac{(j+1)(j+2)}{9}N_{m,j+2}\right)\\
    =&\sum_{j=0}^{\infty}\frac{j!}{(\lam-m)^j}\ti F_{m,j}
\end{align}
where we extend by zero for those $\ti F_{m,j}$ without definition. Applying $(\calL_{\lam,0}^N+\calA)^{-1}$ to ~\eqref{eq:intro-L-sum} with the expression ~\eqref{eq:intro-L-1} suggests the choice of $N_{m,j}$ in ~\eqref{eq:Nnj}. Then we verify directly that $N_{m,j}$ are actually the solutions to  ~\eqref{eq:intro-Nmj}. The bounds of $(R_m,N_m)$  follows directly from the integral expression of $\ti \calH_{m,\ell}, \calH_{m,\ell}$. We point out that we only find one   solution  to ~\eqref{eq:intro-Nmj}, as there is no need to find out all possible solutions for our construction.

After finishing the construction of $(R_m,N_m)_{m\geq 0}$, noting that the infinite sum in ~\eqref{eq:intro-RNte} may not converge, we define $(\psi_M^{\rm{co}}, n_M^{\rm{co}})$ using partial sums, see ~\eqref{eq:core-app}. The superscript `co' stands for `core'. 

\bigskip

\noindent
\textit{Step 1.2.} Construction of corrections in the remote region $r\gtrsim(T-t)^{\frac{1}{2}}$.
The approximate solution $(\psi_M^{\rm{co}}, n_M^{\rm{co}})$ is not ideal due to the slow decay or growth of $N_m$ at infinity, and we have to correct it in the remote region $\ta\eta^2\gtrsim1$.  We first study the behavior of $(\psi_M^{\rm{co}},n_M^{\rm{co}})$ for $\eta\geq1$. The behavior of $\psi_M^{\rm{co}}$ can be ignored due to the exponential decay of $R_m$ at infinity, The leading component of $N_m$ can be extracted easily using the method of characteristics and the method of variation of constant. As shown in Lemma ~\ref{lem:N-infity}, the leading component of $N_m$ is
\begin{align}
    N_m(\log\tau,\eta)\sim\sum_{j=0}^{[\frac{m}{2}]}\eta^{\frac{m-4-6j}{2}}\ti N_{m,1,j}\left(\log\tau+\frac{1}{2}\log\eta\right)+\sum_{j=0}^{[\frac{m}{2}]}\eta^{\frac{m-7-6j}{2}}\ti N_{m,2,j}\left(\log\tau+\frac{1}{2}\log\eta\right),
\end{align}
where $\ti N_{m,1,j}$ and $\ti N_{m,2,j}$ are polynomials. Inserting the expansion of $N_m$ into the expression of $n_M^{\rm{co}}$ and noting that $\log\ta+\frac{1}{2}\log\eta=\frac{1}{2}\log\frac{r}{\sqrt{3}}$, we obtain the expansion of $n_M^{\rm{co}}$ in $(t,r)$ variables, see ~\eqref{eq:N-infty-tr}. This yields that for each $r>0$, $n_M^{\rm{co}}$ and $\p_t n_M^{\rm{co}}$ have limits as $t\to T$. Since such limits tend to infinity as $r\to \infty$, we denote by $ n_{0,M}^{\rm{re}}$ and $n_{1,M}^{\rm{re}}$ the truncated versions of them, see ~\eqref{eq:nre-initial}. Since $\psi_{M}^{\rm{co}}$ decays exponentially at infinity, we take $\psi^{\rm{re}}_M\equiv 0$ and evolve $(n_{0,M}^{\rm{re}},n_{1,M}^{\rm{re}})$ at time $T$ backward along linear wave equation, which yields a correction $n_{M}^{\rm{re}}$ for $r>T-t$. Since $\ta\eta^2\gtrsim1$ is equivalent to $r\gtrsim (T-t)^{\frac{1}{2}}$, for $t$ sufficiently close to $T$, the domain of $n_M^{\rm{re}}$ covers it.

\bigskip

We point out that, if we only want to construct solutions to ~\eqref{eq:Za} with arbitrarily high finite regularity, it suffices to take a specific $M$ sufficiently large and take the approximation solution as
\begin{align}\label{eq:intro-ap}
        \hat \psi^{\rm{ap}}_M(t,r)&=\chi\left(\frac{r}{\sqrt{3}(T-t)^{\frac{1}{2}}}\right)\psi_M^{\rm{co}}(t,r),\\
        \hat n^{\rm{ap}}_M(t,r)&=\chi\left(\frac{r}{\sqrt{3}(T-t)^{\frac{1}{2}}}\right)n_M^{\rm{co}}(t,r)+\left(1-\chi\left(\frac{r}{\sqrt{3}(T-t)^{\frac{1}{2}}}\right)\right)n_M^{\rm{re}}(t,r).
\end{align}
The arguments in Section ~\ref{subsec:E}, Section ~\ref{subsec:linear}, Section ~\ref{subsec:nl} and Section ~\ref{subsec:finite-regularity} are still applicable to $(\hat \psi_{M}^{\rm{ap}},\hat n_M^{\rm{ap}})$. To show the motivation of our choice of $n^{\rm{re}}$ in Lemma~\ref{lem:psi-n-re}, we first describe the perturbative arguments in Section ~\ref{sec:fixed-point} and then turn back to this choice.

\bigskip

\noindent
\textit{Step 1.3.} Construction of the actual solution by perturbative analysis.
Suppose that we have constructed a desired approximate solution $(\psi^{\rm{ap}},n^{\rm{ap}})$ to ~\eqref{eq:Za}, we want to correct it to an actual solution. We write $(\psi,n)=(\psi^{\rm{ap}},n^{\rm{ap}})+(u,v)$, then we obtain the system for $(u,v)$ as
\begin{align}\label{eq:intro-uv}
     \left\{
         \begin{array}{rl}
              i\p_t u+\De_r u-n^{\rm{ap}}u-\psi^{\rm{ap}}v &=uv-\calE_{\psi}^{\rm{ap}},  \\
              \p_{tt}v-\De_r v-2\Re\De_r(\overline{\psi^{\rm{ap}}}u) &=\De_r(|u|^2)-\calE_{n}^{\rm{ap}}, 
         \end{array}
    \right.
\end{align}
where $\calE_{\psi}^{\rm{ap}},\calE_n^{\rm{ap}}$ are the errors generated by  $(\psi^{\rm{ap}},n^{\rm{ap}})$. Inspired by the estimates on the errors in Lemma~\ref{lem:E}, we introduce a function space $\calH_K^s$ that contains the functions $(f,g)$ whose $H^{s+2}\times H^{s+1}$ norm can be bounded by $(T-t)^K$ uniformly in time, see ~\eqref{eq:HKs-norm}. We solve $(u,v)$ in this space. We first study the linear system
\begin{align}
     \left\{
         \begin{array}{rl}
              i\p_t u+\De_r u-n^{\rm{ap}}u-\psi^{\rm{ap}}v &=F,  \\
              \p_{tt}v-\De_r v-2\Re\De_r(\overline{\psi^{\rm{ap}}}u) &=G, 
         \end{array}
    \right.
\end{align}
for suitable $(F,G)$. The difficulties are that the terms $(\psi^{\rm{ap}},n^{\rm{ap}})$ becomes singular as $t\to T^{-}$, and that we require the solutions vanish in high order rather than just being zero at the time $T$.  We apply the method of backwards in time integration along the flow introduced in \cite{M-90} and derive an energy method in  rescaled Sobolev spaces to obtain the existence and uniqueness results,  see Lemma~\ref{lem:UV}, where we require $s$ to be even to avoid taking the square root of certain operator when defining the energy.   The difference of the nonlinear terms in ~\eqref{eq:intro-uv} for $(u_1,v_1)$ and $(u_2,v_2)$ can be estimated easily, see ~\eqref{eq:nl-estimate}. Thus a standard fixed point argument yields the existence and uniqueness of $(u,v)$ in the closed ball of $\calH_K^s$ for any even integer $s$ and $K$ sufficiently large once the approximate solution is sufficiently accurate, i.e. $M$ is sufficiently large for our $(\psi_M^{\rm{ap}},n_M^{\rm{ap}})$, and $t$ is sufficiently close to $T$.

\bigskip

\noindent
\textit{Step 2.} Final choice of the approximate solutions and the construction of smooth blow up solutions.
To get a smooth solution, since ~\eqref{eq:Za} persists high regularity, it is natural to expect that the actual solutions $(\hat \psi_M,\hat n_M)$ obtained from $(\hat \psi_M^{\rm{ap}},\hat n_M^{\rm{ap}})$ coincide for $t$ sufficiently close to $T$. However, the difference of $\hat n_M^{\rm{ap}}$ does not tend to zero in the region $r\sim 1$, which requires us to find a universal correction $n^{\rm{re}}$ in the remote region. Obviously we cannot directly take the infinite sum, and we use the arguments in the proof of Borel's theorem concerning the existence of a smooth function with prescribed Taylor series at $0$, that is, for each term in the sum, we multiply it with a cutoff function $\chi(r/r_m)$ for some $0<r_m\ll1$ with $r_m\to 0$, then the sum is locally finite for any $r>0$ and approximates  $n_{0,M}^{\rm{re}},n_{1,M}^{\rm{re}}$ well, see Lemma~\ref{lem:psi-n-re} for details.   Then we propagate it backward in time at $T$ along the linear wave equation and get the desired $n^{\rm{re}}$. The final approximate solution $(\psi_M^{\rm{ap}},n_M^{\rm{ap}})$ is given by ~\eqref{eq:intro-ap} with $n_M^{\rm{re}}$ replaced by $n^{\rm{re}}$.  The estimates of the difference of the approximate solutions with different $M$ are shown in Lemma ~\ref{lem:ap-diff}.  Then utilizing the fixed-point arguments in Proposition~\ref{prop:fixed-point}, we prove that the  actual solutions $(\psi_M,n_M)$ obtained from $(  \psi_M^{\rm{ap}},   n_M^{\rm{ap}})$ coincide for $t$ sufficiently close to $T$, which gives the desired construction, see Proposition~\ref{prop:smooth} for details. The estimates  ~\eqref{eq:main-rate}, ~\eqref{eq:main-core} and  ~\eqref{eq:main-remote} follows directly from the expression of $(\psi_M^{\rm{ap}},   n_M^{\rm{ap}})$ and that the related $(u,v)\in \calH_K^s$.

\subsection{Notation}

Without further mention, all the functions we consider are radial. In the rest of this paper, We use $H^s$, $W^{s,\infty}$ to denote the Sobolev spaces of radial functions in $\R^3$.

For a function $f$ of variable $r=|x|$ or of variable $\eta=\frac{|x|}{\sqrt{3}(T-t)^{\frac{2}{3}}}$, we regard them as radial functions defined in $\R^3$. We denote by 
\begin{align}
    \na_x f(r)=\left(\frac{x_1}{r}\p_rf(r),\frac{x_2}{r}\p_rf,\frac{x_3}{r}\p_rf\right)
\end{align}
and by
\begin{align}
    \p_x^a f(r)=\p_{x_1}^{a_1}\p_{x_2}^{a_2}\p_{x_3}^{a_3} f(r),\quad |a|=a_1+a_2+a_3
\end{align}
for  $a=(a_1,a_2,a_3)\in \Z_{\geq 0}^3$

For functions $f,g$ defined in $\R^3$, define
\begin{align}
    \La f,g \Ra=\Re\int_{\R^3}\bar{f} g\,\ud x.
\end{align}

We use $\chi(r)=\chi(|x|)$ to denote a smooth radial function in $\R^3$ that satisfies
\begin{align}\label{eq:chi}
    \chi(r)=\begin{cases}
        1,&0\leq r\leq1,\\
        0,&r\geq 2.
    \end{cases},\ 0\leq\chi\leq1.
\end{align}

For $a,b>0$, we write $a\lesssim b$ (resp. $a\gtrsim b$), if there exists a constant $C>0$ depending only on $k$, such that $a\leq Cb$  (resp. $a\geq Cb$).  We write $a\lesssim_p b$ or $a\gtrsim_p b$ if the implicit constant also depends on other parameters $p$.

\noindent
\textbf{Acknowledgments.}   The author would like to thank his advisor, Prof. Lifeng Zhao, for his valuable advice  on this work. Y.Yin is supported by NSFC Grant of China No. 12271497.

\section{Preliminaries}\label{sec:pre}

\subsection{Cauchy theory}
In this section we record the part of the Cauchy theory developed in \cite{OT-92} that is necessary for our construction.

\begin{lem}[\rm{[Theorem 1.1]{\cite{OT-92}}}]
    \label{lem:Cauchy} 
    For $(\psi_0,n_0,n_1)\in H^2\times H^1\times L^2$, then for some $ T>0$ depending on $\|\psi_0\|_{H^2},\|n_0\|_{H^1},\|n_1\|_{L^2}$, there exist  unique solutions $(\psi,n)$ to ~\eqref{eq:Za} with initial data $(\psi(0),n(0),\p_tn(0))=(\psi_0,n_0,n_1)$ such that
    \begin{equation}\label{eq:lem-regularity}
        \begin{aligned}
            \psi\in&C([0,  T];H^2)\cap C^1([0,  T];L^2),\\
            n\in & C([0,  T];H^1)\cap C^1([0,  T];L^2).
        \end{aligned}
    \end{equation}
    Let $s\geq 1$ be an odd integer. If $ (\psi_0,n_0,n_1)\in H^{s+2}\times H^{s+1}\times H^s$, then the solutions $(\psi,n)$ to ~\eqref{eq:Za} given above satisfy 
    \begin{equation}\label{eq:lem-high-regularity}
        \begin{aligned}
            \psi\in&C([0,  T];H^{s+2})\cap C^1([0,  T];H^s),\\
            n\in & C([0,  T];H^{s+1})\cap C^1([0, T];H^s).
        \end{aligned}
    \end{equation}
\end{lem}
\begin{rem}
    The argument  in Proposition 3.3 \cite{OT-92} actually yields that if solutions $(\psi,n)$ to ~\eqref{eq:Za} satisfying ~\eqref{eq:lem-regularity} for some $T>0$ with initial data $ (\psi_0,n_0,n_1)\in H^{s+2}\times H^{s+1}\times H^s$ for some odd integer $s$, then $(\psi,n)$ also satisfy ~\eqref{eq:lem-high-regularity}
\end{rem}

\subsection{Existence of profiles}

In this section we collect the results of the profiles $(\ti P_k,\ti N_k)$ in \cite{Ma-01}. For integers $k\geq 1$, define $\al_k=\frac{1}{3}\sqrt{2k(4k+3)}>0$.

\begin{lem}[\rm{[ Proposition 1, Theorem 1, Proposition 2, Proposition 3,  Theorem 2]\cite{Ma-01}}]
    \label{lem:profile}
    For each integer $k\geq 1$, there exits a pair of radially symmetric solutions $(\ti P_k,\ti N_k)$ to ~\eqref{eq:PN-dominate} with the following properties:
    \begin{enumerate}
        \item 
        \begin{align}
            \al_k<\ti P_k(0)<\al_{k+1},\ \ti N_k(0)=\frac{9\ti P_k^2(0)}{14-9\ti P_k^2(0)}.
        \end{align}
        \item $\ti P_k,\ti N_k\in  C^{\infty}$ are analytic near the origin and can be written as convergent power series of $\eta^2$.
        \item For each $\eta>0$,
        \begin{align}
            \ti P_k(\eta)>0,\ \ti P'_k(\eta)<0,\ \ti N_k(\eta)<0.
        \end{align}
        \item For any integers $j\geq 0$, there exit constants $\de_j>0$, $\ti C_{k,j},C_{k,j}>0$ such that
        \begin{align}\label{eq:PN-decay}
            |\ti P_k^{(j)}(\eta)|\leq \ti C_{k,j}e^{-\de_j\eta},\ |\ti N_k^{(j)}(\eta)|\leq C_{k,j}\frac{1}{(1+\eta)^{2+j}}.
        \end{align}
        \item There holds
        \begin{align}\label{eq:tiNk-P}
            \ti N_k(\eta)=-\frac{3}{2}\frac{\ti P_k^2(0)}{\eta^2}-\frac{9}{8\eta^{\frac{7}{2}}}\int_0^{\eta}s^{\frac{1}{2}}\ti P_k^2(s)\ud s+\frac{9\ti P_k^2(\eta)}{4\eta^2}.
        \end{align}
        \item For any integer $k\geq 1$ and $\eta\geq 0$, there holds
        \begin{align}
            0<\ti P_k(\eta)<\ti P_{k+1}(\eta)\ \ \mand \ \ti N_k(\eta)<0.
        \end{align}
    \end{enumerate}
\end{lem}
\begin{rem}
    Notice that the system ~\eqref{eq:PN-dominate} yields that $\p_{\eta}^2(\eta\ti P_k)=(\ti N_k+1)\eta\ti P_k $. Since $\ti N_k(\eta)=O((1+\eta)^{-2})$, a direct ODE argument yields for $\eta\geq 1$,
        \begin{align}\label{eq:Pk-asymp}
            \ti P_k(\eta) \sim_k \eta^{-1}e^{-\eta}, \ |\p_{\eta}\ti P_k(\eta)|\lesssim_k  \eta^{-1}e^{-\eta},
        \end{align}
        and the constants $\de_j$ in ~\eqref{eq:PN-decay} can be chosen to satisfy $0<\de_j<1$.
\end{rem}

\subsection{Spectral properties}

In this section we establish spectral properties for linear operators appearing in the following sections.

Fix an integer $k\geq 1$, consider the operator
\begin{align}
    \calL_k=\De_\eta-1-\ti N_k:\ {\rm{D}}(\calL_k)\to L^2(\R^3,\C)
\end{align}
with the domain
\begin{align}
    {\rm{D}}(\calL_k)=H^2(\R^3,\C)\subset L^2(\R^3,\C).
\end{align}

\begin{lem}
    \label{lem:Lk}
    There holds
    \begin{align}\label{eq:ker-Lk}
        \ker \calL_k=\mathrm{Span}_{\R} \left(\ti P_k, i\ti P_k\right).
    \end{align}
    Moreover, for any integer $m\geq 0$ and any function $f\in H^m$ such that $\La f, \ti P_k\Ra=\La f, i\ti P_k\Ra=0 $, there exists unique $u\in H^{m+2}$ such that 
    \begin{align}\label{eq:Lku-f}
        \calL_k u=f,\ \La u, \ti P_k\Ra=\La u, i\ti P_k\Ra=0 ,
    \end{align}
    and  
    \begin{align}\label{eq:coercive-Lk}
        \|u\|_{H^{m+2}}\lesssim_{k,m} \|f\|_{H^m}
    \end{align}
\end{lem}

\begin{proof}
    Since $\ti N_k$ is real-valued and is analytic near $0$, Frobenius theorem yields  ~\eqref{eq:ker-Lk}. The inequalities ~\eqref{eq:PN-decay} yields that multiplying $\ti N_k$ is a compact operator  from $H^2$ to  $L^2$, therefore $\calL_k$ is a Fredholm operator. The standard Fredholm theory yields the existence and uniqueness of the function $u\in H^2$ satisfying ~\eqref{eq:Lku-f}.  By Weyl's theorem  we obtain that $\sigma_{\rm{ess}}(\calL_k)=\sigma_{\rm{ess}}(\De_{\eta}-1)=(-\infty,-1]$, therefore $0\in \sigma_{\rm{disc}}(\calL_k)$,  which yields $\|u\|_{L^2}\lesssim \|f\|_{L^2}$. This combined with standard elliptic estimates yields     the estimate ~\eqref{eq:coercive-Lk}. 
\end{proof}

\section{Analysis in the core region}\label{sec:core}

In this section we construct the approximate solution for ~\eqref{eq:Za} in the core region $\{r\lesssim (T-t)^{\frac{1}{2}}\}$.  Fix $k\geq 1$ to be an integer. Recall that $\tau=(T-t)^{\frac{1}{3}}$, $\eta=\frac{r}{\sqrt{3}(T-t)^{\frac{2}{3}}}$.

\begin{prop}
    \label{prop:core}
    For every integer  $s\geq0$, there exists $M_{s}\in \Z_{\geq 0}$ such that for any integer $M\geq M_{s}$, there exists smooth  functions $(\psi^{\rm{co}}_{M} ,n^{\rm{co}}_{M} )$ defined for all $t<T$ sufficiently close to $T$ satisfying the following properties.
    \begin{enumerate}
        \item For any $L>0$ and any integer $s'>0$, there hold
        \begin{equation}\label{eq:core-converge}
            \begin{aligned}
                \tau^3e^{-i\tau^{-1}}\psi^{\rm{co}}_{M}(t,\sqrt{3}\tau^2\eta)&\to\ti P_k(\eta),\\
            3\tau^4n^{\rm{co}}_{M}(t,\sqrt{3}\tau^2\eta)&\to\ti N_k(\eta)
            \end{aligned}
        \end{equation}
        both in $C_{\eta}^{s'}([0,L(T-t)^{-\frac{1}{6}}))$.
        \item Define
        \begin{equation}\label{eq:E-core}
            \begin{aligned}
                \calE^{\rm{co}}_{\psi,M}=&i\p_t\psi^{\rm{co}}_{M}+\De\psi^{\rm{co}}_{M}-n^{\rm{co}}_{M}\psi^{\rm{co}}_{M},\\
            \calE^{\rm{co}}_{n,M}=&\p_{tt}n^{\rm{co}}_{M}-\De n^{\rm{co}}_{M}-\De|\psi^{\rm{co}}_{M}|^2.
            \end{aligned}
        \end{equation}
        Then for any $a_1,a_2\in \Z_{\geq 0}$ with $a_1+a_2\leq s,$  we have
        \begin{align}\label{eq:E-core-decay}
             \sup_{r\leq2\sqrt{3} (T-t)^{1/2}}\left(\left|\p_{t}^{a_1}\p_r^{a_2}\calE^{\rm{co}}_{\psi,M}(t,r)\right|+ \left|\p_{t}^{a_1}\p_r^{a_2}\calE^{\rm{co}}_{n,M}(t,r)\right|\right)\lesssim_{M,s} (T-t)^{\frac{1}{4}M-3-2s}.
        \end{align}
    \end{enumerate}
\end{prop}

The proof of Proposition ~\ref{prop:core} will be shown in the end of this section. 

\bigskip

We consider the ansatz
\begin{align}
    \psi(t,r)=\frac{1}{\tau^3}P\left(\tau,\eta\right)e^{i\tau^{-1}},
    n(t,r)=\frac{1}{3\tau^4} N\left(\tau,\eta\right).
\end{align} 
Let $R(\tau.\eta)=e^{i\te(\tau)+i\frac{1}{2}\tau\eta^2}  P(\tau,\eta)$  with $\te(\tau)\in \R$, then ~\eqref{eq:PN} turns into
\begin{equation}
    \label{eq:RN}
    \begin{aligned}
    \left\{
        \begin{array}{l}
            \De_{\eta}R-R-NR= -\frac{1}{2}\tau^2\eta^2 R  +i\tau^2\p_{\tau}R+\tau
            ^2\te'R,\\
        \frac{2}{9}\left(2\eta^2\p_{\eta}^2N+13\eta\p_{\eta}N+14N\right) =\De_{\eta}(|R|^2)\\
        \quad +\frac{1}{3}\tau^2\left(\p_{\eta}^2N+2\eta^{-1}\p_{\eta}N\right) -\frac{1}{9}\left(\tau^2\p_{\tau}^2N-4\tau\eta\p_{\tau}\p_{\eta}N-10\tau\p_{\tau}N\right).
        \end{array}
        \right.
    \end{aligned}
\end{equation}

We seek for an approximate solution to ~\eqref{eq:RN} of the form 
\begin{equation}
    \label{eq:RNtheta-expansion}
    \begin{aligned}
        \left\{
        \begin{array}{rl}
            R(\tau,\eta)&=\mathop{\sum}\limits_{m=0}^{\infty} \tau^m R_{m}(\log\tau. \eta),\\
            N(\tau,\eta)&=\mathop{\sum}\limits_{m=0}^{\infty} \tau^m N_{m }(\log\tau. \eta),\\
            \te(\tau)&=\mathop{\sum}\limits_{m=0}^{\infty} \tau^m \te_{m }(\log\tau).
        \end{array}
        \right.
    \end{aligned}
\end{equation}
Here $R_m(\log\tau. \eta)$, $N_{m }(\log\tau. \eta)$  and $\te_{m }(\log\tau)$ are polynomials of $\log \tau$, and $N_m,\te_m$ are real-valued. Denote by $\sigma=\log\tau$. For $m<0$, we adopt the convention $R_m=N_m=\te_m\equiv0$.

We take 
\begin{align}\label{eq:initial}
    R_{0}(\sigma,\eta)= \ti P_k(\eta),\ N_{0}(\sigma,\eta)=\ti N_k(\eta),\ \te_0(\sigma)=0.
\end{align}
Inserting ~\eqref{eq:RNtheta-expansion} and ~\eqref{eq:initial}  into ~\eqref{eq:RN} and comparing the coefficients of $\tau$, we obtain that the coefficient of $\tau^0$ vanishes due to ~\eqref{eq:PN-dominate}, and for any integer $m\geq 1$,there holds
\begin{equation}
    \label{eq:RNtheta-iteration}
    \begin{aligned}
        \left\{
        \begin{array}{rl}
             \calL^R R_m-R_0N_m-R_0(\p_{\sigma}+m-1)\te_{m-1}&=F_m^R,\\
             \calL_m^N N_m-2\Re \De_{\eta}\left(\overline{ R_0} R_m\right)&=F_m^N.
        \end{array}
        \right.
    \end{aligned}
\end{equation}
Here 
\begin{equation}
    \begin{aligned}
            \calL^R&=  \De_{\eta}-1-N_0,\\
            \calL_m^N  &=   \frac{2}{9}\left(2\eta^2\p_{\eta}^2+13\eta\p_{\eta}+14\right)  +\frac{1}{9}\big[ \p_{\sigma}^2+(2m-11)\p_{\sigma}+m^2-11m-4\eta\p_{\sigma}\p_{\eta}-4m\eta\p_{\eta}\big],
    \end{aligned}
\end{equation}
and
\begin{equation}
    \begin{aligned}
            F_m^R &=  i(\p_{\sigma}+m-1)R_{m-1} -\frac{1}{2}\eta^2R_{m-2}+\mathop{\sum}\limits_{j=1}^{m-1} N_j R_{m-j}  + \mathop{\sum} \limits_{j=0}^{m-2} R_{m-1-j}(\p_{\sigma}+j)\te_j,   \\
            F_m^N  &= \frac{1}{3}\De_{\eta} N_{m-2}+\mathop{\sum}\limits_{j=1}^{m-1}   \De_{\eta}\left(\overline{R_j} R_{m-j}\right).
    \end{aligned}
\end{equation}

Under the choice ~\eqref{eq:initial}, we have  $\calL^R=\calL_k$.  For $m\geq 1$, we impose the following orthogonal conditions:
\begin{align}
    \La  R_0N_m+R_0(\p_{\sigma}+m-1)\te_{m-1}+F_m^R,\ti P_k \Ra&=0,\label{eq:ortho-1}\\
    \La  R_0N_m+R_0(\p_{\sigma}+m-1)\te_{m-1}+F_m^R,i\ti P_k \Ra&=0.\label{eq:ortho-2}
\end{align}
It is easy to see that ~\eqref{eq:ortho-2} holds for $m=1$ under the choice ~\eqref{eq:initial}.

\begin{lem}\label{lem:ortho12}
    Under the choice of ~\eqref{eq:initial}, for any integer $m_0\geq 2$. if there exists $\{R_m,N_m\}_{0\leq m\leq m_0-1}$ and $\{\te_{m-1}\}_{1\leq m\leq m_0-1}$ of the form ~\eqref{eq:RNtheta-expansion} such that ~\eqref{eq:RNtheta-iteration} holds for $1\leq m\leq {m_0}-1$, then  
     \begin{align}\label{eq:ortho-2-equiv}
        \La  R_0N_{m_0}+R_0(\p_{\sigma}+m_0-1)\te_{{m_0}-1}+F_{m_0}^R,i\ti P_k \Ra =\frac{\p_{\sigma}+{m_0}-1}{2} \left(\sum_{j=1}^{{m_0}} \La  R_{j-1}   , R_{{m_0}-j}\Ra  \right).
    \end{align}
\end{lem}

\begin{proof}[Proof of Lemma ~\ref{lem:ortho12}]
      We immediately have ~\eqref{eq:ortho-1} and ~\eqref{eq:ortho-2} hold for $1\leq m\leq {m_0}-1$. Since $N_m\in \R$ for all $1\leq m\leq {m_0}-1$, using ~\eqref{eq:RNtheta-iteration} for $1\leq j\leq {m_0}-1$ we have
    \begin{align}
        &\La  R_0N_{m_0}+R_0(\p_{\sigma}+{m_0}-1)\te_{{m_0}-1}+F_{m_0}^R,i\ti P_k \Ra\\
        =&\La  R_0N_{m_0}+R_0(\p_{\sigma}+{m_0}-1)\te_{{m_0}-1}+F_{m_0}^R,i R_0 \Ra\\
        =&(\p_{\sigma}+{m_0}-1)\left(\La  R_{{m_0}-1},R_0  \Ra\right)-\frac{1}{2}\La \eta^2R_{{m_0}-2},iR_0 \Ra\\
        &+\sum_{j=1}^{{m_0}-1}\La R_{{m_0}-j}\left( N_j+ (\p_{\sigma}+j-1)\te_{j-1} \right),iR_0    \Ra\\
        =&(\p_{\sigma}+{m_0}-1)\left(\La  R_{{m_0}-1},R_0  \Ra\right)-\frac{1}{2}\La \eta^2R_{{m_0}-2},iR_0 \Ra\\
        &+\sum_{j=1}^{{m_0}-1}\La R_{{m_0}-j}\left(\calL^RR_j-F_j^R \right),i    \Ra.
    \end{align}
    Since $\calL^RR_j-F_j^R=R_0\left(N_j+ (\p_{\sigma}+j-1)\te_{j-1}\right)\in\R$ and $\calL^R$ is self-adjoint, we have
    \begin{align}
        &\sum_{j=1}^{{m_0}-1}\La R_{{m_0}-j}\left(\calL^RR_j-F_j^R \right),i    \Ra\\
        =&-\sum_{j=1}^{{m_0}-1}\La\calL^RR_j-F_j^R ,iR_{{m_0}-j}  \Ra\\
        =&-\frac{1}{2}\sum_{j=1}^{{m_0}-1} \left(\La\calL^RR_j ,iR_{{m_0}-j} \Ra+ \La\calL^RR_{{m_0}-j} ,iR_{ j} \Ra \right) + \sum_{j=1}^{{m_0}-1}\La F_j^R ,iR_{{m_0}-j}  \Ra\\
        =&\sum_{j=1}^{{m_0}-1}\Big< i(\p_{\sigma}+j-1)R_{j-1}-\frac{1}{2}\eta^2R_{j-2}+\sum_{\ell=1}^{j-1}R_{j-\ell}\left(N_{\ell}+(\p_{\sigma}+\ell-1)\te_{\ell-1}\right)  ,iR_{{m_0}-j}\Big>  \\
        =&\sum_{j=1}^{{m_0}-1} \La  (\p_{\sigma}+j-1)R_{j-1}   , R_{{m_0}-j}\Ra  -\frac{1}{2}\La\eta^2R_0,iR_{{m_0}-2}\Ra-\frac{1}{2}\sum_{j=3}^{{m_0}-1}\La\eta^2R_{j-2},iR_{{m_0}-j}\Ra\\
        &+ \sum_{\ell=1}^{{m_0}-2}\sum_{j=\ell+1}^{{m_0}-1}\Big< \left(N_{\ell}+(\p_{\sigma}+\ell-1)\te_{\ell-1}\right) R_{j-\ell} ,iR_{{m_0}-j}\Big>\\
        =&\sum_{j=1}^{{m_0}-1} \La (\p_{\sigma}+j-1) R_{j-1}   , R_{{m_0}-j}\Ra  -\frac{1}{2}\La\eta^2R_0,iR_{{m_0}-2}\Ra.
    \end{align}
    Therefore
    \begin{align}
        \La  R_0N_{m_0}+R_0(\p_{\sigma}+{m_0}-1)\te_{{m_0}-1}+F_{m_0}^R,i\ti P_k \Ra
        &=\sum_{j=1}^{{m_0}} \La (\p_{\sigma}+j-1) R_{j-1}   , R_{{m_0}-j}\Ra  \\
        &=\frac{\p_{\sigma}+{m_0}-1}{2} \left(\sum_{j=1}^{{m_0}} \La  R_{j-1}   , R_{{m_0}-j}\Ra  \right).
    \end{align}

\end{proof}

Now we solve the iteration system ~\eqref{eq:RNtheta-iteration}. 
\begin{lem}
   \label{lem:iteration-solve} 
   Under the choice of ~\eqref{eq:initial}, there exists $\{R_m,N_m,\te_m\}_{m\in \Z_{\geq 1}}$ such that
     \begin{enumerate}
        \item  $R_m$ and $N_m$ are polynomials of $\sigma$
        \begin{align}
            R_m(\sigma,\eta)=\sum_{j=0}^{L_m^R}R_{m,j}(\eta)\sigma^j,\ N_m(\sigma,\eta)=\sum_{j=0}^{L_m^N}N_{m,j}(\eta)\sigma^j,
        \end{align}
        where $0\leq L_m^R,L_m^N<\infty$. $R_{m,j}$ and $N_{m,j}$ are smooth   function and are analytic in $\eta^2$ near $0$. $N_{m,j}$ are real-valued. $\te_m(\sigma)$ are real polynomials of $\sigma$.
        
         \item For any integer $m\geq 1$, ~\eqref{eq:RNtheta-iteration} holds.
         \item For any integer $m\geq 1$, $R_m$ admits a decomposition of the form
         \begin{align}\label{eq:Rm-decomposition}
             R_m(\sigma,\eta)=R_m^{\perp}(\sigma,\eta)+a_m(\sigma)\ti P_{k}(\eta),
         \end{align}
         where $\La R_m^{\perp},\ti P_{k} \Ra=\La R_m^{\perp}, i\ti P_{k} \Ra=0$, and 
         \begin{align}
             a_m(\sigma)=-\frac{1}{2\|R_0\|_{L^2}^2}\sum_{j=1}^{m-1}\La R_{j },R_{m-j} \Ra.
         \end{align}
         \item We have the following decay properties: for any $0<\de<1$ and any integers $s\geq 0$, and any $m\geq 0$ and $j\geq 0$, there hold
         \begin{align}\label{eq:RN-decay}
             |\p_{\eta}^{s} R_{m,j}(\eta)|\lesssim_{m,j,s,\de}e^{-\de\eta},\ |\p_{\eta}^{s} N_{m,j}(\eta)|\lesssim_{m,j,s,\de} (1+\eta)^{-s+\frac{m-4}{2}+\de}.
         \end{align}
     \end{enumerate}
\end{lem}
\begin{proof}[Proof of Lemma ~\ref{lem:iteration-solve}]
     We argue by induction. Inspired by the Frobenius theory, we slightly strengthen the induction hypothesis. Let $\rho>0$ such that $\ti P_k$ and $\ti N_k$ are analytic in $\eta^2$ in a neighborhood of $0$ strictly larger than $\{\eta:0\leq \eta\leq \rho\}$  and their Taylor series at $0$ converge to themselves and  absolutely converge in the same region. Such $\rho$ exists due to Lemma ~\ref{lem:profile}.  We will also require that for $0\leq \eta\leq \rho$, $R_{m,j}$ and $N_{m,j}$ have power series expressions of the form
     \begin{align}
         R_{m,j}(\eta)=&\sum_{\ell=0}^{\infty}R_{m,j,\ell}\eta^{2\ell},\ N_{m,j}(\eta)=\sum_{\ell=0}^{\infty} N_{m,j,\ell}\eta^{2\ell},
     \end{align}
     with
     \begin{align}\label{eq:assumption-strengthen}
         \sum_{\ell=0}^{\infty}(1+\ell)^2\left(|R_{m,j,\ell}|+|N_{m,j,\ell}|\right)\rho^{2\ell}<\infty.
     \end{align}
     In the case $m=1$, it suffices to take $R_1=N_1\equiv 0$. Suppose that we have constructed $R_m,N_m,\te_{m-1}$ for $m\leq {m_0}-1$, where ${m_0}\geq 2$ is an integer.  We first express $\te_{{m_0}-1}$ and $R_{m_0}$ with $N_{m_0}$ and  $\{R_m,N_m,\te_{m-1}\}_{0\leq m\leq {m_0}-1}$, then we will get a nonlocal equations for $N_{m_0}$. Solving this equation finishes the proof.
     
    We firstly consider $\te_{{m_0}-1}$. It is obtained via the orthogonal condition ~\eqref{eq:ortho-1}. Actually, ~\eqref{eq:ortho-1} is equivalent to 
    \begin{align}
        \left(\p_{\sigma}+{m_0}-1\right)\te_{{m_0}-1}\|R_0\|_{L^2}^2&=-\La R_0N_{m_0},R_0 \Ra-\La F_{m_0}^R,R_0 \Ra,
    \end{align}
    then we take
    \begin{align}\label{eq:theta-N}
        \te_{{m_0}-1}(\sigma)=-e^{-({m_0}-1)\sigma}\int_{-\infty}^{\sigma}e^{({m_0}-1)\ti\sigma}\frac{\La R_0N_{m_0}(\ti\sigma),R_0 \Ra+\La F_{m_0}^R(\ti\sigma),R_0 \Ra}{\|R_0\|_{L^2}^2}\ud \ti\sigma
    \end{align}
    such that ~\eqref{eq:ortho-1} holds for $m={m_0}$.

    Then we consider $R_{m_0}$. Since ~\eqref{eq:Rm-decomposition} holds for $1\leq m\leq {m_0}-1$, using ~\eqref{eq:ortho-2-equiv} we obtain ~\eqref{eq:ortho-2} holds for $m={m_0}$. Therefore, using Lemma~\ref{lem:Lk}, there exists unique function, denoted by $R_{m_0}^{\perp}$, such that $\La R_{m_0}^{\perp},\ti P_{k} \Ra=\La R_{m_0}^{\perp}, i\ti P_{k} \Ra=0$ and 
    \begin{align}
        \calL^R R_{m_0}^{\perp}-R_0N_{m_0}-R_0(\p_{\sigma}+{m_0}-1)\te_{{m_0}-1} =F_{m_0}^R.
    \end{align}
    Then we take
    \begin{align}\label{eq:R-N}
        R_{m_0}=R_{m_0}^{\perp}-\frac{1}{2\|R_0\|_{L^2}^2}\sum_{j=1}^{{m_0}-1}\La R_{j },R_{{m_0}-j} \Ra\ti P_k.
    \end{align}
    In the rest of the proof, for $f\in L^2$ such that $\La f,\ti P_k \Ra=\La f,i\ti P_k \Ra=0$, we denote by $(\calL^R)^{-1}f$ the unique function $u\in H^2 $ such that $\calL^R u=f$ and $\La u,\ti P_k\Ra=\La u,i\ti P_k \Ra=0$.

    Inserting ~\eqref{eq:theta-N} and ~\eqref{eq:R-N} into the second equation of ~\eqref{eq:RNtheta-iteration}, we obtain a nonlocal equation for $N_{m_0}$
    \begin{align}\label{eq:Nn}
        \calL_{m_0}^NN_{m_0}+\calA  (N_{m_0})= \ti F_{m_0}^{N},
    \end{align}
    where  
    \begin{align}\label{eq:An-def}
        \calA (f)\coloneqq -2\Re\De_{\eta}\left(\overline{R_0} (\calL^R)^{-1}\left( \ti P_kf-\frac{\La \ti P_k f,\ti P_k \Ra}{\La\ti P_k,\ti P_k\Ra}\ti P_k \right)\right)
    \end{align}
    is an $\R$-linear nonlocal operator for real-valued function $f$, and 
    \begin{align}
        \ti F_{m_0}^N \coloneqq  F_{m_0}^N+2\Re \De_{\eta}\left(\overline{R_0} (\calL^R)^{-1}\left( F_{m_0}^R-\frac{\La   F_{m_0}^R,\ti P_k \Ra}{\La\ti P_k,\ti P_k\Ra}\ti P_k\right) \right)+2a_{m_0}\De_{\eta}(|\ti P_k|^2).
    \end{align}
    only depends on $\{R_m,N_m,\te_{m-1}\}_{0\leq m\leq {m_0}-1}$.

    By induction, $F_{m_0}^R\in L^2$, then the term $\ti F_{m_0}^N$ is well-defined and  admits a decomposition of the form
    \begin{align}
        \ti F_{m_0}^N=\sum_{j=0}^{J_{m_0}} \ti F_{{m_0},j}(\eta)\sigma^j,
    \end{align}
    where $\ti F_{{m_0},j}(\eta)$ are real-valued and only depends on $\{R_m,N_m,\te_{m-1}\}_{0\leq m\leq {m_0}-1}$. To solve ~\eqref{eq:Nn}, we recall the ansatz
    \begin{align}
         \ N_{m_0}=\sum_{j=0}^{L_{m_0}^N}N_{{m_0},j}(\eta)\sigma^j,
    \end{align}
    and denote by 
    \begin{align}
        \calL_{{m_0},0}^N= \frac{1}{9}\left(4\eta^2\p_{\eta}^2+(26-4{m_0})\eta\p_{\eta}+{m_0}^2-11{m_0} +28\right), 
    \end{align}
    then ~\eqref{eq:Nn} turns into
    \begin{align}\label{eq:Nn-1}
        (\calL_{{m_0},0}^N+\calA )N_{{m_0},j}+\frac{j+1}{9}(-4\eta\p_{\eta}+2{m_0}-11)N_{{m_0},j+1}+\frac{(j+1)(j+2)}{9}N_{{m_0},j+2}=\ti F_{{m_0},j}
    \end{align}
    for $j\geq 0$. Here we use the convention that $N_{m,j}=0$ if $j>L_m^N$ and $\ti F_{m,j} =0$ if $j>J_{m }$.  

    We will solve ~\eqref{eq:Nn-1} in a suitable function space. Before we introduce the space, we give an integral expression of $(\calL^R)^{-1}$ using Green function. For real-valued function $f\in L^2$ such that $\La f,\ti P_k \Ra=0$, We have
    \begin{align}
        ((\calL^R)^{-1}f)(\eta)=\calG f(\eta)-\frac{\La \calG f,\ti P_k \Ra}{\La \ti P_k,\ti P_k \Ra}\ti P_k(\eta),
    \end{align}
    where
    \begin{align}
        \calG f(\eta)=&\int_0^{\infty}G(\eta,r)f(r)r^2\ud r,\\
        G(\eta,r)=&-\ti P_k(\eta)\ti P_k(r)\int_1^{\min\{\eta,r\}}\frac{1}{s^2\ti P_k^2(s)}\ud s,\ \eta,r>0.
    \end{align}
    With this expression one can study the relation of the regularity and decay properties between $f$ and $(\calL^R)^{-1}f$.
    
\begin{claim}
    \label{claim:Green}
    There holds
    \begin{align}\label{eq:Green-1}
        \sup_{0<\eta\leq 1}\left(\eta^{-2}\left( |\calG f(\eta)|+|(\eta\p_{\eta}\calG f)(\eta)|\right)\right)\lesssim \sup_{0\leq \eta\leq 1}|f(\eta)|.
    \end{align}
    For any $0<\de<1$ and $s\geq 0$, there holds
    \begin{align}\label{eq:Green-2}
        \sup_{\eta\geq 1}\left(e^{\de \eta}  \left( |\calG f(\eta)|+|(\p_{\eta}\calG f)(\eta)|\right)\right)\lesssim  \sup_{\eta\geq 0} e^{\de \eta} |f(\eta)| ,
    \end{align}
    \begin{align}\label{eq:Green-3}
        \sup_{\eta\geq 1}\left(e^{  \eta}(1+\eta)^{-s-1} \left( |\calG f(\eta)|+|(\p_{\eta}\calG f)(\eta)|\right)\right)\lesssim  \sup_{\eta\geq 0} e^{ \eta}(1+\eta)^{-s}|f(\eta)| .
    \end{align}
\end{claim}

    \begin{proof}[Proof of Claim~\ref{claim:Green}]
                Using $ \La f,\ti P_k \Ra=0 $ we have
        \begin{equation}\label{eq:Green-rewrite}
            \begin{aligned}
                \calG f(\eta)=&-\ti P_k(\eta)\int_0^{\eta}\ti P_k(r)\left(\int_{\eta}^{r}\frac{1}{s^2\ti P_k^2(s)}\ud s\right)f(r)r^2\ud r\\
            =&\int_0^{\eta}\int_0^{s}\frac{\ti P_k(\eta)\ti P_k(r)f(r)r^2}{s^2\ti P_k^2(s)}\ud r\ud s
            \end{aligned}
        \end{equation}
        Then for $0< \eta\leq 1$, since $\ti P_k(0)\neq 0$, we have
        \begin{align}
            |\calG f(\eta)|\lesssim \left|\int_0^{\eta}\int_0^{s}\frac{ |f(r)|r^2}{s^2 }\ud r\ud s\right|\lesssim \eta^2\sup_{0\leq r\leq 1}|f(r)|.
        \end{align}
        Similarly, since
        \begin{align}
            (\p_{\eta}\calG f)(\eta)=\int_0^{\eta}\int_0^{s}\frac{\ti P_k'(\eta)\ti P_k(r)f(r)r^2}{s^2\ti P_k^2(s)}\ud r\ud s+ \int_0^{\eta}\frac{\ti P_k(\eta)\ti P_k(r)f(r)r^2}{\eta^2\ti P_k^2(\eta)}\ud r ,
        \end{align}
        we have
        \begin{align}
            | (\p_{\eta}\calG f)(\eta) |\lesssim  \eta \sup_{0\leq r\leq 1}|f(r)|.
        \end{align}

        For $\eta\geq 1$, using ~\eqref{eq:Pk-asymp},  the Green function can be estimated as follows: for $r>0$,
        \begin{equation}
            \label{eq:G-bdd}
            \begin{aligned}
                |G(\eta,r)|+|\p_{\eta}G(\eta,r)|\lesssim& \frac{e^{-\eta}}{\eta}\frac{e^{-r}}{r}\left|\int_{1}^{\min\{\eta,r\}}e^{2s}\ud s\right|+|\mathbbm{1}_{\eta<r}\frac{\ti P_k(r)}{\ti P_k(\eta)\eta^2} |\lesssim\frac{e^{-|\eta-r|}}{\eta r}.
            \end{aligned}
        \end{equation}
        Then we have  
        \begin{align}
            |\calG f(\eta)|+|(\p_{\eta}\calG f)(\eta) |&\lesssim \int_0^{\infty} \frac{e^{-|\eta-r|}}{\eta r}|f(r)|r^2\ud r\\
            &\lesssim \left(\sup_{\eta\geq 0} e^{\de \eta}|f(\eta)|\right)e^{-\de\eta}\int_0^{\infty}\frac{e^{-|\eta-r|+\de(\eta-r)}}{\eta r}r^2\ud r\\
            &\lesssim \left(\sup_{\eta\geq 0} e^{\de \eta}|f(\eta)|\right)e^{-\de\eta}.
        \end{align}

        Similarly for $s\geq 0$, 
        \begin{align}
            |\calG f(\eta)|+|(\p_{\eta}\calG f)(\eta) |&\lesssim \left(\sup_{\eta\geq 0} e^{ \eta}(1+\eta)^{-s}|f(\eta)|\right)e^{- \eta}(1+\eta)^{s+1}.
        \end{align}
    \end{proof}

    \bigskip

    Using Claim~\ref{claim:Green} we can obtain the decay properties of $\ti F_{{m_0},j}$. By induction hypothesis and ~\eqref{eq:Green-2}, the components in $\ti F_{{m_0},j}$  except those coming from $\frac{1}{3}\De_{\eta}N_{{m_0}-2}$ decay exponentially, and the components coming from $\frac{1}{3}\De_{\eta}N_{{m_0}-2} $ are of the order $O(\eta^{\frac{{m_0}-2-4}{2}+1-2})=O( \eta^{\frac{{m_0}-8}{2} })$. Then we introduce a function space as follows.   
    
    Let $X_{m_0}$ denote a Banach space of complex-valued radial continuous functions in $\R^3$ which have absolutely convergent power series expressions $f(\eta)\coloneqq \sum_{\ell=0}^{\infty}f_{\ell}\eta^{2\ell}$ for $0\leq \eta\leq \rho$ and  be equipped with the  norm
    \begin{align}\label{eq:Xn-def}
        \|f\|_{X_{m_0}}:=\sum_{\ell=0}^{\infty}|f_\ell|\rho^{2\ell}+\sup_{\eta\geq 0}\frac{|f(\eta)|}{\left(1+\eta\right)^{{m_0}+1}}.
    \end{align}
    Since $ \frac{{m_0}-8}{2}<{m_0}+1$, and, by induction hypothesis and ~\eqref{eq:Green-rewrite}, that $\ti F_{{m_0},j}$  admits  an absolutely convergent power series expression in $\eta^2$ in $\{\eta:0\leq\eta\leq \rho\}$,  we have that $\ti F_{{m_0},j}\in X_{m_0}$, for all $0\leq j\leq J_{m_0}$. 
    
    Define the Banach space $\ti X_{m_0}\coloneqq \{f\in X_{m_0}: \eta\p_{\eta}f\in X_{m_0}, (\eta\p_{\eta})^2 f\in X_{m_0}\}$ equipped with the norm
    \begin{align}\label{eq:tiXn-def}
        \|f\|_{\ti X_{m_0}}\coloneqq \|f\|_{X_{m_0}}+\|\eta\p_{\eta}f\|_{X_{m_0}}+\|(\eta\p_{\eta})^2 f\|_{X_{m_0}}.
    \end{align}

    \bigskip
    
    Firstly we study the operator $\calL_{{m_0},0}^N$ on $X_{m_0}$.  We define a family of operators with  a parameter  $\lam\in \{\mu\in\C:\frac{\Re\mu-4}{2}<{m_0}+1\}$ given by
    \begin{align}
        \calL_{\lam,0}^N\coloneqq   \frac{1}{9}\left(4\eta^2\p_{\eta}^2+(26-4\lam)\eta\p_{\eta}+\lam^2-11\lam +28\right), 
    \end{align}
    with the domain ${\rm{D}}(\calL_{\lam,0}^N)=\ti X_{m_0} $. Then we have $\|\calL_{\lam,0}^Ng\|_{X_{m_0}}\lesssim \|g\|_{\ti X_{m_0}}$. 
    
    \begin{claim}\label{claim:Llam}
        Denote by $\calL(\ti X_{m_0},X_{m_0})$ the space of bounded $\C$-linear operators from $\ti X_{m_0}$ to $X_{m_0}$. The map
        \begin{equation}
            \begin{aligned}
                \begin{array}{ccc}
                    T:\{\mu\in\C:\frac{\Re\mu-4}{2}<{m_0}+1\}&\to &\calL(\ti X_{m_0},X_{m_0}), \\
                    \lam &\mapsto& \calL_{\lam,0}^N.
                \end{array}
            \end{aligned}
        \end{equation}
        is holomorphic, and for each $\lam\in \{\mu\in\C:\frac{\Re\mu-4}{2}<{m_0}+1\}$, $\calL_{\lam,0}^N$ is a Fredholm operator.
        
        Moreover, for $\lam\notin \{4\ell+4,4\ell+7\}_{\ell\in \Z_{\geq 0}}$ with $\frac{\Re\lam-4}{2}<{m_0}+1$, $\calL_{\lam,0}^N$ admits a bounded inverse $(\calL_{\lam,0}^N)^{-1}$ from $X_{m_0}$ to $\ti X_{m_0}$.
    \end{claim}
    \begin{proof}[Proof of Claim ~\ref{claim:Llam}]
        By definition we have
        \begin{align}
            \frac{\ud}{\ud{\lam}}\Big|_{\lam=\lam_0}\calL_{\lam,0}^N=\frac{1}{9}\left(-4\eta\p_{\eta}+2\lam_0-11\right)\in \calL(\ti X_{m_0},X_{m_0}),
        \end{align}
        therefore $T$ is holomorphic. 
        
        To prove that $\calL_{\lam,0}^N$ is Fredholm, we only have to prove that $\dim \ker \calL_{\lam,0}^N ,\dim\rm{coker}\calL_{\lam,0}^N<\infty$,  as it follows from $\dim{\rm{coker}}\calL_{\lam,0}^N<\infty$ that ${\rm{Ran}} \calL_{\lam,0}^N$ is closed. The finite dimension of $\dim \ker \calL_{\lam,0}^N $ follows from the standard uniqueness theory of linear ODE. Before we prove $\dim{\rm{coker}}\calL_{\lam,0}^N<\infty$, we prove that $\calL_{\lam,0}^N$ admits a bounded inverse if $\lam\notin\{4\ell+4,4\ell+7\}_{\ell\in \Z_{\geq 0}}$.
    
        For $f\in X_{m_0}$ with $f(\eta)=\sum_{\ell=0}^{\infty}f_{\ell}\eta^{2\ell}$ for $0\leq \eta\leq \rho$, we define 
        \begin{align}\label{eq:g-inner}
            g(\eta)\coloneqq \frac{9}{4}\sum_{\ell=0}^{\infty}\left(2\ell-\frac{\lam-4}{2}\right)^{-1}\left(2\ell-\frac{\lam-7}{2}\right)^{-1}f_{\ell}\eta^{2\ell}
        \end{align}
        for $0\leq \eta\leq \rho$. Since $\lam\notin \{4\ell+4,4\ell+7\}_{\ell\in \Z_{\geq 0}}$, we have $g, \eta\p_{\eta} g, \eta^2\p_{\eta}^2g$ are  continuous  and  admit absolutely convergent power series expressions  in $\{\eta:0\leq\eta\leq \rho\}$. The standard existence and uniqueness theory of linear ODE yields a unique extension of $g$ in $\{\eta:0\leq\eta<\infty\}$, which we still denote by $g$. For $\eta\geq \rho$,  using method of variation of constant we obtain
        \begin{equation}\label{eq:g-extension}
            \begin{aligned}
                g(\eta)=&\frac{-(\lam-7)g(\rho)+2\rho\p_{\eta}g(\rho)}{3}\left(\frac{\eta}{\rho}\right)^{\frac{\lam-4}{2}}+\frac{(\lam-4)g(\rho)-2\rho\p_{\eta}g(\rho)}{3}\left(\frac{\eta}{\rho}\right)^{\frac{\lam-7}{2}}\\
            &+\frac{3}{2}\int_{\rho}^{\eta}\left(\eta^{\frac{\lam-4}{2}} s^{\frac{2-\lam}{2}}-\eta^{\frac{\lam-7}{2}}s^{\frac{5-\lam}{2}}\right)f(s)\ud s.
            \end{aligned}
        \end{equation}
        Since $\frac{\Re\lam-4}{2}<{m_0}+1$,  we have $|g(\eta)|+|\eta\p_{\eta}g(\eta)|+|\eta^2\p_{\eta}^2 g(\eta)|\lesssim \|f\|_{X_{m_0}}(1+\eta)^{{m_0}+1}$ for $\eta\geq \rho$, which combined with the expression of $g$ in $\{\eta:0\leq \eta\leq \rho\}$  yields
        \begin{align}
        \|g\|_{X_{m_0}}+\|\eta\p_{\eta}g\|_{X_{m_0}}+\|\eta^2\p_{\eta}^2 g\|_{X_{m_0}}\lesssim \|f\|_{X_{m_0}}.
        \end{align}

        Finally we prove that $\dim{\rm{coker}}\calL_{\lam,0}^N<\infty$ if there exists $\ell_0\in \Z_{\geq 0}$ such that $\lam=4\ell_0+4$ or $\lam=4\ell_0+7$. In this case we can similarly prove that 
        \begin{align}
            {\rm{Ran}} \calL_{\lam,0}^N=\{ f\in X_{m_0}:f_{\ell_0}=0 \}.
        \end{align}
        Since $\eta^{2\ell_0}\lesssim \eta^{\frac{\lam-4}{2}}\ll\eta^{{m_0}+1}$ as $\eta\to \infty$, write $f=f-f_{\ell_0}\eta^{2\ell_0}+f_{\ell_0}\eta^{2\ell_0}$ we obtain that $\dim{\rm{coker}}\calL_{\lam,0}^N=1$ for $\lam\in \{4\ell+4,4\ell+7\}_{\ell\in \Z_{\geq 0}}$, which finishes the proof. 
    \end{proof}

    \bigskip

    Next we study the operator $\calA $ on $X_{m_0}$. Since $\calA $ only acts on real-valued functions, we extend it as follows. For $f=f_1+if_2\in X_{m_0}$ with $f_1,f_2$ real-valued, we have $f_1,f_2\in X_{m_0}$. We extend the operator $\calA $, keeping the same notation, by
    \begin{align}
        \calA (f_1+if_2)\coloneqq \calA (f_1)+i\calA (f_2).
    \end{align}
    Then $\calA $ is a $\C$-linear operator and maps real-valued functions to real-valued functions.
    
    \begin{claim}\label{claim:An}
        $\calA :X_{m_0}\to X_{m_0}$ is bounded, and is compact from $\ti X_{m_0}$ to $X_{m_0}$.
    \end{claim}
    \begin{proof}[Proof of Claim ~\ref{claim:An}]
               By definition of $\calA $ we only have to consider the real-valued function in $X_{m_0}$.
    
        We first prove that $\calA$ is a bounded linear operator from $X_{m_0}$ to $X_{m_0}$. For $f \in X_{m_0}$, using ~\eqref{eq:Pk-asymp} we have
        \begin{align}
            \left|\ti P_kf-\frac{\La \ti P_k f,\ti P_k \Ra}{\La\ti P_k,\ti P_k\Ra}\ti P_k\right|\lesssim \|f\|_{X_{m_0}}e^{-\eta}(1+\eta)^{m_0}.
        \end{align}
        Therefore using ~\eqref{eq:Green-3} we have
    \begin{align}
        \left|(\calL^R)^{-1}\left(\ti P_kf-\frac{\La \ti P_k f,\ti P_k \Ra}{\La\ti P_k,\ti P_k\Ra}\ti P_k\right)\right|\lesssim  \|f\|_{X_{m_0}}e^{-\eta}(1+\eta)^{{m_0}+1} ,
    \end{align}
    \begin{align}
        \left|\p_{\eta}\left((\calL^R)^{-1}\left(\ti P_kf-\frac{\La \ti P_k f,\ti P_k \Ra}{\La\ti P_k,\ti P_k\Ra}\ti P_k\right)\right)\right|\lesssim \|f\|_{X_{m_0}}e^{-\eta}(1+\eta)^{{m_0}+1} ,
    \end{align}
    \begin{align}
        &\De_{\eta}\left(  (\calL^R)^{-1}\left(\ti P_kf-\frac{\La \ti P_k f,\ti P_k \Ra}{\La\ti P_k,\ti P_k\Ra}\ti P_k\right)\right)\\
        =&\ti P_kf-\frac{\La \ti P_k f,\ti P_k \Ra}{\La\ti P_k,\ti P_k\Ra}\ti P_k +\left(1+\ti N_k\right)(\calL^R)^{-1}\left(\ti P_kf-\frac{\La \ti P_k f,\ti P_k \Ra}{\La\ti P_k,\ti P_k\Ra}\ti P_k \right)\lesssim \|f\|_{X_{m_0}}e^{-\eta}(1+\eta)^{{m_0}+1}.
    \end{align}

    Then using ~\eqref{eq:Pk-asymp} and ~\eqref{eq:initial} we have
    \begin{align}\label{eq:An-pointwise}
        |\calA (f)(\eta)|\lesssim \|f\|_{X_{m_0}}e^{-\eta}(1+\eta)^{{m_0}+1}.
    \end{align}

    Next we calculate the power series expansion of $\calA (f)$ using ~\eqref{eq:Green-rewrite}. Denote by $ g = \ti P_kf-\frac{\La \ti P_k f,\ti P_k \Ra}{\La\ti P_k,\ti P_k\Ra}\ti P_k$,  we have
    \begin{equation}
        \begin{aligned}
            \calA (f)=&-2\Re\De_{\eta}\left( \ti P_k\left( \calG g-\frac{\La \calG g,\ti P_k \Ra}{\La \ti P_k,\ti P_k\Ra}\ti P_k  \right)\right)\\
            =&-2\Re\left(\De_{\eta} \ti P_k \right)\left( \calG g \right)-4\Re \left(\p_{\eta}\ti P_k\right)\left(\p_{\eta}\calG g\right)\\
            &-2\Re \ti P_k g-2\Re\ti P_k(1+\ti N_k)\left( \calG g-\frac{\La \calG g,\ti P_k \Ra}{\La \ti P_k,\ti P_k\Ra}\ti P_k  \right)\\
            &+2\Re\left(\De_{\eta} \ti P_k\right)\left( \frac{\La \calG g,\ti P_k \Ra}{\La \ti P_k,\ti P_k\Ra}\ti P_k  \right)+4\Re\left(\p_{\eta}\ti P_k\right)\left(\frac{\La \calG g,\ti P_k \Ra}{\La \ti P_k,\ti P_k\Ra}\p_{\eta}\ti P_k \right).
        \end{aligned}
    \end{equation}
    Therefore
    \begin{align}
        \|\calA (f)\|_{X_{m_0}}\lesssim \|\calG g\|_{X_{m_0}}+\|\eta^{-1}\p_{\eta}\calG g\|_{X_{m_0}}+\|g\|_{X_{m_0}}.
    \end{align}
    By definition we have $\|g\|_{X_{m_0}}\lesssim \|f\|_{X_{m_0}}$, then it suffices to prove
    \begin{align}
        \|\calG g\|_{X_{m_0}}+\|\eta^{-1}\p_{\eta}\calG g\|_{X_{m_0}}\lesssim \|g\|_{X_{m_0}}.
    \end{align}
    We have proved that $|\calG g(\eta)|\lesssim \|g\|_{X_{m_0}} (1+\eta)^{{m_0}+1}$. To study the power series expansion of $\calG g$ in $\{\eta:0\leq\eta\leq \rho\}$, we may take $\rho$ smaller such that  $\ti P_k^{-2}$ admits an Taylor series in $\eta^2$ that converge to itself and absolutely converges in a strictly larger neighborhood of $\{\eta:0\leq\eta\leq \rho\} $. Write
    \begin{align}
        \frac{1}{\ti P_k^2(\eta)}=\sum_{j=0}^{\infty}\ti p_j\eta^{2j},\ \ti P_k(\eta)=\sum_{j=0}^{\infty}p_j\eta^{2j},\ g(\eta)=\sum_{j=0}^{\infty}g_j\eta^{2j},
    \end{align}
    then we have
    \begin{align}
         \calG g(\eta)=&\int_0^{\eta}\int_0^{s}\frac{\ti P_k(\eta)\ti P_k(r)g(r)r^2}{s^2\ti P_k^2(s)}\ud r\ud s\\
         =&\ti P_k(\eta)\int_0^{\eta}\int_0^{s}\frac{1}{s^2}\left(\sum_{j=0}^{\infty} \ti p_{j }s^{2j}\right) \left(\sum_{j=0}^{\infty} p_{j }r^{2j}\right)\left(\sum_{j=0}^{\infty}g_jr^{2j}\right)r^2\ud r\ud s\\
         =&\ti P_k(\eta)\sum_{j=0}^{\infty}\left(\sum_{\ell_1=0}^{j}\sum_{\ell_2=0}^{j-\ell_1}g_{\ell_1}\ti p_{\ell_2}p_{j-\ell_1-\ell_2}\frac{1}{(2(j-\ell_2)+3)(2j+2)}\right)\eta^{2j+2}\\
         =&\sum_{j=0}^{\infty}\left(\sum_{\ell_1=0}^j\sum_{\ell_2=0}^{j-\ell_1}\sum_{\ell_3=0}^{j-\ell_1-\ell_2}g_{\ell_1}\ti p_{\ell_2}p_{\ell_3}p_{j-\ell_1-\ell_2-\ell_3}\frac{1}{(2(\ell_1+\ell_3)+3)(2(\ell_1+\ell_2+\ell_3)+2)}\right)\eta^{2j+2}.
    \end{align}   
    Denote by
    \begin{align}
        \ti g_{j}=\left(\sum_{\ell_1=0}^j\sum_{\ell_2=0}^{j-\ell_1}\sum_{\ell_3=0}^{j-\ell_1-\ell_2}g_{\ell_1}\ti p_{\ell_2}p_{\ell_3}p_{j-\ell_1-\ell_2-\ell_3}\frac{1}{(2(\ell_1+\ell_3)+3)(2(\ell_1+\ell_2+\ell_3)+2)}\right),
    \end{align}
    then 
    \begin{align}
        \sum_{j=0}^{\infty}|\ti g_{j}|\rho^{2j+2}\lesssim & \rho^2\left(\sum_{j=0}^{\infty}|g_j|\rho^{2j}\right)\left(\sum_{j=0}^{\infty}|\ti p_j|\rho^{2j}\right)\left(\sum_{j=0}^{\infty}| p_j|\rho^{2j}\right)^2\lesssim \|g\|_{X_{m_0}},
    \end{align}
    which implies
    \begin{align}
        \|\calG g\|_{X_{m_0}}\lesssim \|g\|_{X_{m_0}}.
    \end{align}
    Similarly we have $\|\eta^{-1}\p_{\eta}\calG g\|_{X_{m_0}}\lesssim \|g\|_{X_{m_0}}$, therefore $\calA :X_{m_0}\to X_{m_0}$ is bounded.
    
    In the end, we prove that if $\{f_j\}_{j\geq 1}\subset \ti X_{m_0}$ are uniformly bounded, then $\{\calA (f_j)\}_{j\geq 1}\subset X_{m_0}$  has a convergent subsequence in $X_{m_0}$. Denote by $h_j=\calA(f_j)$ for convenience. The boundedness of $\calA $ implies $\{h_j\}$ are uniformly bounded in $X_{m_0}$. Similarly we have
    \begin{align}
        \|\eta\p_{\eta}h_j\|_{X_{m_0}}\lesssim \|f_j\|_{X_{m_0}}+\|\eta\p_{\eta}f_j\|_{X_{m_0}},
    \end{align}
    which implies $\{\eta\p_{\eta }h_j\}$ are also uniformly bounded in $X_{m_0}$. Then we write $h_j(\eta)=\sum_{\ell=0}^{\infty}h_{j,\ell}\eta^{2\ell}$ for $0\leq \eta\leq \rho$, which combined with $\eta\p_{\eta}h_j\in X_{m_0}$ yields $ \p_{\eta}h_j=\sum_{\ell=1}^{\infty}2\ell h_{j,\ell}\eta^{2\ell-1}$, for $0\leq \eta\leq \rho$. Therefore for $0\leq \eta\leq  \rho$, we have
    \begin{align}
        \sum_{\ell=1}^{\infty}2\ell|h_{j,\ell}|\eta^{2\ell-1}= \sum_{\ell=1}^{\infty}2\ell\left(\frac{\eta}{\rho}\right)^{2\ell-1}|h_{j,\ell}|\rho^{2\ell-1}\lesssim  2\rho^{-1}\|\eta\p_{\eta}h_j\|_{X_{m_0}},
    \end{align}
    which  combined with Arzel\`a-Ascoli theorem yields there exists a subsequence of $\{h_j\}$, which we still denote by $\{h_{j}\}$, converge uniformly in $\{\eta:0\leq\eta\leq  \rho\}$. Since $ \eta\p_{\eta}h_j=\sum_{\ell=1}^{\infty}2\ell h_{j,\ell}\eta^{2\ell}$ are uniformly bounded in $X_{m_0}$, for any $\varepsilon>0$ there exists $\ell_0>0$ such that $\sum_{\ell=\ell_0}^{\infty}|h_{j,\ell}|\rho^{2\ell}\leq \frac{\varepsilon}{2}$ for all $j\geq 1$. It follows by induction that $\{h_{j,\ell}\}_{j\geq 1}$ converges for each $\ell \geq 0$. Therefore
    \begin{align}
        \lim_{J\to \infty}\sup_{j_1,j_2\geq J}\left( \sum_{\ell=0}^{\infty}|h_{j_1,\ell}-h_{j_2,\ell}|\rho^{2\ell} \right)=0.
    \end{align}
    Using ~\eqref{eq:An-pointwise} we have 
    \begin{align}
        \lim_{\eta_0\to\infty}\sup_{\eta\geq \eta_0,j\geq 1}\left(\frac{|h_j(\eta)|}{(1+\eta)^{{m_0}+1}}\right)=0.
    \end{align}
    Then for any $\varepsilon>0$, there exists $\eta_0>0$ such that $\sup_{\eta\geq \eta_0,j\geq 1}\left(\frac{|h_j(\eta)|}{(1+\eta)^{{m_0}+1}}\right)<\frac{1}{3}\varepsilon$.  Since $\{\eta\p_{\eta}h_j\}$ are uniformly bounded in $X_{m_0}$, we can apply Arzel\`a-Ascoli theorem to obtain a subsequence converging in $\{\eta:\rho\leq \eta\leq \eta_0\}$. Taking $\varepsilon=\frac{1}{\ell}$ for $\ell\in \Z_{\geq 1}$, a standard diagonal argument yields the existence of the subsequence of $\{ h_j \}$ converging in $X_{m_0}$.
    \end{proof}

    \bigskip

    Combining Claim ~\ref{claim:Llam} and Claim ~\ref{claim:An} we obtain that $\{\calL_{\lam,0}^N+\calA \}_{\Re \lam<2{m_0}+6}$ is a family of Fredholm operators from $\ti X_{m_0}$ to $X_{m_0}$ and the map $\ti T:\lam\mapsto\calL_{\lam,0}^N +\calA $ is holomorphic for $\Re\lam<2{m_0}+6$. For $\lam\notin \{4\ell+4,4\ell+7\}_{\ell\in\Z_{\geq 0}}$, it follows from the proof of Claim~\ref{claim:Llam} that for real $\lam$ satisfying $ \lam\ll -1$, $f\in X_{m_0}$ and $g=((\calL_{\lam,0}^N)^{-1}f)\in \ti X_{m_0}$ with power series expansions $f(\eta)=\sum_{\ell=0}^{\infty}f_{\ell}\eta^{2\ell}$, $g(\eta)=\sum_{\ell=0}^{\infty}g_{\ell}\eta^{2\ell}$, we have
    \begin{align}
        \sum_{\ell=0}^{\infty}|g_{\ell}|\rho^{2\ell}+\sup_{0\leq\eta\leq \rho}\frac{|g(\eta)|}{(1+\eta)^{{m_0}+1}}\lesssim \sum_{\ell=0}^{\infty}|g_{\ell}|\rho^{2\ell} \lesssim \frac{1}{1+|\lam|^2} \sum_{\ell=0}^{\infty}|f_{\ell}|\rho^{2\ell}\lesssim \frac{1}{1+|\lam|^2}\|f\|_{X_{m_0}}.
    \end{align}
    For $ \eta\geq \rho$, using ~\eqref{eq:g-extension}, there holds 
    \begin{align}
        \frac{|g(\eta)|}{(1+\eta)^{{m_0}+1}}\lesssim &(1+|\lam|)|g(\rho)|+|\rho\p_{\eta}g(\rho)|+\frac{1}{1+|\lam|^2}\|f\|_{X_{m_0}}\\
        \lesssim & \frac{1+|\lam|}{1+|\lam|^2}\sum_{\ell=0}^{\infty}|f_{\ell}|\rho^{2\ell}+\frac{1}{1+|\lam|}\sum_{\ell=0}^{\infty}|f_{\ell}|\rho^{2\ell}+\frac{1}{1+|\lam|^2}\|f\|_{X_{m_0}}\\
        \lesssim & \frac{1}{1+|\lam|}\|f\|_{X_{m_0}}.
    \end{align}
    Therefore we have $ \|(\calL_{\lam,0}^N)^{-1}\|_{X_{m_0}\to  X_{m_0}}\lesssim (1+|\lam|)^{-1}\to 0$  as $ \lam\to -\infty$. Hence there exists a real $\lam\ll-1$  such that $\|\calA (\calL_{\lam,0}^N)^{-1}\|_{X_{m_0}\to  X_{m_0}}<\frac{1}{2}$, which yields $\calL_{\lam,0}^N+  \calA $ has a bounded inverse from $X_{m_0}$ to $\ti X_{m_0}$ at this $\lam$. Thus the standard analytic Fredholm theory yields that there exists a  discrete set $\Sigma\subset\{\lam:\Re\lam<2{m_0}+6\}$ without accumulation points in $\{\lam:\Re\lam<2{m_0}+6\}$ such that $(\calL_{\lam,0}^N+\calA )^{-1}$ is holomorphic in $\{\lam:\Re\lam<2{m_0}+6, \lam\notin \Sigma\}$, and that every $\lam\in\Sigma$ with $\Re\lam<2{m_0}+6$  is a pole of finite order.  For ${m_0}\geq 2$, let $q_{m_0}\geq 0$ be the order of the pole of $(\calL_{\lam,0}^N+\calA )^{-1}$ at ${m_0}$. For $\lam$ near ${m_0}$, we have
    \begin{align}\label{eq:Llam+An-inverse}
        (\calL_{\lam,0}^N+\calA )^{-1}=\sum_{\ell=0}^{\infty}(\lam-{m_0})^{\ell}\ti\calH_{{m_0},\ell}+\sum_{\ell=1}^{q_{m_0}}\frac{1}{(\lam-{m_0})^{\ell}}\calH_{{m_0},\ell},    
    \end{align}
    where  $\ti \calH_{{m_0},\ell}$ and $\calH_{{m_0},\ell}$ lie in $\calL(X_{m_0},\ti X_{m_0})$.  The equation ~\eqref{eq:Nn-1} can be rewritten into
    \begin{align}\label{eq:Nnj-rewrite}
        &\left(\calL_{\lam,0}^N+\calA \right)\Bigg|_{\lam={m_0}}(j!N_{{m_0},j})+ \frac{\ud}{\ud \lam}\Bigg|_{\lam={m_0}}\left(\calL_{\lam,0}^N+\calA \right)((j+1)!N_{{m_0},j+1})\\
        &\quad + \frac{1}{2}\frac{\ud^2}{\ud \lam^2}\Bigg|_{\lam={m_0}}\left(\calL_{\lam,0}^N+\calA \right)((j+2)!N_{{m_0},j+2})
        =j!\ti F_{{m_0},j}.
    \end{align}
    Then we take $L_{m_0}^N=J_{m_0}+q_{m_0}$, and  for $0\leq j\leq J_{m_0}+q_{m_0}$,  we take
    \begin{align}\label{eq:Nnj}
        N_{{m_0},j}=\frac{1}{j!}\left(\sum_{\ell=0}^{ J_{m_0}-j }(j+\ell)!\ti \calH_{{m_0},\ell} \ti F_{{m_0},j+\ell} +\sum_{\ell=1}^{\min\{q_{m_0},j\}}(j-\ell)!\calH_{{m_0},\ell} \ti F_{{m_0},j-\ell} \right)\in \ti X_{m_0}. 
    \end{align}
    Since $(\calL_{\lam,0}^N+\calA)(\calL_{\lam,0}^N+\calA)^{-1}=\Id$ for $\lam\neq m_0$ in a small neighborhood of $m_0$, we have
    \begin{align}\label{eq:LL-1}
        &\left((\calL_{m_0,0}^N+\calA)+(\lam-m_0)\frac{\ud}{\ud \lam}\Bigg|_{\lam=m_0}(\calL_{\lam,0}^N+\calA)+\frac{1}{2}(\lam-m_0)^2\frac{\ud^2}{\ud \lam^2}\Bigg|_{\lam=m_0}(\calL_{\lam,0}^N+\calA)\right)\\
        &\times \left(\sum_{\ell=0}^{\infty}(\lam-{m_0})^{\ell}\ti\calH_{{m_0},\ell}+\sum_{\ell=1}^{q_{m_0}}\frac{1}{(\lam-{m_0})^{\ell}}\calH_{{m_0},\ell}\right)=\Id.
    \end{align}
    Comparing the coefficients of $(\lam-m_0)^{\ell}$ in ~\eqref{eq:LL-1} yields that ~\eqref{eq:Nnj} solves ~\eqref{eq:Nn-1}. It follows from definition that $N_{{m_0},j}$ are real-valued.

    \bigskip

    Now we prove that the $N_{{m_0},j}$ we obtained in ~\eqref{eq:Nnj} and the related $R_{{m_0},j}$ obtained via ~\eqref{eq:R-N} satisfy ~\eqref{eq:RN-decay}. It suffices to study the behavior of $\ti \calH_{{m_0},\ell}$ and $\calH_{{m_0},\ell}$.
    \begin{claim}\label{claim:Hn}
        For any $0<\de<1$, $0\leq j\leq J_{m_0}$, $\eta\geq 0$ and  integers $s\geq 0$, we have 
        \begin{equation}\label{eq:Hn}
            \begin{aligned}
                 \left|\p_{\eta}^{s}\ti \calH_{{m_0},\ell}\ti F_{{m_0},j}(\eta)\right|&\lesssim_{\de,{m_0},j,\ell}  (1+\eta)^{\frac{{m_0}-4}{2}-s+\de},\ \ell\geq 0,\\
            \left|\p_{\eta}^{s} \calH_{{m_0},\ell}\ti F_{{m_0},j}(\eta)\right|&\lesssim_{\de,{m_0},j,\ell}  (1+\eta)^{\frac{{m_0}-4}{2}-s+\de},\ 1\leq \ell\leq q_{m_0}.
            \end{aligned}
        \end{equation}
    \end{claim}
    \begin{proof}[Proof of Claim ~\ref{claim:Hn}]
        It suffices to prove ~\eqref{eq:Hn} for $\eta\geq 1$. For $\ti F_{{m_0},j}$, by induction hypothesis, we obtain that 
        \begin{align}
            |\p_{\eta}^s\ti F_{{m_0},j}|\lesssim_{{m_0},j,\de} (1+\eta)^{\frac{{m_0}-2-4}{2}-2-s+\de}.
        \end{align}

        Since the poles of $(\calL_{\lam,0}^N+\calA )^{-1}$ are discrete and have no accumulation points in $\{\lam:\Re\lam<2{m_0}+6\}$, we can take $0<r<\frac{\de}{2}$ such that $\calL_{\lam,0}^N+\calA $ has a bounded inverse from $X_{m_0}$ to $\ti X_{m_0}$  for $|\lam-{m_0}|=r$. Then we have
        \begin{align}
            \ti \calH_{{m_0},\ell}\ti F_{{m_0},j}&=\frac{1}{2\pi i}\int_{|\lam-{m_0}|=r}\frac{(\calL_{\lam,0}^N+\calA )^{-1}\ti F_{{m_0},j}}{(\lam-{m_0})^{\ell+1}}\ud \lam,\ \ell\geq 0,\\
             \calH_{{m_0},\ell}\ti F_{{m_0},j}&=\frac{1}{2\pi i}\int_{|\lam-{m_0}|=r} (\lam-{m_0})^{\ell-1}(\calL_{\lam,0}^N+\calA )^{-1}\ti F_{{m_0},j} \ud \lam,\ 1\leq \ell\leq q_{m_0}.\\
        \end{align}
        Let $u_{\lam,{m_0},j}=(\calL_{\lam,0}^N+\calA)^{-1}\ti F_{{m_0},j}\in \ti X_{m_0}$, since for $|\lam-{m_0}|=r$ we have $\lam\notin \{4\ell+4,4\ell+7\}_{\ell\in \Z_{\geq 0}}$, we have
        \begin{align}
            u_{\lam,{m_0},j}=(\calL_{\lam,0}^{N})^{-1}\left(\ti F_{{m_0},j}-\calA  u_{\lam,{m_0},j} \right).
        \end{align}
        Recall the definition of $\calA $ in ~\eqref{eq:An-def} and the proof of Claim ~\ref{claim:Green} and Claim ~\ref{claim:An}, it follows similarly that
        \begin{align}
            |\calA  u_{\lam,{m_0},j}(\eta)|+|\p_{\eta}\left(\calA  u_{\lam,{m_0},j}\right)(\eta)|\lesssim e^{-\eta}(1+\eta)^{{m_0}+1}.
        \end{align}
        Therefore using the expression ~\eqref{eq:g-inner}, ~\eqref{eq:g-extension} and that $ \frac{{m_0}-8}{2}+\de<\frac{\Re\lam-4}{2}$, we have for $\eta\geq 1$,
        \begin{align}
            |u_{ \lam,{m_0},j} (\eta)|+|\eta\p_{\eta}u_{\lam,{m_0},j}(\eta)|\lesssim (1+\eta)^{\frac{\Re\lam-4}{2}}.
        \end{align}
        Using 
        \begin{align}
            \frac{4}{9}\eta^2\p_{\eta}^2u_{\lam,{m_0},j}=\ti F_{{m_0},j}-\calA  u_{\lam,{m_0},j}-\frac{1}{9}\left((26-4\lam)\eta\p_{\eta}+\lam^2-11\lam+28\right)u_{\lam,{m_0},j}
        \end{align}
        and that the derivatives of $\calA  u_{\lam,{m_0},j}$ can be estimated recursively in the same way, we have  for $\eta\geq 1$ and any integer $s\geq 0$,
        \begin{align}
            |\p_{\eta}^s u_{\lam,{m_0},j}(\eta)|\lesssim (1+\eta)^{\frac{\Re\lam-4}{2}-s} .
        \end{align}
        Note that the estimates above hold uniformly in $\lam$, therefore for $\eta\geq 1$  and integers $s\geq 0$,
        \begin{align}
             \left|\p_{\eta}^s\left(\ti \calH_{{m_0},\ell}\ti F_{{m_0},j}\right)(\eta)\right|&\lesssim\frac{1}{2\pi }\int_{|\lam-{m_0}|=r}\frac{\left|\p_{\eta}^s u_{\lam,{m_0},j}(\eta)\right|}{|\lam-{m_0}|^{\ell+1}}|\ud \lam|\lesssim_{\de}  (1+\eta)^{\frac{{m_0}-4}{2}-s+\de},\ \ell\geq 0\\
             \left|\p_{\eta}^s\left(  \calH_{{m_0},\ell}\ti F_{{m_0},j}\right)(\eta)\right|&\lesssim\frac{1}{2\pi }\int_{|\lam-{m_0}|=r}|\lam-{m_0}|^{\ell-1}\left|\p_{\eta}^s u_{\lam,{m_0},j}(\eta)\right|   |\ud  \lam|\lesssim_{\de}  (1+\eta)^{\frac{{m_0}-4}{2}-s+\de},\ 1\leq\ell\leq q_{m_0}.
        \end{align}
    \end{proof}

    The estimates in ~\eqref{eq:RN-decay} for $N_{{m_0},j}$ follows directly from ~\eqref{eq:Hn}. For $R_{{m_0},j}$, the estimates for $s=0,1$ follows from  $R_0=\ti P_k$, ~\eqref{eq:R-N} and Claim ~\ref{claim:Green}. The estimates for general $s$ follows from the identity
    \begin{align}
        \p_{\eta}^2 R_{m_0}=&R_0N_{m_0}-\frac{\La R_0N_{m_0},R_0 \Ra}{\|R_0\|_{L^2}^2}R_0+F_{m_0}^R-\frac{\La F_{m_0}^R,R_0 \Ra}{\|R_0\|_{L^2}^2}R_0\\
        &+\left(-\frac{2}{\eta}\p_{\eta}+1+N_0\right)R_{m_0}.
    \end{align}
    It follows directly from the construction that $N_{m_0}$ and $\te_{{m_0}-1}$ are real-valued.  

    It remains to prove that $R_{{m_0},j}$ and $N_{{m_0},j}$ satisfies ~\eqref{eq:assumption-strengthen} with $m$ replaced by ${m_0}$, which follows directly from $N_{{m_0},j}\in \ti X_{m_0}$ and ~\eqref{eq:R-N}. This finishes the proof.
     
     \end{proof}

    \bigskip

    In the end of this section, we prove Proposition ~\ref{prop:core} using Lemma~\ref{lem:iteration-solve}.

    \begin{proof}[Proof of Proposition ~\ref{prop:core}]
        Fix $s\in\Z_{\geq 0}$, take an integer $M\gg 1$ to be determined later. Take
        \begin{equation}\label{eq:core-app}
            \begin{aligned}
                 \psi^{\rm{co}}_{M}(t,r)=&\frac{1}{\tau^3}e^{i\tau^{-1}-i\frac{1}{2}\tau\eta^2-i\sum_{m=0}^{M-1}\tau^m\te_m(\log\tau)}\sum_{m=0}^M \tau^m R_m(\log\tau,\eta),\\
                 n^{\rm{co}}_{M}(t,r)=&\frac{1}{3\tau^4}\sum_{m=0}^{M}\tau^mN_m(\log\tau.\eta),
            \end{aligned}
        \end{equation}
         where  $\tau=(T-t)^{\frac{1}{3}}$, $\eta=\frac{r}{\sqrt{3}(T-t)^{\frac{2}{3}}}$, and $\{R_m,N_m,\te_{m}\}_{m\in\Z_{\geq 0}}$ are constructed in ~\eqref{eq:initial} and  Lemma ~\ref{lem:iteration-solve}. Inserting ~\eqref{eq:core-app} into ~\eqref{eq:E-core}, and recall that $\sigma=\log\tau$, we obtain
         \begin{equation}
             \begin{aligned}
                 \calE^{\rm{co}}_{\psi,M}=& e^{i\tau^{-1}-i\frac{1}{2}\tau\eta^2-i\sum_{m=0}^{M-1}\tau^m\te_m(\log\tau)}\\
                 &\times\Bigg[-\frac{i}{3}\tau^{M-6}(\p_{\sigma}+M)R_M+\frac{\tau^{M-6}}{6}\eta^2(R_{M-1}+\tau R_M)\\
                 &\quad \quad -\frac{1}{3}\sum_{\substack{ 0\leq m,\ell\leq M\\ m+\ell\geq M+1}} \tau^{m+\ell-7}R_mN_{\ell}-\frac{1}{3}\sum_{\substack{0\leq m\leq M\\ 0\leq \ell\leq M-1\\ m+\ell\geq M }}\tau^{m+\ell-6}R_m(\p_{\sigma}+\ell)\te_{\ell}\Bigg],\\
                 \calE^{\rm{co}}_{n,M}=&-\frac{\tau^{M-9}}{9}\De_{\eta}N_{M-1}-\frac{\tau^{M-8}}{9}\De_{\eta}N_M
                 -\frac{1}{3}\sum_{\substack{0\leq m,\ell\leq M\\ m+\ell\geq M+1}}\tau^{m+\ell-10}\De_{\eta}\left(\overline{R_m}R_{\ell}\right).
             \end{aligned}
         \end{equation}
         In the region $r\leq 2\sqrt{3}(T-t)^{\frac{1}{2}}$, we have $0\leq\eta\leq 2\tau^{-\frac{1}{2}}$.  Using ~\eqref{eq:RN-decay} with $\de=\frac{3}{4}$ for $R_{m,j}$ and $\de=\frac{1}{2}$ for $N_{m,j}$ and the fact that
        $\p_t=-\frac{1}{3\tau^3}(\tau\p_{\tau}-2\eta\p_{\eta})$, $\p_{r}=\frac{1}{\sqrt{3}\tau^2}\p_{\eta}$, for any $a_1, a_2\in \Z_{\geq 0}$ with $a_1+a_2\leq s$,  we have
        \begin{equation}
            \begin{aligned}
                |\p_t^{a_1}\p_r^{a_2}\calE^{\rm{co}}_{\psi,M}(t,r)|\lesssim_{M,s} & \tau^{M-6-4s}(1+|\log\tau|)^{J_{M,s}} e^{-\frac{\eta}{2}},\\
                |\p_t^{a_1}\p_r^{a_2}\calE^{\rm{co}}_{n,M}(t,r)|\lesssim_{M,s} & \tau^{M-9-3s}(1+|\log\tau|)^{J_{M,s}} (1+\eta)^{\frac{M-8}{2} },
            \end{aligned}
        \end{equation}
         where $J_{M,s}\in \Z_{\geq 0}$. Therefore
         \begin{equation}
             \begin{aligned}
                 \sup_{\substack{ a_1, a_2\in \Z_{\geq 0} \\ a_1+a_2\leq s}}\sup_{r\leq2\sqrt{3} (T-t)^{1/2}}\left(\left|\p_{t}^{a_1}\p_r^{a_2}\calE^{\rm{co}}_{\psi,M}(t,r)\right|+ \left|\p_{t}^{a_1}\p_r^{a_2}\calE^{\rm{co}}_{n,M}(t,r)\right|\right)\\
                 \lesssim_{M,s}  \left(\tau^{M-6-4s}+\tau^{\frac{3}{4}M-7-3s}\right)(1+|\log\tau|)^{J_{M,s}}.
             \end{aligned}
         \end{equation}
        Then for $M\geq 12+12s$, the related $\calE^{\rm{co}}_{\psi,M},\calE^{\rm{co}}_{n,M}$ satisfies ~\eqref{eq:E-core-decay}. It remains to prove ~\eqref{eq:core-converge}. Using ~\eqref{eq:RN-decay}, ~\eqref{eq:initial} and ~\eqref{eq:core-app} we obtain that for any $L>0$ and $s'\in \Z_{\geq 0}$, there exists a integer $\ti J_{s' }$ such that
        \begin{align}
            &\|\tau^3e^{-i\tau^{-1}}\psi^{\rm{co}}_{M}(t,\sqrt{3}\tau^2\eta)-\ti P_k(\eta) \|_{C^{s' }([0,L(T-t)^{-\frac{1}{6}}])}+\|3\tau^4n^{\rm{co}}_{M}(t,\sqrt{3}\tau^2\eta)-\ti N_k(\eta)\|_{C^{s' }([0,L(T-t)^{-\frac{1}{6}}])}\\
            \lesssim &\ \tau(1+|\log\tau|)^{\ti J_{s' }}\to 0,\mas \  t\to T,
        \end{align}
        which finishes the proof
         
    \end{proof}

\section{Construction of the approximate solution}\label{sec:remote}

In this section we construct the approximate solution in the whole region and estimate the error.

\subsection{Correction in the remote region}

We first present a lemma that describes the behavior of the approximate solution $(\psi^{\rm{co}}_{M},n^{\rm{co}}_{M})$ in the region $\{\eta\geq 1\}$, recalling that $\tau=(T-t)^{\frac{1}{3}}$, $\eta=\frac{r}{\sqrt{3}(T-t)^{\frac{2}{3}}}$, $\sigma=\log\tau$.

\begin{lem}
    \label{lem:N-infity}
    Let $\{R_m,N_m,\te_{m}\}_{m\geq 0}$  be given by ~\eqref{eq:initial} and Lemma ~\ref{lem:iteration-solve}. Then for every $m\geq 0$ and $\eta\geq 1$, there holds
    \begin{equation}\label{eq:Nm-expansion-infty}
        \begin{aligned}
            N_{m}(\sigma,\eta)=\sum_{j=0}^{[\frac{m}{2}]}\eta^{\frac{m-4-6j}{2}}\ti N_{m,1,j}\left(\sigma+\frac{1}{2}\log\eta\right)+\sum_{j=0}^{[\frac{m}{2}]}\eta^{\frac{m-7-6j}{2}}\ti N_{m,2,j}\left(\sigma+\frac{1}{2}\log\eta\right)+\ti N_{m,3}(\sigma,\eta),
        \end{aligned}
    \end{equation}
    where $\ti N_{m,1,j}$ and $\ti N_{m,2,j}$ are polynomials and $\ti N_{m,3}$ satisfies
    \begin{align}\label{eq:tiNm3-decay}
        |\p_{\sigma}^{a_1}\p_{\eta}^{a_2} \ti N_{m,3}(\sigma,\eta)|\lesssim_{m,\de,a_1,a_2} (1+|\sigma|)^{\hat J_m}e^{-\de\eta}
    \end{align}
    for all $a_1,a_2\in \Z_{\geq 0}$, $0<\de<2$ and some $\hat J_m\in \Z_{\geq 0}$.
    Moreover, for $\ti N_{m,1,j}$ and $\ti N_{m,2,j}$ with $1\leq j\leq [\frac{m}{2}]$ we have the following relations:
    \begin{equation}\label{eq:tiN12-recursion}
        \begin{aligned}
             \ti N_{m ,1,j}(\zeta)&=\frac{1}{6j(2j-1)}\left(\frac{1}{2}\p_{\zeta}+\frac{m-6j}{2}\right)\left(\frac{1}{2}\p_{\zeta}+\frac{m-6j+2}{2}\right)\ti N_{m-2,1,j-1}(\zeta),\\
        \ti N_{m,2,j}(\zeta)&=\frac{1}{6j(2j+1)}\left(\frac{1}{2}\p_{\zeta}+\frac{m-6j-3}{2}\right)\left(\frac{1}{2}\p_{\zeta}+\frac{m-6j-1}{2}\right)\ti N_{m -2,2,j-1}(\zeta)
        \end{aligned}
    \end{equation}
\end{lem}

\begin{proof}[Proof of Lemma~\ref{lem:N-infity}]
    We argue by induction. The case $m=0$ follows from ~\eqref{eq:tiNk-P} and ~\eqref{eq:PN-decay} once we write $\int_0^{\eta}=\int_0^{\infty}-\int_{\eta}^{\infty}$ in ~\eqref{eq:tiNk-P}. The case $m=1$ follows from that $R_1=N_1\equiv 0$ as shown in the proof of Lemma~\ref{lem:iteration-solve}.  Suppose that ~\eqref{eq:Nm-expansion-infty} holds for $0\leq m\leq m_0-1$ where $m_0\geq 2$ is an integer, Denote by $\hat N(\zeta,\eta)=N(\zeta-\frac{1}{2}\log\eta,\eta)$, and  $\hat R(\zeta,\eta)=R(\zeta-\frac{1}{2}\log\eta,\eta)$, then in variables $(\zeta,\eta)$, the operators $\calL_m^N$ and $\De_{\eta}$ turn into
    \begin{align}
        \hat\calL_{m}^N=\frac{4}{9}\left( \eta\p_{\eta}+\frac{7-m}{2} \right)\left( \eta\p_{\eta}+\frac{4-m}{2} \right),\ \hat\De_{\eta}=\p_{\eta}^2+\frac{2}{\eta}\p_{\eta}+\frac{1}{\eta}\p_{\zeta}\p_{\eta}+\frac{1}{2\eta^2}\p_{\zeta}+\frac{1}{4\eta^2}\p_{\zeta}^2.
    \end{align}

    Using the method of variation of constants and the standard uniqueness theory of linear ODE, we obtain the expression of $N_{m_0}$ for $\eta\geq 1$:
    \begin{equation}
        \begin{aligned}
            \hat N_{m_0}(\zeta,\eta)=&\frac{-({m_0}-7)\hat N_{m_0}(\zeta,1)+2\p_{\eta}\hat N_{m_0}(\zeta,1)}{3}\eta^{\frac{{m_0}-4}{2}}\\
            &+\frac{({m_0}-4)\hat N_{m_0}(\zeta,1)-2\p_{\eta}\hat N_{m_0}(\zeta,1)}{3}\eta^{\frac{{m_0}-7}{2}}\\
            &+\frac{3}{2}\int_1^{\eta} \left(\eta^{\frac{m_0-4}{2}}s^{\frac{2-{m_0}}{2}}-\eta^{\frac{{m_0}-7}{2}}s^{\frac{5-{m_0}}{2}}\right) \\
            &\quad\quad\quad\quad\quad\quad\times\left(\frac{1}{3}\hat \De_{s}\hat N_{{m_0}-2}+\sum_{j=0}^{m_0}\hat\De_{s}\left(\overline{\hat R_j}\hat R_{m_0-j}\right)\right)(\zeta,s)\ud s.
        \end{aligned}
    \end{equation}
    Then in the original variables $(\sigma,\eta)$ we have
    \begin{equation}\label{eq:Nm-iteration-vc}
        \begin{aligned}
            N_{m_0}(\sigma,\eta)=&\frac{(2\eta\p_{\eta}-\p_{\sigma}-m_0+7)N_{m_0}(\sigma+\frac{1}{2}\log\eta,1)}{3}\eta^{\frac{m_0-4}{2}}\\
            &+\frac{(-2\eta\p_{\eta}+\p_{\sigma}+m_0-4)N_{m_0}(\sigma+\frac{1}{2}\log\eta,1)}{3}\eta^{\frac{m_0-7}{2}}\\
            &+\frac{3}{2}\int_1^{\eta} \left(\eta^{\frac{m_0-4}{2}}s^{\frac{2-{m_0}}{2}}-\eta^{\frac{{m_0}-7}{2}}s^{\frac{5-{m_0}}{2}}\right) \\
            &\quad\quad\quad\quad\quad\quad\times\left(\frac{1}{3}  \De_{s}  N_{{m_0}-2}+\sum_{j=0}^{m_0} \De_{s}\left(\overline{  R_j}  R_{m_0-j}\right)\right)(\sigma+\frac{1}{2}\log\frac{\eta}{s},s)\ud s.
        \end{aligned}
    \end{equation}
    Inserting ~\eqref{eq:Nm-expansion-infty} into ~\eqref{eq:Nm-iteration-vc} we obtain
    \begin{equation}
        \begin{aligned}
            N_{m_0}(\sigma,\eta)=&\frac{(2\eta\p_{\eta}-\p_{\sigma}-m_0+7)N_{m_0}(\sigma+\frac{1}{2}\log\eta,1)}{3}\eta^{\frac{m_0-4}{2}}\\
            &+\frac{(-2\eta\p_{\eta}+\p_{\sigma}+m_0-4)N_{m_0}(\sigma+\frac{1}{2}\log\eta,1)}{3}\eta^{\frac{m_0-7}{2}}\\
            &+\frac{1}{2}\sum_{j=0}^{[\frac{m_0}{2}]-1}\int_1^{\eta} \left(\eta^{\frac{m_0-4}{2}}s^{\frac{2-{m_0}}{2}}-\eta^{\frac{{m_0}-7}{2}}s^{\frac{5-{m_0}}{2}}\right)  s^{\frac{m_0-6-6j}{2}-2}\\
            & \times \left[\frac{1}{4} \p_{\zeta}^2+\frac{m_0-5-6j}{2}\p_{\zeta}+\frac{m_0-6-6j}{2}\frac{m_0-4-6j}{2} \right]\ti N_{m_0-2,1,j}(\sigma+\frac{1}{2}\log\eta)    \ud s\\
            &+\frac{1}{2}\sum_{j=0}^{[\frac{m_0}{2}]-1}\int_1^{\eta} \left(\eta^{\frac{m_0-4}{2}}s^{\frac{2-{m_0}}{2}}-\eta^{\frac{{m_0}-7}{2}}s^{\frac{5-{m_0}}{2}}\right)  s^{\frac{m_0-9-6j}{2}-2}\\
            & \times \left[\frac{1}{4} \p_{\zeta}^2+\frac{m_0-8-6j}{2}\p_{\zeta}+\frac{m_0-9-6j}{2}\frac{m_0-7-6j}{2} \right]\ti N_{m_0-2,2,j}(\sigma+\frac{1}{2}\log\eta)    \ud s\\
            &+   \frac{1}{2} \left(\int_1^{\infty}-\int_{\eta}^{\infty} \right)\left(\eta^{\frac{m_0-4}{2}}s^{\frac{2-{m_0}}{2}}-\eta^{\frac{{m_0}-7}{2}}s^{\frac{5-{m_0}}{2}}\right) \left(\De_{s} \ti N_{m_0-2,3} \right)(\sigma+\frac{1}{2}\log\frac{\eta}{s},s)\ud s \\
            &+   \frac{3}{2}\sum_{j=0}^{m_0}\left(\int_1^{\infty}-\int_{\eta}^{\infty}  \right)\left(\eta^{\frac{m_0-4}{2}}s^{\frac{2-{m_0}}{2}}-\eta^{\frac{{m_0}-7}{2}}s^{\frac{5-{m_0}}{2}}\right) \left(\De_{s}\left(\overline{  R_j}  R_{m_0-j}\right)\right)(\sigma+\frac{1}{2}\log\frac{\eta}{s},s)\ud s.
        \end{aligned}
    \end{equation}
    A direct computation yields that $N_{m_0}$  satisfies ~\eqref{eq:Nm-expansion-infty} and ~\eqref{eq:tiN12-recursion}.

\end{proof}

Note that $\sigma+\frac{1}{2}\log\eta=\frac{1}{2}\log(\tau^2\eta)=\frac{1}{2}\log\frac{r}{\sqrt{3}}$ and
\begin{align}\label{eq:N-infty-tr}
    &\frac{1}{3\tau^4}\sum_{m=0}^M\tau^m\left(\sum_{j=0}^{[\frac{m}{2}]}\eta^{\frac{m-4-6j}{2}}\ti N_{m,1,j}\left(\sigma+\frac{1}{2}\log\eta\right)+\sum_{j=0}^{[\frac{m}{2}]}\eta^{\frac{m-7-6j}{2}}\ti N_{m,2,j}\left(\sigma+\frac{1}{2}\log\eta\right) \right)\\
    =&\frac{1}{3 }\sum_{m=0}^M\left(\sum_{j=0}^{[\frac{m}{2}]}\tau^{ 6j}(\tau^2\eta)^{\frac{m-4-6j}{2}}\ti N_{m,1,j}\left(\sigma+\frac{1}{2}\log\eta\right)+\sum_{j=0}^{[\frac{m}{2}]}\tau^{3+6j}(\tau^2\eta)^{\frac{m-7-6j}{2}}\ti N_{m,2,j}\left(\sigma+\frac{1}{2}\log\eta\right) \right)\\
    =&\frac{1}{3 }\sum_{m=0}^M\sum_{j=0}^{[\frac{m}{2}]}(T-t)^{ 2j}\left(\frac{r}{\sqrt{3}}\right)^{\frac{m-4-6j}{2}}\ti N_{m,1,j}\left(\frac{1}{2}\log\frac{r}{\sqrt{3}}\right)   \\
    &+\frac{1}{3 }\sum_{m=0}^M \sum_{j=0}^{[\frac{m}{2}]}(T-t)^{1+2j}\left(\frac{r}{\sqrt{3}}\right)^{\frac{m-7-6j}{2}}\ti N_{m,2,j}\left(\frac{1}{2}\log\frac{r}{\sqrt{3}}\right) .
\end{align}
 which suggests us to define
\begin{equation}\label{eq:nre-initial}
    \begin{aligned}
        n_{0,M}^{\rm{re}}(r)&=\frac{1}{3}\chi(r)\sum_{m=0}^M\left(\frac{r}{\sqrt{3}}\right)^{\frac{m-4}{2}}\ti N_{m,1,0}\left(\frac{1}{2}\log\frac{r}{\sqrt{3}}\right),\\
        n_{1,M}^{\rm{re}}(r)&=-\frac{1}{3}\chi(r)\sum_{m=0}^M\left(\frac{r}{\sqrt{3}}\right)^{\frac{m-7}{2}}\ti N_{m,2,0}\left(\frac{1}{2}\log\frac{r}{\sqrt{3}}\right),
    \end{aligned}
\end{equation} 
where $r>0$ and $\chi$ is the smooth cutoff function given in ~\eqref{eq:chi}. 

Note that $\psi_M^{\rm{co}}$ decays exponentially in ~\eqref{eq:RN-decay}, we take  $\psi_M^{\rm{re}}\equiv0$ and denote by $(n_M^{\rm{re}},\p_t n_M^{\rm{re}})$ the solution to 
\begin{align}\label{eq:nMre}
    \begin{cases}
         (\p_t^2-\De_r)n_M^{\rm{re}}=0,\\
         (n_M^{\rm{re}}(T),\p_tn_M^{\rm{re}}(T))=(n_{0,M}^{\rm{re}},n_{1,M}^{\rm{re}}),
    \end{cases}\quad (t,r)\in\{(t,r):r>T-t\}.
\end{align}
For $r>T-t$, the solution can be expressed as
\begin{equation}\label{eq:nre-express}
    \begin{aligned}
        n_M^{\rm{re}}(t,r)=&\frac{(r+T-t)n_{0,M}^{\rm{re}}(r+T-t)+(r-T+t)n_{0,M}^{\rm{re}}(r-T+t)}{2r}\\
    &-\frac{1}{2r}\int_{r-T+t}^{r+T-t}sn_{1,M}^{\rm{re}}(s)\ud s.
    \end{aligned}
\end{equation}

Next we introduce a  correction $n^{\rm{re}}(t,r)$ defined uniformly for all $M\geq 1$. An immediate result is Lemma~\ref{lem:ap-diff}, which plays an important role in obtaining a smooth blow up solutions to ~\eqref{eq:Za}. See Proposition~\ref{prop:smooth} for details.

\begin{lem}
    \label{lem:psi-n-re}
    There exists a strictly decreasing sequence $\{r_m\}_{m\in\Z_{\geq 0}}$ such that the functions
    \begin{equation}\label{eq:nre-initial-1}
    \begin{aligned}
        n_{0}^{\rm{re}}(r)&=\frac{1}{3}\chi(r)\sum_{m=0}^\infty\chi\left(\frac{r}{r_m}\right)\left(\frac{r}{\sqrt{3}}\right)^{\frac{m-4}{2}}\ti N_{m,1,0}\left(\frac{1}{2}\log\frac{r}{\sqrt{3}}\right),\\
        n_{1}^{\rm{re}}(r)&=-\frac{1}{3}\chi(r)\sum_{m=0}^\infty\chi\left(\frac{r}{r_m}\right)\left(\frac{r}{\sqrt{3}}\right)^{\frac{m-7}{2}}\ti N_{m,2,0}\left(\frac{1}{2}\log\frac{r}{\sqrt{3}}\right)
    \end{aligned}
\end{equation} 
    are smooth for $r>0$, supported in $\{r\leq 2\}$. For $r>T-t$, define
    \begin{equation}\label{eq:nre-express-1}
    \begin{aligned}
        n^{\rm{re}}(t,r)=&\frac{(r+T-t)n_{0}^{\rm{re}}(r+T-t)+(r-T+t)n_{0}^{\rm{re}}(r-T+t)}{2r}\\
    &-\frac{1}{2r}\int_{r-T+t}^{r+T-t}sn_{1}^{\rm{re}}(s)\ud s.
    \end{aligned}
\end{equation}
    Then we have the following estimates.
    \begin{enumerate}
        \item For any integers $M,a\geq 0$, there holds 
        \begin{equation}\label{eq:nre-diff}
            \begin{aligned}
                |\p_r^a (n_0^{\rm{re}}-n_{0,M}^{\rm{re}})(r)|&\lesssim_{M,a } r^{\frac{M-4}{2}-a },\\ 
                |\p_r^a (n_1^{\rm{re}}-n_{1,M}^{\rm{re}})(r)|&\lesssim_{M,a } r^{\frac{M-7}{2}-a }.
            \end{aligned}
        \end{equation}
        \item For any integers $M,a,b\geq 0$  and any $0<T-t<\frac{r}{2}$, there holds
        \begin{equation}\label{eq:nre-diff-1}
            \begin{aligned}
                |\p_t^a\p_r^b (n^{\rm{re}}-n_{M}^{\rm{re}})(t,r)|&\lesssim_{M,a,b} r^{\frac{M-6}{2}-a-b}.
            \end{aligned}
        \end{equation}
    \end{enumerate}
\end{lem}

\begin{proof}[Proof of Lemma~\ref{lem:psi-n-re}]
    Take $0<2r_1<r_0<\frac{1}{2}$. For each $m\geq 2$, choose $r_m<\frac{r_{m-1}}{2}$ sufficiently small such that for any integer $a$ satisfying $0\leq a\leq m$ and every $r>0$, there hold
    \begin{equation}\label{eq:rm-choice}
        \begin{aligned}
            \left|\p_r^a\left(\chi\left(\frac{r}{r_m}\right)\left(\frac{r}{\sqrt{3}}\right)^{\frac{m-4}{2}}\ti N_{m,1,0}\left(\frac{1}{2}\log\frac{r}{\sqrt{3}}\right)\right)\right|\leq &r^{\frac{m-5}{2}-a},\\
        \left|\p_r^a\left(\chi\left(\frac{r}{r_m}\right)\left(\frac{r}{\sqrt{3}}\right)^{\frac{m-7}{2}}\ti N_{m,2,0}\left(\frac{1}{2}\log\frac{r}{\sqrt{3}}\right)\right)\right|\leq &r^{\frac{m-8}{2}-a}.
        \end{aligned}
    \end{equation}
    Since for each $r>0$, the sums in ~\eqref{eq:nre-initial-1} are finite, $n_0^{\rm{re}},n_1^{\rm{re}}$ are smooth for $r>0$ and are supported in $\{r\leq 2\}$. Fix the integers $M\geq 0$, for $0<r\leq r_M$ we have
    \begin{align}
        n_0^{\rm{re}}(r)-n_{0,M}^{\rm{re}}(r)&=\frac{1}{3}\sum_{m=M+1}^{\infty}\chi\left(\frac{r}{r_m}\right)\left(\frac{r}{\sqrt{3}}\right)^{\frac{m-4}{2}}\ti N_{m,1,0}\left(\frac{1}{2}\log\frac{r}{\sqrt{3}}\right),\\
        n_1^{\rm{re}}(r)-n_{1,M}^{\rm{re}}(r)&=-\frac{1}{3}\sum_{m=M+1}^{\infty}\chi\left(\frac{r}{r_m}\right)\left(\frac{r}{\sqrt{3}}\right)^{\frac{m-7}{2}}\ti N_{m,2,0}\left(\frac{1}{2}\log\frac{r}{\sqrt{3}}\right),
    \end{align}
    which combined with ~\eqref{eq:rm-choice} yields for $0<r\leq  r_M$, ~\eqref{eq:nre-diff} hold. 
    For $0<2(T-t)<r\leq \frac{2}{3}r_M$, we also have
    \begin{align}
         n^{\rm{re}}(t,r)-n_M^{\rm{re}}(t,r)=&\frac{r+T-t}{6r} \sum_{m=M+1}^{\infty}\chi\left(\frac{r+T-t}{r_m}\right)\left(\frac{r+T-t}{\sqrt{3}}\right)^{\frac{m-4}{2}}\ti N_{m,1,0}\left(\frac{1}{2}\log\frac{r+T-t}{\sqrt{3}}\right) \\
         &+\frac{r-T+t}{6r} \sum_{m=M+1}^{\infty}\chi\left(\frac{r-T+t}{r_m}\right)\left(\frac{r-T+t}{\sqrt{3}}\right)^{\frac{m-4}{2}}\ti N_{m,1,0}\left(\frac{1}{2}\log\frac{r-T+t}{\sqrt{3}}\right)\\
    &+\frac{1}{6r}\int_{r-T+t}^{r+T-t}s\left( \sum_{m=M+1}^{\infty}\chi\left(\frac{s}{r_m}\right)\left(\frac{s}{\sqrt{3}}\right)^{\frac{m-7}{2}}\ti N_{m,2,0}\left(\frac{1}{2}\log\frac{s}{\sqrt{3}}\right)\right)\ud s.
    \end{align}
    which combined with ~\eqref{eq:rm-choice} yields for $2(T-t)<r\leq \frac{2}{3}r_M$, ~\eqref{eq:nre-diff-1} hold.
    
    For $r>\frac{2}{3}r_M$, ~\eqref{eq:nre-diff} and ~\eqref{eq:nre-diff-1}  hold trivially.

\end{proof}

For $t<T$ sufficiently close to $T$, we have $(T-t)^{\frac{1}{2}}\gg T-t$. We define
\begin{equation}\label{eq:ap-def}
    \begin{aligned}
        \psi^{\rm{ap}}_M(t,r)&=\chi\left(\frac{r}{\sqrt{3}(T-t)^{\frac{1}{2}}}\right)\psi_M^{\rm{co}}(t,r),\\
        n^{\rm{ap}}_M(t,r)&=\chi\left(\frac{r}{\sqrt{3}(T-t)^{\frac{1}{2}}}\right)n_M^{\rm{co}}(t,r)+\left(1-\chi\left(\frac{r}{\sqrt{3}(T-t)^{\frac{1}{2}}}\right)\right)n^{\rm{re}}(t,r).
    \end{aligned}
\end{equation}
Then for $\frac{r}{\sqrt{3}(T-t)^{\frac{1}{2}}}\geq 1 $, we have $r\gg T-t$, thus $n_M^{\rm{ap}}$ is well-defined. We denote by $\chi_1(t,r)=\chi\left(\frac{r}{\sqrt{3}(T-t)^{\frac{1}{2}}}\right)$ for convenience.  

\begin{lem}\label{lem:ap-diff}
    For any integers $s\geq 0$, $M_1>M_2\geq 4s+24$ and $t<T$ sufficiently close to $T$, there holds
    \begin{equation}
        \begin{aligned}
            &\|\psi_{M_1}^{\rm{ap}}-\psi_{M_2}^{\rm{ap}}\|_{H^{s+2}}+\|\p_t(\psi_{M_1}^{\rm{ap}}-\psi_{M_2}^{\rm{ap}})\|_{H^{s}}+\|n_{M_1}^{\rm{ap}}-n_{M_2}^{\rm{ap}}\|_{H^{s+1}}+\|\p_t(n_{M_1}^{\rm{ap}}-n_{M_2}^{\rm{ap}})\|_{H^{s}}\\
            \lesssim_{M_1} & (T-t)^{\frac{M_2}{4}-2s-4}.
        \end{aligned}
    \end{equation}

\end{lem}
\begin{proof}[Proof of Lemma~\ref{lem:ap-diff}]
    We have $\psi_{M_1}^{\rm{ap}}-\psi_{M_2}^{\rm{ap}}=\chi_1(\psi_{M_1}^{\rm{co}}-\psi_{M_2}^{\rm{co}}) $, $n_{M_1}^{\rm{ap}}-n_{M_2}^{\rm{ap}}=\chi_1(n_{M_1}^{\rm{co}}-n_{M_2}^{\rm{co}}) $, the proof follows directly from  ~\eqref{eq:core-app}, ~\eqref{eq:RN-decay} and ~\eqref{eq:Nm-expansion-infty}
\end{proof}

\subsection{Estimates on the error}\label{subsec:E}

Recalling the definition of $(\psi_M^{\rm{ap}},n_M^{\rm{ap}}) $ in ~\eqref{eq:ap-def}, we define
        \begin{equation}\label{eq:E}
            \begin{aligned}
                \calE^{\rm{ap}}_{\psi,M}=&i\p_t\psi^{\rm{ap}}_{M}+\De\psi^{\rm{ap}}_{M}-n^{\rm{ap}}_{M}\psi^{\rm{ap}}_{M},\\
            \calE^{\rm{ap}}_{n,M}=&\p_{tt}n^{\rm{ap}}_{M}-\De n^{\rm{ap}}_{M}-\De|\psi^{\rm{ap}}_{M}|^2.
            \end{aligned}
        \end{equation}

\begin{lem}\label{lem:E}
    For every integer  $s\geq0$, there exists $M_{s}\in \Z_{\geq 0}$ such that for any integer $M\geq M_{s}$ and for all $t<T$ sufficiently close to $T$, we have 
    \begin{align}\label{eq:E-supp}
        \supp \calE^{\rm{ap}}_{\psi,M},\ \supp \calE^{\rm{ap}}_{n,M}\subset \{(t,r):r\leq2\sqrt{3}(T-t)^{\frac{1}{2}}\},
    \end{align}
    and the following estimates hold for any $a_1,a_2\in \Z_{\geq 0}$ with $a_1+a_2\leq s,$:
    \begin{equation}\label{eq:E-estimate}
        \begin{aligned}
            \|\p_t^{a_1}\p_r^{a_2}\calE^{\rm{ap}}_{\psi,M}\|_{L^{\infty}}+\|\p_t^{a_1}\p_r^{a_2}\calE^{\rm{ap}}_{\psi,M}\|_{L^2}+\|\p_t^{a_1}\p_r^{a_2}\calE^{\rm{ap}}_{n,M}\|_{L^{\infty}}+\|\p_t^{a_1}\p_r^{a_2}\calE^{\rm{ap}}_{n,M}\|_{L^2}\lesssim_{M,s}(T-t)^{\frac{1}{4}M-4-2s}.
        \end{aligned}
    \end{equation}
\end{lem}

\begin{proof}[Proof of Lemma~\ref{lem:E}]
    ~\eqref{eq:E-supp} follows directly from ~\eqref{eq:ap-def}, ~\eqref{eq:E} and the finite speed of propagation of wave equation. A direct calculation combining with ~\eqref{eq:nre-express-1} yields
    \begin{equation}\label{eq:E-expression}
        \begin{aligned}
            \calE_{\psi,M}^{\rm{ap}}=&\chi_1\calE_{\psi,M}^{\rm{co}}+(i\p_t\chi_1+\De_r\chi_1)\psi_{M}^{\rm{co}}+2\p_r \chi_1\cdot\p_r \psi_M^{\rm{co}}+\chi_1(1-\chi_1)\psi_M^{\rm{co}}(n_M^{\rm{co}}-n^{\rm{re}}),\\
            \calE_{n,M}^{\rm{ap}}=&\chi_1\calE_{n,M}^{\rm{co}}+(\p_t^2\chi_1+2\p_t\chi_1\p_t-\De_r\chi_1-2\p_r \chi_1\cdot\p_r )(n_M^{\rm{co}}-n^{\rm{re}})+\chi_1(1-\chi_1)\De_r(|\psi_M^{\rm{co}}|^2)\\
            &-(\De_r(\chi_1^2))|\psi_M^{\rm{co}}|^2-2(\p_r \chi_1^2)(\p_r |\psi_M^{\rm{co}}|^2).
        \end{aligned}
    \end{equation}

We estimate $n_M^{\rm{co}}-n^{\rm{re}} $ in the region $\{(t,r):\sqrt{3}(T-t)^{\frac{1}{2}}\leq r\leq 2\sqrt{3}(T-t)^{\frac{1}{2}}\}$. We write $n_M^{\rm{co}}-n^{\rm{re}}=n_M^{co}-n_M^{\rm{re}}+n_M^{\rm{re}}-n^{\rm{re}}$ and estimate $ n_M^{co}-n_M^{\rm{re}}$ first. In this region $n_M^{\rm{re}}(t,r)$ has the form ~\eqref{eq:nre-express}, which combined with  ~\eqref{eq:Nm-expansion-infty} and ~\eqref{eq:nre-initial} yields
\begin{equation}\label{eq:n-difference}
    \begin{aligned}
        n_M^{\rm{co}}-n_M^{\rm{re}} &=n_M^{\rm{co}}-n_{0,M}^{\rm{re}}+n_{0,M}^{\rm{re}}-n_M^{\rm{re}} \\
        &=\frac{1}{3 }\sum_{m=0}^M\sum_{j=1}^{[\frac{m}{2}]}(T-t)^{ 2j}\left(\frac{r}{\sqrt{3}}\right)^{\frac{m-4-6j}{2}}\ti N_{m,1,j}\left(\frac{1}{2}\log\frac{r}{\sqrt{3}}\right) \\
        &\quad + \frac{1}{3 }\sum_{m=0}^M \sum_{j=0}^{[\frac{m}{2}]}(T-t)^{1+2j}\left(\frac{r}{\sqrt{3}}\right)^{\frac{m-7-6j}{2}}\ti N_{m,2,j}\left(\frac{1}{2}\log\frac{r}{\sqrt{3}}\right) \\
        &\quad +\frac{1}{3}\sum_{m=0}^M(T-t)^{\frac{m-4}{3}} \ti N_{m,3}\left(\frac{1}{3}\log (T-t), \frac{r}{\sqrt{3}(T-t)^{\frac{2}{3}}} \right)\\
        &\quad +\frac{r+T-t}{2r}\left(n_{0,M}^{\rm{re}}(r)-n_{0,M}^{\rm{re}}(r+T-t)\right)\\
        &\quad +\frac{r-T+t}{2r}\left(n_{0,M}^{\rm{re}}(r)-n_{0,M}^{\rm{re}}(r-T+t)\right)\\
        &\quad +\frac{1}{2r}\int_{r-T+t}^{r+T-t}sn_{1,M}^{\rm{re}}(s)\ud s.
    \end{aligned}
\end{equation}
For the term in the last line of ~\eqref{eq:n-difference}, we use   ~\eqref{eq:nre-initial}, that $\ti N_{m,2,0}$ are  polynomials,   Taylor expansion and the identity $\p_r^{2j}(rf(r))=r\De^j_r f(r)$ to obtain
\begin{equation}\label{eq:n-diff-1}
    \begin{aligned}
        &\frac{1}{2r}\int_{r-T+t}^{r+T-t}sn_{1,M}^{\rm{re}}(s)\ud s\\
        =& -\frac{1}{3}\frac{1}{2r}\int_{r-T+t}^{r+T-t}\left(\sum_{m=0}^M\left(  \sum_{j=0}^{2[\frac{M-m}{2}]+1}\frac{1}{j!}\p_r^j\left(r\left(\frac{r}{\sqrt{3}}\right)^{\frac{m-7}{2}}\ti N_{m,2,0}\left(\frac{1}{2}\log\frac{r}{\sqrt{3}}\right) \right)(r)(s-r)^j  \right)+\calR^{\rm{ap}}_{1,M,0}(s)\right)\ud s\\
        =&-\frac{1}{3} \sum_{m=0}^M\sum_{j=0}^{[\frac{M-m}{2}]}\frac{3^{-\frac{m-7}{4}}}{(2j+1)!}\left( \De_r^j \left(r^{\frac{m-7}{2}}\ti N_{m,2,0}\left(\frac{1}{2}\log\frac{r}{\sqrt{3}}\right)\right)\right)(r) (T-t)^{2j+1}+\calR^{\rm{ap}}_{1,M}(t,r),
    \end{aligned}
\end{equation}
where $\calR^{\rm{ap}}_{1,M}(t,r) =-\frac{1}{6r}\int_{r-T+t}^{r+T-t}\calR^{\rm{ap}}_{1,M,0}(s)\ud s$ satisfies
\begin{equation}\label{eq:R1M-decay}
    \begin{aligned}
        |\p_t^{a_1}\p_r^{a_2}\calR^{\rm{ap}}_{1,M}(t,r)|\lesssim_{M,a_1,a_2} (T-t)^{\frac{M}{4}-a_1-a_2}(1+|\log(T-t)|)^{J_{1,M}},  
    \end{aligned}
\end{equation}
for some $J_{1,M}\in \Z_{\geq 0}$ and any $a_1,a_2\in \Z_{\geq 0}$.

Recall the notation $\zeta=\frac{1}{2}\log\frac{r}{\sqrt{3}}$. For any $\lam\in \R$ and function $f(\zeta)$, we have
\begin{align}
    \De_r(r^{\lam}f(\zeta))=r^{\lam-2}(\frac{1}{2}\p_{\zeta}+\lam)(\frac{1}{2}\p_{\zeta}+\lam+1)f(\zeta)
\end{align}
Then using ~\eqref{eq:tiN12-recursion} we have
\begin{equation}
    \begin{aligned}
        &\De_r^j \left(r^{\frac{m-7}{2}}\ti N_{m,2,0}\left(\frac{1}{2}\log\frac{r}{\sqrt{3}}\right)\right)(r)\\
        =&3^j(2j+1)!r^{\frac{m-7-4j}{2}}\ti N_{m+2j,2,j}\left(\frac{1}{2}\log\frac{r}{\sqrt{3}}\right),
    \end{aligned}
\end{equation}
which implies
\begin{equation}\label{eq:n-diff-2}
    \begin{aligned}
        &-\frac{1}{3} \sum_{m=0}^M\sum_{j=0}^{[\frac{M-m}{2}]}\frac{3^{-\frac{m-7}{4}}}{(2j+1)!}\left( \De_r^j \left(r^{\frac{m-7}{2}}\ti N_{m,2,0}\left(\frac{1}{2}\log\frac{r}{\sqrt{3}}\right)\right)\right)(r) (T-t)^{2j+1}\\
        =&-\frac{1}{3} \sum_{m=0}^M\sum_{j=0}^{[\frac{M-m}{2}]} 3^{j-\frac{m-7}{4}} r^{\frac{m-7-4j}{2}}\ti N_{m+2j,2,j}\left(\frac{1}{2}\log\frac{r}{\sqrt{3}}\right)(T-t)^{2j+1}\\
        =&-\frac{1}{3}\sum_{m=0}^M\sum_{j=0}^{[\frac{m}{2}]} 3^{-\frac{m-6j-7}{4}}r^{\frac{m-6j-7}{2}}\ti N_{m,2,j}\left(\frac{1}{2}\log\frac{r}{\sqrt{3}}\right)(T-t)^{2j+1}.
    \end{aligned}
\end{equation}

    For the terms in the second and the third to the last line of ~\eqref{eq:n-difference}, recalling ~\eqref{eq:nre-initial} and that $\ti N_{m,1,0}$ are  polynomials, using Taylor expansion we have
    \begin{equation}\label{eq:n-diff-3}
        \begin{aligned}
            &\frac{r+T-t}{2r}\left(n_{0,M}^{\rm{re}}(r)-n_{0,M}^{\rm{re}}(r+T-t)\right)+\frac{r-T+t}{2r}\left(n_{0,M}^{\rm{re}}(r)-n_{0,M}^{\rm{re}}(r-T+t)\right)\\
            =&\frac{1}{2r}(rn_{0,M}^{\rm{re}}(r)-(r+T-t)n_{0,M}^{\rm{re}}(r+T-t) )\\
            &+\frac{1}{2r}(rn_{0,M}^{\rm{re}}(r)-(r-T+t)n_{0,M}^{\rm{re}}(r-T+t) )\\
            =&-\frac{1}{3r}\sum_{m=0}^M\sum_{j=1}^{[\frac{M-m}{2}]}\frac{1}{(2j)!}\p_r^{2j}\left(r\left(\frac{r}{\sqrt{3}}\right)^{\frac{m-4}{2}}\ti N_{m,1,0}\left(\frac{1}{2}\log\frac{r}{\sqrt{3}}\right)\right)(r)(T-t)^{2j}+\calR^{\rm{ap}}_{0,M}(t,r)\\
            =&-\frac{1}{3}\sum_{m=0}^M\sum_{j=1}^{[\frac{M-m}{2}]}\frac{1}{(2j)!}\De_r^j\left( \left(\frac{r}{\sqrt{3}}\right)^{\frac{m-4}{2}}\ti N_{m,1,0}\left(\frac{1}{2}\log\frac{r}{\sqrt{3}}\right)\right)(r)(T-t)^{2j}+\calR^{\rm{ap}}_{0,M}(t,r),
        \end{aligned}
    \end{equation}
where $\calR^{\rm{ap}}_{0,M}$ satisfies
\begin{equation}\label{eq:R0M-decay}
    \begin{aligned}
        |\p_t^{a_1}\p_r^{a_2}\calR^{\rm{ap}}_{0,M}(t,r)|\lesssim_{M,a_1,a_2} (T-t)^{\frac{M}{4}-a_1-a_2-1}(1+|\log(T-t)|)^{J_{0,M}},  
    \end{aligned}
\end{equation}
for some $J_{0,M}\in \Z_{\geq 0}$ and any $a_1,a_2\in \Z_{\geq 0}$.

Using ~\eqref{eq:tiN12-recursion} we have
\begin{equation}
    \begin{aligned}
        &\De_r^j\left( r^{\frac{m-4}{2}}\ti N_{m,1,0}\left(\frac{1}{2}\log\frac{r}{\sqrt{3}}\right)\right)(r)\\
        =&3^j(2j)!r^{\frac{m-4-4j}{2}}\ti N_{m+2j,1,j}\left(\frac{1}{2}\log\frac{r}{\sqrt{3}}\right),
    \end{aligned}
\end{equation}
which implies
\begin{equation}\label{eq:n-diff-4}
    \begin{aligned}
        &-\frac{1}{3}\sum_{m=0}^M\sum_{j=1}^{[\frac{M-m}{2}]}\frac{1}{(2j)!}\De_r^j\left( \left(\frac{r}{\sqrt{3}}\right)^{\frac{m-4}{2}}\ti N_{m,1,0}\left(\frac{1}{2}\log\frac{r}{\sqrt{3}}\right)\right)(r)(T-t)^{2j}\\
        =&-\frac{1}{3}\sum_{m=0}^M\sum_{j=1}^{[\frac{M-m}{2}]} 3^{-\frac{m-4-4j}{4}} r^{\frac{m-4-4j}{2}}\ti N_{m+2j,1,j}\left(\frac{1}{2}\log\frac{r}{\sqrt{3}}\right)(T-t)^{2j}\\
        =&-\frac{1}{3}\sum_{m=0}^M\sum_{j=1}^{[\frac{m }{2}]} 3^{-\frac{m-6j-4}{4}} r^{\frac{m-6j-4}{2}}\ti N_{m ,1,j}\left(\frac{1}{2}\log\frac{r}{\sqrt{3}}\right)(T-t)^{2j}
    \end{aligned}
\end{equation}

Inserting ~\eqref{eq:n-diff-1}, ~\eqref{eq:n-diff-2}, ~\eqref{eq:n-diff-3} and ~\eqref{eq:n-diff-4} into ~\eqref{eq:n-difference}, we obtain that in the region $\{(t,r):\sqrt{3}(T-t)^{\frac{1}{2}}\leq r\leq2\sqrt{3}(T-t)^{\frac{1}{2}}\}$ there holds
\begin{equation}\label{eq:n-difference-1}
    \begin{aligned}
        n_M^{\rm{co}}-n_M^{\rm{re}} &=\frac{1}{3}\sum_{m=0}^M(T-t)^{\frac{m-4}{3}} \ti N_{m,3}\left(\frac{1}{3}\log (T-t), \frac{r}{\sqrt{3}(T-t)^{\frac{2}{3}}} \right)\\
        &\quad +\calR^{\rm{ap}}_{0,M}(t,r)+\calR^{\rm{ap}}_{1,M}(t,r).
    \end{aligned}
\end{equation}

For $n_M^{\rm{re}}-n^{\rm{re}}$, ~\eqref{eq:nre-diff-1} yields that
\begin{align}\label{eq:nre-diff-2}
    \|n_M^{\rm{re}}-n^{\rm{re}} \|_{W_{t,r}^{s,\infty}(\sqrt{3}(T-t)^{\frac{1}{2}}\leq r\leq 2\sqrt{3}(T-t)^{\frac{1}{2}})}\lesssim (T-t)^{\frac{M}{4}-s-2}
\end{align}

Recalling the expression ~\eqref{eq:E-expression}, the terms $\chi_1\calE^{\rm{co}}_{\psi,M}$ and $\chi_1\calE^{\rm{co}}_{n,M}$ can be estimated using ~\eqref{eq:E-core-decay}. By definition of $\chi_1$, the terms containing $\psi_M^{\rm{co}}$ are supported in the region $\{(t,r):\sqrt{3}(T-t)^{\frac{1}{2}}\leq r\leq2\sqrt{3}(T-t)^{\frac{1}{2}}\}$, which combined with ~\eqref{eq:core-app} and ~\eqref{eq:RN-decay} yields the desired estimates. The terms containing $n_M^{\rm{co}}-n^{\rm{re}}$ can be estimated using  ~\eqref{eq:tiNm3-decay}, ~\eqref{eq:R0M-decay}, ~\eqref{eq:R1M-decay},  ~\eqref{eq:n-difference-1} and ~\eqref{eq:nre-diff-2}. Then we obtain the $L^{\infty}$ estimates in ~\eqref{eq:E-estimate} for $M\geq 12s+24$. The $L^2$ estimates follows from the $L^{\infty}$ estimates and ~\eqref{eq:E-supp}. The proof is finished.

\end{proof}

\section{Construction of smooth solutions}\label{sec:fixed-point}
In this section we finish the proof of Theorem ~\ref{thm:main}.  We first construct solutions $(\psi,n)$ to ~\eqref{eq:Za} in the space $H^{s+2}\times H^{s+1}$ for any prescribed $s$ in Section ~\ref{subsec:linear}, Section ~\ref{subsec:nl} and Section ~\ref{subsec:finite-regularity}. Then we use the results in Lemma~\ref{lem:ap-diff} and the fixed point argument in Proposition ~\ref{prop:fixed-point} to prove that these solutions are the same for $t$ sufficiently close to $T$, which combined with the high regularity persistence in Lemma~\ref{lem:Cauchy} finishes the proof.

\subsection{Linear estimates}\label{subsec:linear}
Consider 
\begin{equation}\label{eq:uv}
    \begin{aligned}
    \begin{cases}
        \psi=\psi_M^{\rm{ap}}+u,\\
        n=n_M^{\rm{ap}}+v,
    \end{cases}
    \end{aligned}
\end{equation}
 with $u(T)=v(T)=\p_t v(T)=0$. Inserting ~\eqref{eq:uv} into ~\eqref{eq:Za}, we obtain
 \begin{equation}
     \begin{aligned}
     \left\{
         \begin{array}{rl}
              i\p_t u+\De_r u-n_M^{\rm{ap}}u-\psi_{M}^{\rm{ap}}v &=uv-\calE_{\psi,M}^{\rm{ap}},  \\
              \p_{tt}v-\De_r v-2\Re\De_r(\overline{\psi_{M}^{\rm{ap}}}u) &=\De_r(|u|^2)-\calE_{n,M}^{\rm{ap}}, 
         \end{array}
    \right.
     \end{aligned}
 \end{equation}
which leads us to consider the linear system
\begin{equation}\label{eq:UV}
    \begin{aligned}
     \left\{
         \begin{array}{rl}
              i\p_t U+\De_r U-n_M^{\rm{ap}}U-\psi_{M}^{\rm{ap}}V &=uv-\calE_{\psi,M}^{\rm{ap}},  \\
              \p_{tt}V-\De_r V-2\Re\De(\overline{\psi_{M}^{\rm{ap}}}U) &=\De(|u|^2)-\calE_{n,M}^{\rm{ap}}.
         \end{array}
    \right.
    \end{aligned}
\end{equation}
We introduce the following function space. Take an integer $s\gg1$ sufficiently large. Take $T_1<T$ such that Proposition ~\ref{prop:core} and Lemma ~\ref{lem:E} hold for $T_1\leq t<T$. For $K\geq 1$ and an interval $I\subset[T_1,T)$, 
we denote by $\calH_K^s(I)$ the space contains all pairs of  functions $(f(t,r),g(t,r))$   satisfying 
\begin{equation}
        \begin{aligned}
            f\in& C_{\rm{loc}}(I;H^{s+2}(\R^3,\C))\cap C_{\rm{loc}}^1(I;H^{s}(\R^3,\C)),\\
            g\in&C_{\rm{loc}}(I;H^{s+1}(\R^3,\R))\cap C_{\rm{loc}}^1(I;H^{s}(\R^3,\R)).
        \end{aligned}
\end{equation}
The space $\calH_K^s(I)$ is defined as a real Banach space equipped with the norm
\begin{equation}\label{eq:HKs-norm}
    \begin{aligned}
        \|(f,g)\|_{\calH_K^s(I)}=&\sup_{t\in I}(T-t)^{-K}\left(\|f\|_{H^{s+2}}+\|\p_t f\|_{H^{s}}+\|g\|_{H^{s+1}}+\|\p_t g\|_{H^{s}}\right).
    \end{aligned}    
\end{equation}
Specifically, we write $\calH_K^s=\calH_K^s([T_1,T))$. For $(u,v)\in \calH_K^s$, combining with ~\eqref{eq:E-estimate} and $s\geq 2>\frac{3}{2}$ we have
\begin{equation}\label{eq:source-estimate}
    \begin{aligned}
        \|uv\|_{H^{s+1}}\lesssim &\|u\|_{H^{s+1}}\|v\|_{H^{s+1}}\lesssim (T-t)^{2K}\|(u,v)\|_{\calH_K^s}^2,\\
        \|\p_t(uv)\|_{H^s}\lesssim & \|u\|_{H^s}\|\p_t v\|_{H^s}+\|\p_tu\|_{H^{s}}\|v\|_{H^s}\lesssim (T-t)^{2K}\|(u,v)\|_{\calH_K^s}^2,\\
        \|\calE_{\psi,M}^{\rm{ap}}\|_{H^{s+1}}\lesssim & (T-t)^{\frac{M}{4}-6-2s}  ,\\
        \|\p_t\calE_{\psi,M}^{\rm{ap}}\|_{H^{s}}\lesssim & (T-t)^{\frac{M}{4}-6-2s} ,\\
        \|\De_r(|u|^2)\|_{H^s}\lesssim& \|u\|_{H^s}\|\De_r u\|_{H^s}+\|\na_x u\|_{H^s}^2\lesssim(T-t)^{2K}\|(u,v)\|_{\calH_K^s}^2 ,\\
        \|\calE_{n,M}^{\rm{ap}}\|_{H^{s}}\lesssim& (T-t)^{\frac{M}{4}-6-2s}.
    \end{aligned}
\end{equation}

\begin{lem}\label{lem:UV}
    For every even integer $s\geq 2$, there exists $K_s>0$ such that for every $K\geq K_s$, there exists $M_{K,s}$ such that for every $M\geq M_{K,s}$, such that there exists $T_{1,s,K,M}<T$, such that for every $T_1$ satisfying $T_{1,s,K,M}<T_1<T$,  the following properties hold: For radial functions $(F,G)$ such that
    \begin{equation}\label{eq:FG-assumption}
        \begin{aligned}
             F\in& C_{\rm{loc}}([T_1,T);H^{s+1}(\R^3,\C))\cap C_{\rm{loc}}^1([T_1,T);H^{s}(\R^3,\C)),\\
                G\in & C_{\rm{loc}}([T_1,T);H^{s}(\R^3,\R)) ,\\
        \end{aligned}
    \end{equation}
    and that there exists a constant $C_0>0$ with
    \begin{equation}\label{eq:FG}
        \begin{aligned}
             \|F(t)\|_{H^{s+1}}+\|\p_t F(t)\|_{H^s}+\|G(t)\|_{H^s}\leq C_0 (T-t)^{2K},
        \end{aligned}
    \end{equation}
     there exists a unique pair of solutions $(U,V)\in \calH_K^s$ to 
    \begin{equation}\label{eq:UV-FG}
        \begin{aligned}
            \left\{
         \begin{array}{rl}
              i\p_t U+\De_r U-n_M^{\rm{ap}}U-\psi_{M}^{\rm{ap}}V &=F,  \\
              \p_{tt}V-\De_r V-2\Re\De_r(\overline{\psi_{M}^{\rm{ap}}}U) &=G.
         \end{array}
    \right.
        \end{aligned}
    \end{equation}
    with the following estimates:
    \begin{equation}\label{eq:UV-estimate}
        \begin{aligned}
           \|(U,V)\|_{\calH_K^s}\lesssim C_0(T-T_1)^{K-s-2}.
        \end{aligned}
    \end{equation}
    
\end{lem}

\begin{proof}[Proof of Lemma~\ref{lem:UV}]
    We use the classical method of backwards in time integration of the flow. Recall that $\tau=(T-t)^{\frac{1}{3}}$ and $\eta=\frac{r}{\sqrt{3}(T-t)^{\frac{2}{3}}}$. 

    \noindent
    \textit{Step 1.}   In this step we introduce new unknowns and rearrange the system ~\eqref{eq:UV-FG}. Firstly define
    \begin{equation}\label{eq:ti-def}
        \begin{aligned}
            (\ti U,\ti F,\ti \psi_{M}^{\rm{ap}},\ti \calE_{\psi,M}^{\rm{ap}})&=(  e^{-i\tau^{-1}}U,e^{-i\tau^{-1}} F,  e^{-i\tau^{-1}}\psi_{M}^{\rm{ap}},  e^{-i\tau^{-1}}\calE_{\psi,M}^{\rm{ap}}),\\
            (\ti V,\ti G,\ti n_{M}^{\rm{ap}},\ti \calE_{n,M}^{\rm{ap}})&=(  V,  G,  n_{M}^{\rm{ap}},  \calE_{n,M}^{\rm{ap}}).
        \end{aligned}
    \end{equation}
    Then ~\eqref{eq:UV-FG} turns into
    \begin{equation}\label{eq:UV-FG-1}
        \begin{aligned}
               \left\{
         \begin{array}{rl}
              i\p_t \ti U+\calL_{k,\tau}\ti U-\ti \calV_M\ti U-\ti \psi_{M}^{\rm{ap}}\ti V &=\ti F,  \\
              \p_{tt} \ti V-\De_r\ti V-2\Re\De_r(\overline{\ti \psi_{M}^{\rm{ap}}}\ti U) &=\ti G,
         \end{array}
    \right.
        \end{aligned}
    \end{equation}
    where
    \begin{equation}
        \begin{aligned}
            \calL_{k,\tau}= \De_r-\frac{1}{3\tau^4} -\frac{1}{3\tau^4} \ti N_k\left(\eta\right),\  
            \ti \calV_M= \ti n_M^{\rm{ap}}(t,r)-\frac{1}{3\tau^4} \ti N_k\left(\eta\right).
        \end{aligned}
    \end{equation}
    Note that for fixed $t$,
    \begin{align}\label{eq:Lk-Lkveps}
        (\calL_{k,\tau}f)(\sqrt{3}\tau^2\eta)=\frac{1}{3\tau^4}(\calL_k(f(\sqrt{3}\tau^2\cdot)))(\eta).
    \end{align}
    
    Next we decompose $\ti U$ into a finite dimensional component in ${\rm{Span}}_{\R}\{\ti \psi_M^{\rm{ap}},i\ti \psi_M^{\rm{ap}} \}$ and a infinite dimensional component orthogonal to $\ker \calL_{k,\tau}$, and then derive the evolution equations for them.  
    
    For function $f$ we define
    \begin{align}\label{eq:PM-def}
        \calP_M  f=f-\al_f\ti \psi_M^{\rm{ap}}-\be_fi\ti  \psi_M^{\rm{ap}},
    \end{align}
    where $\al_f$ and $\be_f$ are uniquely determined such that  $\La \calP_M f,\tau^{-3}\ti P_k (\eta)\Ra=\La \calP_M f,i\tau^{-3}\ti P_k (\eta) \Ra=0$, i.e.
    \begin{equation}\label{eq:albe-ortho}
        \begin{aligned}
            \begin{cases}
                \al_f\La \ti \psi_M^{\rm{ap}},\tau^{-3}\ti P_k (\eta)\Ra+\be_f  \La i\ti \psi_M^{\rm{ap}},\tau^{-3}\ti P_k (\eta)\Ra&=\La f,\tau^{-3}\ti P_k (\eta)\Ra,\\
                \al_f\La \ti \psi_M^{\rm{ap}},i\tau^{-3}\ti P_k (\eta)\Ra+\be_f  \La i\ti \psi_M^{\rm{ap}},i\tau^{-3}\ti P_k (\eta)\Ra&=\La f,i\tau^{-3}\ti P_k (\eta)\Ra.
            \end{cases}
        \end{aligned}
    \end{equation}
    The existence and uniqueness of $\al_f,\be_f$ follows from the identities
    \begin{align}
        \label{eq:albe-existence}
        \La \ti \psi_M^{\rm{ap}},\tau^{-3}\ti P_k (\eta)\Ra^2=27\|\ti P_k(\cdot)\|_{L^2}^4+o_{t\to T}(1),\  \La i\ti \psi_M^{\rm{ap}},\tau^{-3}\ti P_k (\eta)\Ra^2=o_{t\to T}(1).
    \end{align}
    Define $\al(t)=\al_{\ti U(t)}$, $\be(t)=\be_{\ti U(t)}$. A direct calculation of ~\eqref{eq:UV-FG-1} yields
    \begin{equation}\label{eq:PU}
        \begin{aligned}
            \left(i\p_t+\calL_{k,\tau}-\ti \calV_M\right)\calP_M \ti U&=\ti\psi_{M}^{\rm{ap}}\ti V+\ti F-(\al+i\be)\ti \calE_{\psi,M}^{\rm{ap}}-(-\p_t\be+i\p_t\al)\ti\psi_M^{\rm{ap}}
        \end{aligned}
    \end{equation}
    Differentiating the orthogonal conditions ~\eqref{eq:albe-ortho} for $f=\ti U$ by $t$  and using ~\eqref{eq:PU} yields the ode system for $\al(t),\be(t)$:
    \begin{equation}\label{eq:mod}
        \begin{aligned}
            \left\{
            \begin{array}{l}
                 \p_t\al\La \ti \psi_M^{\rm{ap}}, \tau^{-3}\ti P_k (\eta)\Ra+\p_t\be\La i\ti \psi_M^{\rm{ap}}, \tau^{-3}\ti P_k (\eta) \Ra= \La \calP_M \ti U,\p_t\left(\tau^{-3}\ti P_k (\eta)\right) \Ra \\
                 \quad +\La \ti\calV_M\calP_M \ti U+\ti\psi_M^{\rm{ap}}\ti V+\ti F ,i\tau^{-3}\ti P_k (\eta)  \Ra
                 -\La (\al+i\be)\ti \calE_{\psi,M}^{\rm{ap}} ,i\tau^{-3}\ti P_k (\eta)  \Ra,\\
                 \p_t\al\La \ti \psi_M^{\rm{ap}}, i\tau^{-3}\ti P_k (\eta)\Ra+\p_t\be\La i\ti \psi_M^{\rm{ap}}, i\tau^{-3}\ti P_k (\eta) \Ra= \La \calP_M \ti U,\p_t\left(i\tau^{-3}\ti P_k (\eta)\right) \Ra\\
                 \quad -\La \ti\calV_M\calP_M \ti U+\ti\psi_M^{\rm{ap}}\ti V+\ti F ,\tau^{-3}\ti P_k (\eta)  \Ra
                 +\La (\al+i\be)\ti \calE_{\psi,M}^{\rm{ap}} ,\tau^{-3}\ti P_k (\eta)  \Ra.
            \end{array}
            \right.
        \end{aligned}   
    \end{equation}
    Then we derive the equation for $\calP_M\ti U$ without $\p_t\al,\p_t\be$ from ~\eqref{eq:PU}. Using ~\eqref{eq:mod} we have
    \begin{equation}
        \begin{aligned}
            \left\{
            \begin{array}{l}
                  \La (\p_t\al+i\p_t\be)\ti \psi_M^{\rm{ap}}, \tau^{-3}\ti P_k (\eta)\Ra = -\La (\p_t\calP_M) \ti U, \tau^{-3}\ti P_k (\eta)\Ra \\
                 \quad-\La i\left( \ti\calV_M\calP_M \ti U+\ti\psi_M^{\rm{ap}}\ti V+\ti F-(\al+i\be)\ti \calE_{\psi,M}^{\rm{ap}} \right), \tau^{-3}\ti P_k (\eta)  \Ra,\\
                 \La  (\p_t\al+i\p_t\be)\ti \psi_M^{\rm{ap}}, i\tau^{-3}\ti P_k (\eta)\Ra =- \La (\p_t\calP_M) \ti U, i\tau^{-3}\ti P_k (\eta)\Ra\\
                 \quad-\La i\left(\ti\calV_M\calP_M \ti U+\ti\psi_M^{\rm{ap}}\ti V+\ti F -(\al+i\be)\ti \calE_{\psi,M}^{\rm{ap}}\right),i\tau^{-3}\ti P_k (\eta)  \Ra,
            \end{array}
            \right.
        \end{aligned}
    \end{equation}
    which combined with ~\eqref{eq:albe-existence},  $\calP_M^2=\calP_M$ and the fact that $\ti \psi_M^{\rm{ap}}\in\ker \calP_M $ yields
    \begin{equation}\label{eq:mod-1}
        \begin{aligned}
            (\p_t\al+i\p_t\be)\ti \psi_M^{\rm{ap}}=&-(\Id-\calP_M)(\p_t\calP_M) \ti U-i(\Id-\calP_M)\left(\ti\calV_M\calP_M \ti U+\ti\psi_M^{\rm{ap}}\ti V+\ti F-(\al+i\be)\ti \calE_{\psi,M}^{\rm{ap}} \right)\\
            =&-(\p_t\calP_M)\calP_M\ti U-i(\Id-\calP_M)\left(\ti\calV_M\calP_M \ti U+\ti\psi_M^{\rm{ap}}\ti V+\ti F-(\al+i\be)\ti \calE_{\psi,M}^{\rm{ap}} \right).
        \end{aligned}
    \end{equation}
    Inserting ~\eqref{eq:mod-1} into ~\eqref{eq:PU} yields
    \begin{equation}\label{eq:PU-1}
        \begin{aligned}
            \left(i\p_t+\calL_{k,\tau}\right)\calP_M \ti U&= \ti \calV_M \calP_M \ti U+\ti\psi_{M}^{\rm{ap}}\ti V+\ti F-(\al+i\be)\ti \calE_{\psi,M}^{\rm{ap}} \\
            &\quad +i\left((\p_t\calP_M)\calP_M\ti U+i(\Id-\calP_M)\left(\ti\calV_M\calP_M \ti U+\ti\psi_M^{\rm{ap}}\ti V+\ti F-(\al+i\be)\ti \calE_{\psi,M}^{\rm{ap}} \right)\right)\\
            &=i (\p_t\calP_M)\calP_M\ti U+\calP_M\left(\ti\calV_M\calP_M \ti U+\ti\psi_M^{\rm{ap}}\ti V+\ti F-(\al+i\be)\ti \calE_{\psi,M}^{\rm{ap}}\right).
        \end{aligned}
    \end{equation}
    The term  $\ti\psi_M^{\rm{ap}}\ti V+\ti F$ is of lower regularity than $\ti U$, and we introduce a new unknown to absorb it.
    
    For $f(r)\in L^2 $ such that $\La f ,\tau^{-3}\ti P_k(\eta)\Ra=\La f ,i\tau^{-3}\ti P_k(\eta)\Ra=0$, we denote by $(\calL_{k,\tau})^{-1}f$ the unique solution $u(r)\in H^2$ such that $\calL_{k,\tau}u=f$ and $\La u ,\tau^{-3}\ti P_k(\eta)\Ra=\La u ,i\tau^{-3}\ti P_k(\eta)\Ra=0$. The existence and uniqueness of $u$ are guaranteed by Lemma~\ref{lem:Lk} and ~\eqref{eq:Lk-Lkveps}. Define
    \begin{align}\label{eq:W-def}
        \ti W=\calP_M\ti U-(\calL_{k,\tau})^{-1}\left(\calP_{M}\left( \ti\psi_M^{\rm{ap}}\ti V+\ti F\right)\right).
    \end{align}
    Then $\calP_M\ti W=\ti W$, and ~\eqref{eq:PU-1} reduces to
    \begin{equation}\label{eq:W-1}
        \begin{aligned}
            (i\p_t+\calL_{k,\tau})\ti W=&-i\p_t\left((\calL_{k,\tau})^{-1}\left(\calP_{M}\left( \ti\psi_M^{\rm{ap}}\ti V+\ti F\right)\right) \right)\\
            &+i(\p_t\calP_M)\left(\ti W+(\calL_{k,\tau})^{-1}\left(\calP_{M}\left( \ti\psi_M^{\rm{ap}}\ti V+\ti F\right)\right)\right)\\
            &+\calP_M\left(\ti \calV_M\left(\ti W+(\calL_{k,\tau})^{-1}\left(\calP_{M}\left( \ti\psi_M^{\rm{ap}}\ti V+\ti F\right)\right)\right)\right)\\
            &-\calP_M\left((\al+i\be)\ti\calE_{\psi,M}^{\rm{ap}}\right).
        \end{aligned}
    \end{equation}
    The second line of ~\eqref{eq:UV-FG-1} reduces to
    \begin{equation}\label{eq:V-1}
        \begin{aligned}
            \p_{tt} \ti V-\De_r\ti V &=2\Re\De_r\left(\overline{\ti \psi_{M}^{\rm{ap}}}\left(\ti W+(\calL_{k,\tau})^{-1}\left(\calP_{M}\left( \ti\psi_M^{\rm{ap}}\ti V+\ti F\right)\right)\right)\right) +2\al\De_r\left|\ti\psi_M^{\rm{ap}}\right|^2+\ti G.
        \end{aligned}
    \end{equation}

    \noindent
    \textit{Step 2.} In this section we prove that for any $T_1<T_2<T$, the solution $(U,V)$ to ~\eqref{eq:UV-FG} exists for $T_1\leq t\leq T_2$. Using ~\eqref{eq:ti-def}, ~\eqref{eq:UV-FG-1}, ~\eqref{eq:mod}, ~\eqref{eq:W-1} and ~\eqref{eq:V-1}, the system ~\eqref{eq:UV-FG} reduces to the system of $(\ti W,\ti V,\al,\be)$ as follows.
    \begin{equation}\label{eq:UV-reduced}
        \begin{aligned}
            \left\{
            \begin{array}{rl}
            (i\p_t+\De_r)\ti W= & \calF^W_{\rm{L}}+\calF^W_0,\\
                 (\p_{tt}-\De_r)\ti V=& \calF^V_{\rm{L}}+\calF^V_0,\\
                 \p_t\al=& \calF^{\al}_{\rm{L}}+\calF^{\al}_0,\\
                 \p_t\be=& \calF^{\be}_{\rm{L}}+\calF^{\be}_0.
            \end{array}
            \right.
        \end{aligned}
    \end{equation}
    Where $\calF_{*}^{*}$ are given by
    \begin{align}
            \calF^{W}_{\rm{L}}=&\left(\frac{1}{3\tau^4} +\frac{1}{3\tau^4} \ti N_k\left(\eta\right)\right)\ti W
            -i\p_t\left((\calL_{k,\tau})^{-1}\left(\calP_{M}\left( \ti\psi_M^{\rm{ap}}\ti V \right)\right) \right)\\
            &+i(\p_t\calP_M)\left(\ti W+(\calL_{k,\tau})^{-1}\left(\calP_{M}\left( \ti\psi_M^{\rm{ap}}\ti V \right)\right)\right)\\
            &+\calP_M\left(\ti \calV_M\left(\ti W+(\calL_{k,\tau})^{-1}\left(\calP_{M}\left( \ti\psi_M^{\rm{ap}}\ti V \right)\right)\right)\right)\\
            &-\calP_M\left((\al+i\be)\ti\calE_{\psi,M}^{\rm{ap}}\right),\\
            \calF^{W}_0=&-i\p_t\left((\calL_{k,\tau})^{-1}\left(\calP_{M} \ti F \right) \right)
            +i(\p_t\calP_M)\left( (\calL_{k,\tau})^{-1}\left(\calP_{M} \ti F \right)\right)\\
            &+\calP_M\left(\ti \calV_M\left( (\calL_{k,\tau})^{-1}\left(\calP_{M}     \ti F \right)\right)\right),\\
            \calF^{V}_{\rm{L}}=&2\Re\De_r\left(\overline{\ti \psi_{M}^{\rm{ap}}}\left(\ti W+(\calL_{k,\tau})^{-1}\left(\calP_{M}\left( \ti\psi_M^{\rm{ap}}\ti V \right)\right)\right)\right) +2\al\De_r\left|\ti\psi_M^{\rm{ap}}\right|^2 ,\\
            \calF^V_{0}=&2\Re\De_r\left(\overline{\ti \psi_{M}^{\rm{ap}}}\left( (\calL_{k,\tau})^{-1}\left(\calP_{M}  \ti F \right)\right)\right) + \ti G,\\
            \calF^{\al}_{\rm{L}}=&\frac{1}{\La \ti \psi_M^{\rm{ap}},\tau^{-3}\ti P_k (\eta)\Ra^2+\La i\ti \psi_M^{\rm{ap}}, \tau^{-3}\ti P_k (\eta)\Ra^2}\bigg[\La \ti \psi_M^{\rm{ap}},\tau^{-3}\ti P_k (\eta)\Ra\\
                &\quad \times\bigg(\La \ti W+(\calL_{k,\tau})^{-1}\left(\calP_M\left(\ti \psi_{M}^{\rm{ap}}\ti V \right)\right),\p_t\left(\tau^{-3}\ti P_k (\eta)\right) \Ra \\
                 &\quad +\La \ti\calV_M\left(\ti W+(\calL_{k,\tau})^{-1}\left(\calP_M\left(\ti \psi_{M}^{\rm{ap}}\ti V \right)\right)\right)+\ti\psi_M^{\rm{ap}}\ti V ,i\tau^{-3}\ti P_k (\eta)  \Ra\\
                 &\quad -\La (\al+i\be)\ti \calE_{\psi,M}^{\rm{ap}} ,i\tau^{-3}\ti P_k (\eta)  \Ra\bigg) -\La i\ti \psi_M^{\rm{ap}}, \tau^{-3}\ti P_k (\eta)\Ra\\
                 &\quad \times \bigg(\La \ti W+(\calL_{k,\tau})^{-1}\left(\calP_M\left(\ti \psi_{M}^{\rm{ap}}\ti V \right)\right),\p_t\left(i\tau^{-3}\ti P_k (\eta)\right) \Ra\\
                 &\quad -\La \ti\calV_M\left(\ti W+(\calL_{k,\tau})^{-1}\left(\calP_M\left(\ti \psi_{M}^{\rm{ap}}\ti V \right)\right)\right)+\ti\psi_M^{\rm{ap}}\ti V  ,\tau^{-3}\ti P_k (\eta)  \Ra\\
                 &\quad +\La (\al+i\be)\ti \calE_{\psi,M}^{\rm{ap}} ,\tau^{-3}\ti P_k (\eta)  \Ra\bigg)\bigg],\\
                 \calF_0^{\al}=&\frac{1}{\La \ti \psi_M^{\rm{ap}},\tau^{-3}\ti P_k (\eta)\Ra^2+\La i\ti \psi_M^{\rm{ap}}, \tau^{-3}\ti P_k (\eta)\Ra^2}\bigg[\La \ti \psi_M^{\rm{ap}},\tau^{-3}\ti P_k (\eta)\Ra\\
                &\quad \times\bigg(\La  (\calL_{k,\tau})^{-1}\left(\calP_M \ti F \right),\p_t\left(\tau^{-3}\ti P_k (\eta)\right) \Ra \\
                 &\quad +\La \ti\calV_M\left( (\calL_{k,\tau})^{-1}\left(\calP_M  \ti F \right)\right) +\ti F ,i\tau^{-3}\ti P_k (\eta)  \Ra \bigg)\\
                 &\quad  -\La i\ti \psi_M^{\rm{ap}}, \tau^{-3}\ti P_k (\eta)\Ra\\
                 &\quad \times \bigg(\La  (\calL_{k,\tau})^{-1}\left(\calP_M \ti F\right),\p_t\left(i\tau^{-3}\ti P_k (\eta)\right) \Ra\\
                 &\quad -\La \ti\calV_M\left( (\calL_{k,\tau})^{-1}\left(\calP_M \ti F\right)\right) +\ti F ,\tau^{-3}\ti P_k (\eta)  \Ra\bigg)\bigg],\\
                             \calF^{\be}_{\rm{L}}=&\frac{1}{\La \ti \psi_M^{\rm{ap}},\tau^{-3}\ti P_k (\eta)\Ra^2+\La i\ti \psi_M^{\rm{ap}}, \tau^{-3}\ti P_k (\eta)\Ra^2}\bigg[\La \ti \psi_M^{\rm{ap}},\tau^{-3}\ti P_k (\eta)\Ra\\
                &\quad \times\bigg(\La \ti W+(\calL_{k,\tau})^{-1}\left(\calP_M\left(\ti \psi_{M}^{\rm{ap}}\ti V \right)\right),\p_t\left(i\tau^{-3}\ti P_k (\eta)\right) \Ra \\
                 &\quad -\La \ti\calV_M\left(\ti W+(\calL_{k,\tau})^{-1}\left(\calP_M\left(\ti \psi_{M}^{\rm{ap}}\ti V \right)\right)\right)+\ti\psi_M^{\rm{ap}}\ti V , \tau^{-3}\ti P_k (\eta)  \Ra\\
                 &\quad +\La (\al+i\be)\ti \calE_{\psi,M}^{\rm{ap}} , \tau^{-3}\ti P_k (\eta)  \Ra\bigg) +\La i\ti \psi_M^{\rm{ap}}, \tau^{-3}\ti P_k (\eta)\Ra\\
                 &\quad \times \bigg(\La \ti W+(\calL_{k,\tau})^{-1}\left(\calP_M\left(\ti \psi_{M}^{\rm{ap}}\ti V \right)\right),\p_t\left( \tau^{-3}\ti P_k (\eta)\right) \Ra\\
                 &\quad +\La \ti\calV_M\left(\ti W+(\calL_{k,\tau})^{-1}\left(\calP_M\left(\ti \psi_{M}^{\rm{ap}}\ti V \right)\right)\right)+\ti\psi_M^{\rm{ap}}\ti V  ,i\tau^{-3}\ti P_k (\eta)  \Ra\\
                 &\quad -\La (\al+i\be)\ti \calE_{\psi,M}^{\rm{ap}} ,i\tau^{-3}\ti P_k (\eta)  \Ra\bigg)\bigg],\\
                 \calF_0^{\be}=&\frac{1}{\La \ti \psi_M^{\rm{ap}},\tau^{-3}\ti P_k (\eta)\Ra^2+\La i\ti \psi_M^{\rm{ap}}, \tau^{-3}\ti P_k (\eta)\Ra^2}\bigg[\La \ti \psi_M^{\rm{ap}},\tau^{-3}\ti P_k (\eta)\Ra\\
                &\quad \times\bigg(\La  (\calL_{k,\tau})^{-1}\left(\calP_M \ti F \right),\p_t\left(i\tau^{-3}\ti P_k (\eta)\right) \Ra \\
                 &\quad -\La \ti\calV_M\left( (\calL_{k,\tau})^{-1}\left(\calP_M\  \ti F \right)\right) +\ti F , \tau^{-3}\ti P_k (\eta)  \Ra \bigg)\\
                 &\quad  +\La i\ti \psi_M^{\rm{ap}}, \tau^{-3}\ti P_k (\eta)\Ra\\
                 &\quad \times \bigg(\La  (\calL_{k,\tau})^{-1}\left(\calP_M \ti F\right),\p_t\left( \tau^{-3}\ti P_k (\eta)\right) \Ra\\
                 &\quad +\La \ti\calV_M\left( (\calL_{k,\tau})^{-1}\left(\calP_M \ti F\right)\right) +\ti F ,i\tau^{-3}\ti P_k (\eta)  \Ra\bigg)\bigg].
    \end{align}
    The terms $\calF_{\rm{L}}^{*}$ are $\R$-linear in $(\ti W,\ti V,\al,\be)$, while the terms $\calF_0^{*}$ are independent of them.
    
    The relation between $(U,V)$ and $(\ti W,\ti V,\al,\be)$ is given by
    \begin{equation}
        \begin{aligned}
            U=e^{i\tau^{-1}}\left( \ti W+(\calL_{k,\tau})^{-1}\left(\calP_M\left(\ti \psi_{M}^{\rm{ap}}\ti V+\ti F\right)\right)+(\al+i\be) \ti \psi_{M}^{\rm{ap}}  \right),
            \quad V=\ti V.
        \end{aligned}
    \end{equation}
    \begin{claim}\label{claim:UV-local}
        For any $T_1<T_2<T$ and any $\ti W_0,\ti V_0,\ti V_1,\al_0,\be_0$ such that
        \begin{equation}
            \begin{aligned}
                \calP_M|_{T_2}\ti W_0=\ti W_0,\ \ti W_0\in H^{s+2}(\R^3,\C), \ \ti V_0\in H^{s+1}(\R^3,\R),\ \ti V_1\in H^{s }(\R^3,\R),\ \al_0,\be_0\in \R,
            \end{aligned}
        \end{equation}
        there exists a unique pair of solutions $(\ti W,\ti V,\al,\be)$ to ~\eqref{eq:UV-reduced} such that 
        \begin{equation}\label{eq:WVab-initial}
            \begin{aligned}
                 \ti W(T_2)=\ti W_0,\ \ti V(T_2)=\ti V_0,\ \p_t\ti V(T_2)=\ti V_1,\ \al(T_2)=\al_0,\ \be(T_2)=\be_0,
            \end{aligned}
        \end{equation}
        \begin{equation}\label{eq:WVab-regularity}
            \begin{aligned}
                \ti W\in& C([T_1,T_2];H^{s+2}(\R^3,\C))\cap C^1([T_1,T_2];H^{s}(\R^3,\C)),\\
                \ti V\in & C([T_1,T_2];H^{s+1}(\R^3,\R))\cap C^1([T_1,T_2];H^{s}(\R^3,\R)),\\
                &\al(t)\in C^1([T_1,T_2]),\ \be(t)\in  C^1([T_1,T_2]),
            \end{aligned}
        \end{equation}
        and $\calP_M|_t\ti W(t)=\ti W(t)$ for all $T_1\leq t\leq T_2$.
    \end{claim}
    \begin{proof}[Proof of Claim~\ref{claim:UV-local}]
    Define the solution map  $\vec \Phi=(\Phi_1,\Phi_2,\Phi_3,\Phi_4)$ as follows.
    \begin{align}
        \Phi_1(\ti W ,\ti V  ,\al ,\be )=&e^{i(t-T_2)\De_r}\ti W_0+i\int_{t}^{T_2}e^{i(t-t' )\De_r}(\calF_{\rm{L}}^W+\calF_0^{W})(t' )\ud t' ,\\
        \Phi_2(\ti W ,\ti V  ,\al ,\be)=&\cos((t-T_2) \sqrt{-\De_r})\ti V_0+\frac{\sin((t-T_2) \sqrt{-\De_r})}{ \sqrt{-\De_r}}\ti V_1\\
        &-\int_t^{T_2}\frac{\sin((t-t' ) \sqrt{-\De_r})}{ \sqrt{-\De_r}}(\calF_{\rm{L}}^V+\calF_0^{V})(t')\ud t' ,\\
        \Phi_3(\ti W ,\ti V , \al ,\be)=&\al_0-\int_t^{T_2}(\calF_{\rm{L}}^{\al}+\calF_{0}^{\al})(t')\ud t',\\
        \Phi_4(\ti W ,\ti V , \al ,\be)=&\be_0-\int_t^{T_2}(\calF_{\rm{L}}^{\be}+\calF_{0}^{\be})(t')\ud t'.
    \end{align}

    Define the Banach space for radial functions $\ti W,\ti V$  and $\al,\be$
    \begin{equation}
        \begin{aligned}
            \|(\ti W,\ti V,\al,\be)\|_{X_{T_1,T_2}} =&\|\ti W\|_{C([T_1,T_2];H^{s+2}(\R^3,\C)) }\\
            &+\|\ti V\|_{C([T_1,T_2];H^{s+1}(\R^3,\R)) }+\|\ti V\|_{  C^1 ([T_1,T_2];H^{s}(\R^3,\R))}\\
            &+\|\al\|_{C ([T_1,T_2])}+\|\be\|_{C([T_1,T_2])}.
        \end{aligned}
    \end{equation}
    We prove that $\vec \Phi$ is contractive in a closed ball of $X_{T_1,T_2}$ for $T_2-T_1$ sufficiently small.

    We first give the estimates of $\calF_{*}^*$. We estimate the bound of $\p_t((\calL_{k,\tau})^{-1}\calP_M)$ as follows:
    Let $(\calL_{k,\tau})^{-1}\calP_M f=u$, we have
    \begin{equation}\label{eq:LP-t}
        \begin{aligned}
            &\left(\p_t((\calL_{k,\tau})^{-1}\calP_M)\right)f\\
            =&\p_tu-(\calL_{k,\tau})^{-1}\calP_M\p_t f\\
            =&(\calL_{k,\tau})^{-1}\left(\calL_{k,\tau}\p_tu-\calP_M\p_t f\right)+\frac{\La \p_t u, \tau^{-3}\ti P_k(\eta)\Ra}{\La \tau^{-3}\ti P_k(\eta),\tau^{-3}\ti P_k(\eta)\Ra}\tau^{-3}\ti P_k(\eta)+\frac{\La \p_t u, i\tau^{-3}\ti P_k(\eta)\Ra}{\La i\tau^{-3}\ti P_k(\eta),i\tau^{-3}\ti P_k(\eta)\Ra}i\tau^{-3}\ti P_k(\eta)\\
            =&(\calL_{k,\tau})^{-1}\left(-(\p_t\calL_{k,\tau})(\calL_{k,\tau})^{-1}\calP_M f+(\p_t\calP_M) f\right)\\
            &-\frac{\La u, \p_t(\tau^{-3}\ti P_k(\eta))\Ra}{\La \tau^{-3}\ti P_k(\eta),\tau^{-3}\ti P_k(\eta)\Ra}\tau^{-3}\ti P_k(\eta)-\frac{\La  u, i\p_t(\tau^{-3}\ti P_k(\eta))\Ra}{\La i\tau^{-3}\ti P_k(\eta),i\tau^{-3}\ti P_k(\eta)\Ra}i\tau^{-3}\ti P_k(\eta).
        \end{aligned}
    \end{equation}
    Then there holds
    \begin{equation}\label{eq:operator-L-1-P}
        \begin{aligned}
            \|\left(\p_t((\calL_{k,\tau})^{-1}\calP_M)\right)f \|_{H^{s+2}}\lesssim \tau^{-2s-3}\|f\|_{H^s}.
        \end{aligned}
    \end{equation}
    For fixed $T_1<T_2<T$, we have $(T-T_2)^{\frac{1}{3}}<\tau<(T-T_1)^{\frac{1}{3}}$. 
    Then using ~\eqref{eq:RN-decay}, ~\eqref{eq:Lk-Lkveps}, ~\eqref{eq:coercive-Lk},\\ ~\eqref{eq:ap-def}, ~\eqref{eq:FG-assumption}, ~\eqref{eq:operator-L-1-P} we have 
    \begin{equation}
        \begin{aligned}
            \|\calF_{\rm{L}}^W\|_{H^{s+2}}\leq & C_{T_2} \left(\|\ti W\|_{H^{s+2}}+\|\p_t\ti V\|_{H^s}+\|\ti V\|_{H^{s+1}}+|\al|+|\be|\right),\\
            \|\calF_{\rm{L}}^V\|_{H^{s}}\leq & C_{T_2}\left(\|\ti W\|_{H^{s+2}}+\|\ti V\|_{H^{s+1}}+|\al|\right),\\
            |\calF_{\rm{L}}^{\al}|+|\calF_{\rm{L}}^{\be}|\leq & C_{T_2} \left(\|\ti W\|_{L^2}+\|\p_t\ti V\|_{L^2}+\|\ti V\|_{L^2}+|\al|+|\be|\right),\\
            \|\calF_{0}^W\|_{H^{s+2}}\leq & C_{T_2} \left(\|F\|_{H^{s+1}}+\|\p_t F\|_{H^s} \right),\\
            \|\calF_{0}^V\|_{H^{s}}\leq & C_{T_2}\|F\|_{H^{s+1}}+\|G\|_{H^{s}} ,\\
            |\calF_{0}^{\al}|+|\calF_{0}^{\be}|\leq & C_{T_2} \|F\|_{L^2} ,
        \end{aligned}
    \end{equation}
    where $C_{T_2}>0$ increases as $T_2$ increases to  $T$.

    Denote by
    \begin{equation}
        \begin{aligned}
            C_1=&\|\vec \Phi(\vec 0)\|_{X_{T_1,T_2}}\\
            \leq&4\left(\|\ti W_0\|_{H^{s+2}}+\|\ti V_0\|_{H^{s+1}}+\|\ti V_1\|_{H^s}+|\al_0|+|\be_0|\right)\\
            &+4C_{T_2}(T_2-T_1)\sup_{T_1\leq t\leq T_2 }\left(\|F(t)\|_{H^{s+1}}+\|\p_t F(t)\|_{H^s} +\|G(t)\|_{H^s}  \right).
        \end{aligned}
    \end{equation}
    Using that $\frac{\sin \sqrt{-\De_r}}{ \sqrt{-\De_r}}$ is bounded from $H^s$ to $H^{s+1}$, there holds
    \begin{equation}
        \begin{aligned}
            \|\vec \Phi(\ti W,\ti V,\al,\be)\|_{X_{T_1,T_2}}\leq C_1+4C_{T_2}(T_2-T_1)\| (\ti W,\ti V,\al,\be )\|_{X_{T_1,T_2}}.
        \end{aligned}
    \end{equation}
    Since $\calF_{\rm{L}}^*$ is $\R$-linear in $(\ti W,\ti V,\al,\be)$, and $\calF_0^*$ is independent of $(\ti W,\ti V,\al,\be)$, we have
    \begin{equation}
        \begin{aligned}
            &\| \vec\Phi(\ti W_1,\ti V_1,\al_1,\be_1)-\vec\Phi(\ti W_2,\ti V_2,\al_2,\be_2) \|_{X_{T_1,T_2}}\\
            \leq & 4C_{T_2}(T_2-T_1)\| (\ti W_1-\ti W_2,\ti V_1-\ti V_2,\al_1-\al_2,\be_1-\be_2 )\|_{X_{T_1,T_2}}.
        \end{aligned}
    \end{equation}
    For $T_2-T_1<\frac{1}{8C_{T_2}}$, $\vec\Phi$ is contractive on the ball $B(\vec 0,2C_2)$ for any $C_2>C_1$, which yields the existence and uniqueness of $(\ti W,\ti V,\al,\be)$ in $X_{T_1,T_2}$. For general $T_1<T_2<T$, since $C_{T_2}>0$ increases as $T_2$ increases to  $T$,  we  divide the interval $[T_1,T_2]$ into a finite number of small intervals, each of which has length smaller that $\frac{1}{8C_{T_2}}$, and inductively apply the result proved above, which proves the existence and uniqueness of $(\ti W,\ti V,\al,\be)$ in $X_{T_1,T_2}$. The properties ~\eqref{eq:WVab-regularity} follow directly from the system ~\eqref{eq:UV-reduced}. For any operator $\calP$ such that $\calP_M\calP=\calP$, using $\calP_M^2=\calP_M$ we have  
    \begin{equation}\label{eq:PMP}
        \begin{aligned}
            \calP_M\p_t\calP=&\p_t(\calP_M\calP)-(\p_t\calP_M)\calP\\
            =&\p_t\calP-(\calP_M\p_t\calP_M+\p_t\calP_M\calP_M)\calP\\
            =&\p_t\calP-(\p_t\calP_M)\calP-\calP_M(\p_t\calP_M)\calP_M\calP\\
            =&\p_t\calP-(\p_t\calP_M)\calP-\calP_M( \p_t\calP_M-\calP_M\p_t\calP_M) \calP\\
            =&\p_t\calP-(\p_t\calP_M)\calP.
        \end{aligned}
    \end{equation}
   Apply $(\Id-\calP_M)$ to the equation satisfied by $\ti W$, using ~\eqref{eq:PMP} we obtain
   \begin{equation}
       \begin{aligned}
           i(\Id-\calP_M)\p_t\ti W=&-i(\p_t \calP_M)(\calL_{k,\tau})^{-1}\calP_M\left( \ti \psi_{M}^{\rm{ap}}\ti V+\ti F\right)\\
           &+i(\p_t\calP_M)\calP_M\ti W+i(\p_t\calP_M)\calP_M(\calL_{k,\tau})^{-1}\calP_M\left( \ti \psi_{M}^{\rm{ap}}\ti V+\ti F\right)\\
           =&i(\p_t\calP_M)\calP_M\ti W,
       \end{aligned}
   \end{equation}
   which yields
    \begin{equation}
        \begin{aligned}
            \p_t\left(\left(\rm{Id}-\calP_M\right)\ti W\right)=&-(\p_t\calP_M)\ti W+(\p_t\calP_M)\calP_M\ti W\\
            =&-(\p_t\calP_M)(\Id-\calP_M)\ti W,
        \end{aligned}
    \end{equation}
    which combined with Gronwall's inequality and that $\calP_M|_{T_2}\ti W_0=\ti W_0$ yields $\calP_M\ti W=\ti W$.

    \end{proof}
    
    \noindent
    \textit{Step 3.} In this step we establish  energy estimates for ~\eqref{eq:UV-reduced}. Denote by $\calF_{\rm{L},1}^W=\calF^W_{\rm{L}}-\frac{1}{3\tau^4}\left( 1+ \ti N_k\left(\eta\right)\right)\ti W$, We rearrange the system
     \begin{equation}\label{eq:UV-reduced-1}
        \begin{aligned}
            \left\{
            \begin{array}{rl}
            (i\p_t+\calL_{k,\tau})\ti W= & \calF_{\rm{L},1}^W+\calF^W_0,\\
                 (\p_{tt}-\De_r)\ti V=& \calF^V_{\rm{L}}+\calF^V_0,\\
                 \p_t\al=& \calF^{\al}_{\rm{L}}+\calF^{\al}_0,\\
                 \p_t\be=& \calF^{\be}_{\rm{L}}+\calF^{\be}_0.
            \end{array}
            \right.
        \end{aligned}
    \end{equation}
    
    Note that the potential $-\frac{1}{3\tau^4} -\frac{1}{3\tau^4} \ti N_k\left(\eta\right)$ of $\calL_{k,\tau}$ becomes singular as $t\to T$, we  introduce the following rescaled Sobolev norms
    \begin{equation}
        \begin{aligned}
            \|f\|_{H^s_{\tau}}=3^{\frac{3}{4}}\tau^3\|f(\sqrt{3}\tau^2\cdot)\|_{H^s}  ,\ \|f\|_{W_{\tau}^{s,\infty}}= \|f(\sqrt{3}\tau^2\cdot)\|_{W^{s,\infty}}.
        \end{aligned}
    \end{equation}
    
    We introduce the energy functional 
    \begin{equation}
        \begin{aligned}
            \bs E=\tau^{-2}\|(-3\tau^{4}\calL_{k,\tau})^{\frac{s+2}{2}}\ti W\|^2_{L^2}+\|\ti V\|_{H^s_{\tau}}^2+\tau^{6}\left(\|\p_t\ti V\|_{H_{\tau}^s}^2+\|\na_x \ti V\|_{H_{\tau}^s}^2\right)+\tau^{-2}|\al|^2+|\be|^2.
        \end{aligned}
    \end{equation}
    Using ~\eqref{eq:coercive-Lk}, ~\eqref{eq:PN-decay} and the fact that $\calP_M\ti W=\ti W$, we have $\|(-3\tau^{4}\calL_{k,\tau})^{\frac{s+2}{2}}\ti W\|_{L^2}\sim\|\ti W\|_{H_{\tau}^{s+2}}$, with implicit constants independent of $\tau$.
    
    \begin{claim}\label{claim:energy-estimate}
        For every even integer $s\geq 2$, there exists $K_s>0$ such that for every $K\geq K_s$, there exists $M_{K,s}$ such that for every $M\geq M_{K,s}$, such that there exists $T_{1,s,K,M}<T$, such that for every $T_1$ satisfying $T_{1,s,K,M}<T_1<T$,  every $T_2$ satisfying $T_1\leq T_2<T$, and every solution $(\ti W,\ti V,\al,\be)$ to ~\eqref{eq:UV-reduced-1} define on $T_1\leq t\leq T_2$ given in Claim~\ref{claim:UV-local}, there exists a constant $C_s$ depending only on $s$ such that
        \begin{equation}\label{eq:Et-estimate}
            \begin{aligned}
                \left|\frac{\ud }{\ud t}\bs E\right|\leq C_s\frac{1}{T-t}\bs E+C_sC_0^2\tau^{12K}=C_s\frac{1}{T-t}\bs E+C_sC_0^2(T-t)^{4K}.
            \end{aligned}
        \end{equation}
    \end{claim}
    \begin{proof}[Proof of Claim ~\ref{claim:energy-estimate}]
        A direct computation yields
        \begin{equation}\label{eq:Et}
            \begin{aligned}
                \frac{\ud }{\ud t}\bs E=&\frac{2}{3}\tau^{-5}\|(-3\tau^{4}\calL_{k,\tau})^{\frac{s+2}{2}}\ti W\|^2_{L^2}+2\tau^{-2}\La (-3\tau^{4}\calL_{k,\tau})^{\frac{s+2}{2}}\p_t\ti W,(-3\tau^{4}\calL_{k,\tau})^{\frac{s+2}{2}} \ti W\Ra\\
                &+2\tau^{-2}\La \left(\p_t(-3\tau^{4}\calL_{k,\tau})^{\frac{s+2}{2}}\right) \ti W, (-3\tau^{4}\calL_{k,\tau})^{\frac{s+2}{2}}\ti W\Ra\\
                &-\frac{4}{3}\tau^{-3}\sum_{|a|\leq s}|a|3^{|a|}\tau^{4|a|}\|\p_x^a\ti V\|_{L^2}^2+2\sum_{|a|\leq s}3^{|a|}\tau^{4|a|}\La \p_t\p_x^{a}\ti V,\p_x^{a}\ti V\Ra\\
                &-\tau^3\sum_{|a|\leq s}\left(2+\frac{4}{3}|a|\right)3^{|a|}\tau^{4|a|}\left(\|\p_x^a\p_t\ti V\|_{L^2}^2+\|\p_x^a\na_x \ti V\|_{L^2}^2\right)\\
                &+2\tau^6\left(\sum_{|a|\leq s}3^{|a|}\tau^{4|a|}\left(\La \p_t\p_x^{a}\p_t\ti V,\p_x^{a}\p_t\ti V\Ra+\La \p_t\p_x^{a}\na_x \ti V,\p_x^{a}\na_x \ti V\Ra\right)\right)\\
                &+\frac{2}{3}\tau^{-5}|\al|^2+2\tau^{-2}\al\p_t\al+2\be\p_t\be.
            \end{aligned}
        \end{equation}
        The terms where the time derivatives do not hit on $(\ti W,\ti V,\p_t\ti V,\al,\be)$ can be estimated easily to obtain the bound $\tau^{-3}\bs E$. We estimate the remaining terms using ~\eqref{eq:UV-reduced-1} and Cauchy inequality to obtain
        \begin{equation}\label{eq:E1}
            \begin{aligned}
                \left|\frac{\ud }{\ud t}\bs E\right|\lesssim &  \tau^{-3} \bs E+\tau^3\Big( \tau^{-2}\left(\|\calF_{\rm{L,1}}^W\|^2_{H^{s+2}_{\tau}}+\|\calF_0^W\|^2_{H^{s+2}_{\tau}}\right)+\|\p_t\ti V\|_{H_{\tau}^s}^2\\
                &+\tau^6\left(  \|\calF_{\rm{L}}^V\|^2_{H^{s}_{\tau}}+\|\calF_0^V\|^2_{H^{s}_{\tau}}   \right)+\tau^{-2}\left(|\calF_{\rm{L}}^{\al}|^2+|\calF_{0}^{\al}|^2 \right)+\left(|\calF_{\rm{L}}^{\be}|^2+|\calF_{0}^{\be}|^2 \right)\Big).
            \end{aligned}
        \end{equation}

        \smallskip
        
        We first estimate $\|\calF_{\rm{L,1}}^W\|^2_{H^{s+2}_{\tau}}+\|\calF_0^W\|^2_{H^{s+2}_{\tau}} $ and $\|\calF_{\rm{L}}^V\|^2_{H^{s}_{\tau}}+\|\calF_0^V\|^2_{H^{s}_{\tau}}$.  Using ~\eqref{eq:LP-t} and ~\eqref{eq:Lk-Lkveps} we obtain that 
        \begin{align}
           \|\left(\p_t((\calL_{k,\tau})^{-1}\calP_M)\right)f \|_{H^{s+2}_{\tau}}\lesssim \tau \|f\|_{H_{\tau}^s}.
        \end{align}
        Using ~\eqref{eq:ap-def}, ~\eqref{eq:nre-initial}, ~\eqref{eq:nre-express}, ~\eqref{eq:nre-diff-1} and that $N_1=0$ we obtain
        \begin{align}
            \|\ti\calV_M\|_{W_{\tau}^{s+2,\infty}}\leq&\left\|  \frac{\chi_1}{3\tau^4}\sum_{m=2}^M\tau^mN_m(\log\tau,\eta)   \right\|_{W_{\tau}^{s+2,\infty}}+\left\|\left(1-\chi_1\right)\left(\frac{1}{3\tau^4}\ti N_k(\eta)-n_M^{\rm{re}}\right)\right\|_{W_{\tau}^{s+2,\infty}}\\
            &+\left\|\left(1-\chi_1\right)\left(n_M^{\rm{re}}-n^{\rm{re}}\right)\right\|_{W_{\tau}^{s+2,\infty}}\\
            \lesssim & \frac{1}{\tau^3},
        \end{align}
        where the implicit constants only depend on $s$ if we take $T_1$ sufficiently close to $T$.
        Then using ~\eqref{eq:ap-def}, ~\eqref{eq:PN-decay} and the condition ~\eqref{eq:FG} we have
        \begin{align}\label{eq:FWV}
            \|\calF_{\rm{L,1}}^W\|_{H^{s+2}_{\tau}}+\|\calF_0^W\|_{H^{s+2}_{\tau}}  \lesssim & \tau^{-2}\|\ti V\|_{H_{\tau}^{s}}+\tau\|\p_t\ti V\|_{H_{\tau}^{s}}+\tau^{-3}\|\ti W\|_{H_{\tau}^{s+2}}\\
            &+\left(|\al|+|\be|\right)\tau^{\frac{1}{4}M-8-2s}+C_0\tau^{6K}\\
            \|\calF_{\rm{L}}^V\|_{H^{s}_{\tau}}+\|\calF_0^V\|_{H^{s}_{\tau}}  \lesssim & \tau^{-7}\|\ti W\|_{H_{\tau}^{s+2}}+\tau^{-6}\|\ti V\|_{H_{\tau}^{s}}+\tau^{-7}|\al|+C_0\left(\tau^{-3}+1\right)\tau^{6K}.
        \end{align}

        \smallskip
        
        Next we estimate $|\calF_{\rm{L}}^{\al}|^2+|\calF_{0}^{\al}|^2 $and $|\calF_{\rm{L}}^{\be}|^2+|\calF_{0}^{\be}|^2$. We consider the numerator of $\calF_{\rm{L}}^{\al}$. The components containing $\ti W+(\calL_{k,\tau})^{-1}\left(\calP_M\left(\ti\psi_M^{\rm{ap}}\ti V\right)\right)$ or $(\al+i\be)\ti\calE_{\psi,M}^{\rm{ap}}$ can be estimated directly to obtain the bound $\tau^{-3}\|\ti W\|_{L^2}+\tau^{-2}\|\ti V\|_{L^2}+\left(|\al|+|\be|\right)\tau^{\frac{M}{4}-4}$. The remaining component is
        \begin{align}
                \La \ti \psi_M^{\rm{ap}},\tau^{-3}\ti P_k (\eta)\Ra    \La  \ti\psi_M^{\rm{ap}}\ti V ,i\tau^{-3}\ti P_k (\eta)  \Ra  +\La i\ti \psi_M^{\rm{ap}}, \tau^{-3}\ti P_k (\eta)\Ra   \La  \ti\psi_M^{\rm{ap}}\ti V  ,\tau^{-3}\ti P_k (\eta)  \Ra.   
        \end{align}
         Denote by $\te_M^{\rm{co}}=\sum_{m=0}^{M-1}\tau^m\te_m(\log\tau)$ and $\ti\psi_{M,1}^{\rm{ap}}=\tau^{-3}e^{ -i\te_M^{\rm{co}}(\tau)} \ti P_k (\eta) $, we rearrange this component as 
         \begin{align}
             &\La \ti \psi_M^{\rm{ap}},\tau^{-3}\ti P_k (\eta)\Ra    \La  \ti\psi_M^{\rm{ap}}\ti V ,i\tau^{-3}\ti P_k (\eta)  \Ra  +\La i\ti \psi_M^{\rm{ap}}, \tau^{-3}\ti P_k (\eta)\Ra   \La  \ti\psi_M^{\rm{ap}}\ti V  ,\tau^{-3}\ti P_k (\eta)  \Ra\\
             =&\La \ti \psi_M^{\rm{ap}},\tau^{-3}\ti P_k (\eta)\Ra    \La  \ti\psi_M^{\rm{ap}}\ti V ,i\tau^{-3}\ti P_k (\eta)  \Ra  +\La i\ti \psi_M^{\rm{ap}}, \tau^{-3}\ti P_k (\eta)\Ra   \La  \ti\psi_M^{\rm{ap}}\ti V  ,\tau^{-3}\ti P_k (\eta)  \Ra\\
             &-\La \ti\psi_{M,1}^{\rm{ap}},\tau^{-3}\ti P_k (\eta)\Ra    \La  \ti\psi_{M,1}^{\rm{ap}}\ti V ,i\tau^{-3}\ti P_k (\eta)  \Ra  -\La i\ti\psi_{M,1}^{\rm{ap}}, \tau^{-3}\ti P_k (\eta)\Ra   \La \ti\psi_{M,1}^{\rm{ap}}\ti V  ,\tau^{-3}\ti P_k (\eta)  \Ra\\
             =&\La \ti \psi_M^{\rm{ap}}-\ti\psi_{M,1}^{\rm{ap}},\tau^{-3}\ti P_k (\eta)\Ra    \La  \ti\psi_M^{\rm{ap}}\ti V ,i\tau^{-3}\ti P_k (\eta)  \Ra +\La \ti \psi_{M,1}^{\rm{ap}},\tau^{-3}\ti P_k (\eta)\Ra    \La  \left(\ti\psi_M^{\rm{ap}}-\ti\psi_{M,1}^{\rm{ap}}\right)\ti V ,i\tau^{-3}\ti P_k (\eta)  \Ra \\
             &+\La i\left(\ti \psi_M^{\rm{ap}}-\ti\psi_{M,1}^{\rm{ap}}\right),\tau^{-3}\ti P_k (\eta)\Ra    \La  \ti\psi_M^{\rm{ap}}\ti V , \tau^{-3}\ti P_k (\eta)  \Ra +\La i\ti \psi_{M,1}^{\rm{ap}},\tau^{-3}\ti P_k (\eta)\Ra    \La  \left(\ti\psi_M^{\rm{ap}}-\ti\psi_{M,1}^{\rm{ap}}\right)\ti V , \tau^{-3}\ti P_k (\eta)  \Ra. 
         \end{align}
        Since $R_1=0$, we have
        \begin{align}
            \ti \psi_M^{\rm{ap}}-\ti\psi_{M,1}^{\rm{ap}}=&\left(\chi_1-1\right)e^{-i\tau^{-1}}\psi_M^{\rm{co}}\\
            &+\frac{1}{\tau^3}\left(e^{-i\frac{1}{2}\tau\eta^2}-1\right)e^{-i\te_M^{\rm{co}}(\tau)}\sum_{m=0}^M\tau^mR_m(\log\tau,\eta)\\
            &+ \frac{1}{\tau^3}e^{-i\te_M^{\rm{co}}(\tau)}\sum_{m=2}^M \tau^mR_m(\log\tau,\eta),
        \end{align}
        which yields
        \begin{align}
            \|\ti \psi_M^{\rm{ap}}-\ti\psi_{M,1}^{\rm{ap}}\|_{L^2}\lesssim \tau,
        \end{align}
        where the implicit constant  only depends on $s$ if we take $T_1$ sufficiently close to $T$.
        Then using Cauchy inequality we have
        \begin{equation}
            \begin{aligned}
                &\left|\La \ti \psi_M^{\rm{ap}},\tau^{-3}\ti P_k (\eta)\Ra    \La  \ti\psi_M^{\rm{ap}}\ti V ,i\tau^{-3}\ti P_k (\eta)  \Ra  +\La i\ti \psi_M^{\rm{ap}}, \tau^{-3}\ti P_k (\eta)\Ra   \La  \ti\psi_M^{\rm{ap}}\ti V  ,\tau^{-3}\ti P_k (\eta)  \Ra\right| \lesssim  \tau^{-2}\|\ti V\|_{L^2}.
            \end{aligned}
        \end{equation}
        Therefore
        \begin{align}\label{eq:Fal}
            |\calF_{\rm{L}}^{\al}|\lesssim \tau^{-3}\|\ti W\|_{L^2}+\tau^{-2}\|\ti V\|_{L^2}+\left(|\al|+|\be|\right)\tau^{\frac{M}{4}-4}.
        \end{align}
        For $ \calF_{\rm{L}}^{\be} $, there is no such cancellation, and we estimate it directly to obtain the bound
        \begin{align}\label{eq:Fbe}
            |\calF_{\rm{L}}^{\be}|\lesssim \tau^{-3}\|\ti W\|_{L^2}+\tau^{-3}\|\ti V\|_{L^2}+\left(|\al|+|\be|\right)\tau^{\frac{M}{4}-4}.
        \end{align}
        The terms $\calF_0^{\al}$ and $\calF_0^{\be}$ can be estimated directly using ~\eqref{eq:FG} to obtain the bound
        \begin{align}\label{eq:Falbe0}
            |\calF_0^{\al}|+|\calF_0^{\be}|\lesssim& C_0\tau^{6K}.
        \end{align}

        Combining the estimates ~\eqref{eq:E1}, ~\eqref{eq:FWV}, ~\eqref{eq:Fal}, ~\eqref{eq:Fbe}, ~\eqref{eq:Falbe0} yields that when $K$ is sufficiently large, $M\geq 12s+8K+48$, and $T_1$ sufficiently close to $T$, there holds
        \begin{align}
            \left|\frac{\ud }{\ud t}\bs E\right|\leq C_s\frac{1}{T-t}\bs E+C_sC_0^2\tau^{12K}=C_s\frac{1}{T-t}\bs E+C_sC_0^2(T-t)^{4K},
        \end{align}
        for some constant $C_s$ only depending on $s$.
        
    \end{proof}

    \noindent
    \textit{Step 4.} In this step we complete the construction and prove the uniqueness. Take $T_1<T$ sufficiently close to $T$ such that Claim~\ref{claim:UV-local} and Claim~\ref{claim:energy-estimate} hold. Take a strictly increasing sequence $\{T_{2,\ell}\}_{\ell\geq 1}$ such that $T_{2,\ell}\to T^-$. For each $T_{2,\ell}$, we consider the initial data at $T_{2,\ell}$ as
    \begin{align}
        \ti W_{0,\ell}=-(\calL_{k,\tau})^{-1}\left(\calP_M\ti F\right)\Big|_{T_{2,\ell}},\ \ti V_{0,\ell}=0,\ \ti V_{1,\ell}=0,\ \al_{0,\ell}=0,\ \be_{0,\ell}=0.
    \end{align}
    Then Claim ~\ref{claim:UV-local} yields the existence and uniqueness of $(\ti W_{\ell},\ti V_{\ell},\al_{\ell},\be_{\ell})$ on $[T_1,T_{2,\ell}]$ for each $\ell\geq 1$. Consequently, the relative $(U_{\ell},V_{\ell})$ satisfies $U_{\ell}(T_{2,\ell})=V_{\ell}(T_{2,\ell})=0$ and solves ~\eqref{eq:UV-FG}. 
    \begin{claim}\label{claim:UV-conv}
        The pairs of solutions $\{(U_{\ell},V_{\ell})\}_{\ell\geq 1}$ converge in 
        \begin{equation}
            \begin{aligned}
                & \Big( C_{\rm{loc}}([T_1,T);H^{s+2}(\R^3,\C))\cap C_{\rm{loc}}^1([T_1,T);H^{s}(\R^3,\C))\Big)\\
            \times& \Big( C_{\rm{loc}}([T_1,T);H^{s+1}(\R^3,\R))\cap C_{\rm{loc}}^1([T_1,T);H^{s}(\R^3,\R))\Big)
            \end{aligned}
        \end{equation}
    \end{claim}
    \begin{proof}[Proof of Claim~\ref{claim:UV-conv}]
        Fix $T_2\in[T_1,T)$. Take $\ell_0\in \Z_{\geq 1}$ such that $T_{2,\ell}>T_2$ for any $\ell\geq \ell_0$. Denote by $\bs E_{\ell}$ the related energy for $(\ti W_{\ell},\ti V_{\ell},\al_{\ell},\be_{\ell})$. The estimate ~\eqref{eq:Et-estimate} yields that  for $4K>C_s+3$, there holds
        \begin{align}\label{eq:Eell}
            \bs E_{\ell}(t)&\lesssim_{s,K,M } \left(\frac{T-t}{T-T_{2,\ell}}\right)^{C_s}\bs E_{\ell}(T_{2,\ell}) +C_0^2(T-t)^{4K+1}\\
            &\lesssim_{s,K} C_0^2(T-t)^{C_s}(T-T_{2,\ell})^{-C_s-2+4K}+C_0^2(T-t)^{4K+1} .
        \end{align}
        
        For $\ell_0\leq\ell_1<\ell_2<\infty$, we denote by $f_{\ell_1,\ell_2}=f_{\ell_1}-f_{\ell_2}$ for $f=U,V,\ti W,\ti V,\al,\be$. Then on $[T_1,T_{2,\ell_1}]$ there holds
        \begin{equation}\label{eq:UV-reduced-2}
        \begin{aligned}
            \left\{
            \begin{array}{rl}
            (i\p_t+\calL_{k,\tau})\ti W_{\ell_1,\ell_2}= & \calF_{\rm{L},1,\ell_1,\ell_2}^W,\\
                 (\p_{tt}-\De_r)\ti V_{\ell_1,\ell_2}=& \calF^V_{\rm{L},\ell_1,\ell_2},\\
                 \p_t\al_{\ell_1,\ell_2}=& \calF^{\al}_{\rm{L},\ell_1,\ell_2},\\
                 \p_t\be_{\ell_1,\ell_2}=& \calF^{\be}_{\rm{L},\ell_1,\ell_2}.
            \end{array}
            \right.
        \end{aligned}
    \end{equation}
        where $ \calF_{\rm{L},1,\ell_1,\ell_2}^W, \calF^V_{\rm{L},\ell_1,\ell_2},\calF^{\al}_{\rm{L},\ell_1,\ell_2},\calF^{\be}_{\rm{L},\ell_1,\ell_2}$ are the difference of the related $\calF_{\rm{L},\ell_j}^*$. Denote by 
        \begin{align}
            \bs E_{\ell_1,\ell_2}=&\tau^{-2}\|(-3\tau^{4}\calL_{k,\tau})^{\frac{s+2}{2}}\ti W_{\ell_1,\ell_2}\|^2_{L^2}+\|\ti V_{\ell_1,\ell_2}\|_{H^s_{\tau}}^2+\tau^{6}\left(\|\p_t\ti V_{\ell_1,\ell_2}\|_{H_{\tau}^s}^2+\|\na_x \ti V_{\ell_1,\ell_2}\|_{H_{\tau}^s}^2\right)\\
            &+\tau^{-2}|\al_{\ell_1,\ell_2}|^2+|\be_{\ell_1,\ell_2}|^2.
        \end{align}
        Using the same argument as in the proof of Claim~\ref{claim:energy-estimate} we obtain
        \begin{align}
            \left|\frac{\ud }{\ud t}\bs E_{\ell_1,\ell_2}\right|\leq C_s \frac{1}{T-t}\bs E_{\ell_1,\ell_2},
        \end{align}
        which combined with ~\eqref{eq:Eell} for $\ell_2$ yields
        \begin{align}
            \bs E_{\ell_1,\ell_2}(t)&\lesssim_{s,K,M } \left(\frac{T-t}{T-T_{2,\ell_1}}\right)^{C_s}\bs E_{\ell_1}(T_{2,\ell_1})+\left(\frac{T-t}{T-T_{2,\ell_1}}\right)^{C_s}\bs E_{\ell_2}(T_{2,\ell_1}) \\
            &\lesssim C_0^2(T-t)^{C_s}\left((T-T_{2,\ell_1})^{-C_s-2+4K}+(T-T_{2,\ell_2})^{4K-2-C_s} \right),
        \end{align}
        which tends to zero uniformly for $t\in[T_1,T_2]$ as $\ell_1,\ell_2\to \infty$. Note that the term $F$ appearing in the expression of $\ti W$ cancels when taking the difference of $\ti W_{\ell_1}-\ti W_{\ell_2}$ and $\tau=(T-t)^{\frac{1}{3}}$ is uniformly bounded away from zero for $t\in [T_1,T_2]$, we have 
        \begin{align}
             &\|U_{\ell_1,\ell_2}(t)\|_{H^{s+2}(\R^3,\C)}^2  +\|\p_t U_{\ell_1,\ell_2}(t)\|_{H^{s}(\R^3,\C)}^2 
             +\|V_{\ell_1,\ell_2}(t)\|_{H^{s+1}(\R^3,\R)}^2+\|\p_t V_{\ell_1,\ell_2}(t)\|_{H^{s}(\R^3,\R)}^2\\
             \lesssim &  \| \ti W_{\ell_1,\ell_2}(t)\|^2_{H^{s+2}}+\|\ti V_{\ell_1,\ell_2}(t)\|_{H^s_{\tau}}^2+ \|\p_t\ti V_{\ell_1,\ell_2}(t)\|_{H_{\tau}^s}^2+\|\na_x \ti V_{\ell_1,\ell_2}(t)\|_{H_{\tau}^s}^2 + |\al_{\ell_1,\ell_2}(t)|^2+|\be_{\ell_1,\ell_2}(t)|^2\\
             \lesssim & \bs E_{\ell_1,\ell_2}(t).
        \end{align}
         Then $\{(U_{\ell},V_{\ell})\}_{\ell\geq 1}$ is Cauchy and hence converges in
         \begin{equation}
            \begin{aligned}
                & \Big( C([T_1,T_2];H^{s+2}(\R^3,\C))\cap C^1([T_1,T_2];H^{s}(\R^3,\C))\Big)\\
            \times& \Big( C([T_1,T_2];H^{s+1}(\R^3,\R))\cap C^1([T_1,T_2];H^{s}(\R^3,\R))\Big).
            \end{aligned}
        \end{equation} 
        For $T_1<T_{2}<T_2'<T$, the limits are the same  on $[T_1,T_2]$ by definition, which yields the desired result.
    \end{proof}
    
    Denote the limit of $\{(U_{\ell},V_{\ell})\}_{\ell\geq 1}$ by $(U_{\infty},V_{\infty})$. Denote by $(\ti W_{\infty},\ti V_{\infty},\al_{\infty},\be_{\infty})$ and $\bs E_{\infty}$ the related functions and quantities. Then $(U_{\infty},V_{\infty})$ solves ~\eqref{eq:UV-FG} and the inequality
    \begin{align}
        \bs E_{\infty}(t)\lesssim_{s,K,M } C_0^2(T-t)^{4K+1}
    \end{align}
    holds uniformly for $t\in[T_1,T)$. Then we have
    \begin{equation}
        \begin{aligned}
            &(T-t)^{-K}\left(\|U_{\infty}\|_{H^{s+2}}+\|\p_t U_{\infty}\|_{H^{s}}+\|V_{\infty}\|_{H^{s+1}}+\|\p_t V_{\infty}\|_{H^{s}}\right)\\
            \lesssim &\  \tau^{-3K-2s-4}\bs E^{\frac{1}{2}}_{\infty}(t)+  C_0(T-t)^{K-\frac{3}{4}-s}\\
            \lesssim & C_0(T-t)^{3K-\frac{2}{3}s}+C_0(T-t)^{K-s-\frac{4}{3}}.
        \end{aligned}
    \end{equation}
    Therefore $(U_{\infty},V_{\infty})\in\calH_K^s$ for $K\geq 2s+2$, and ~\eqref{eq:UV-estimate} is proved.

    Finally we prove the uniqueness. It suffices to prove that the solution $(U,V)$ to ~\eqref{eq:UV-FG} with $F=G=0$ vanish. Consider the related $(\ti W,\ti V,\al,\be)$ and $\bs E$. Using the same argument as in Claim~\ref{claim:energy-estimate} we have
    \begin{align}
        \left|\frac{\ud}{\ud t}\bs E\right|\leq C_s\frac{1}{T-t}\bs E.
    \end{align}
    Fix $t\in [T_1,T)$ and take $T_2$ such that $t<T_2<T$,  we have
    \begin{align}
        \bs E(t)\leq &  \left(\frac{T-t}{T-T_2}\right)^{C_s}\bs E(T_2)\\
        \lesssim & \|(U,V)\|_{\calH_K^s}^2(T-t)^{C_s}(T-T_2)^{2K-C_s-1},
    \end{align}
    where the implicit constants are independent of $T_2$. Then let $T_2\to T^-$ we obtain that for $2K\geq C_s+2$, $E(t)\equiv 0$. Thus the proof is finished.
\end{proof}

\bigskip

Lemma ~\ref{lem:UV} leads to the map defined for $(F,G)$ satisfying ~\eqref{eq:FG} as $\calT(F,G)\coloneqq(U,V)$, where $(U,V)$ is the solution to ~\eqref{eq:UV-FG} given by Lemma~\ref{lem:UV}. Since $(uv,\De(|u|^2))$ and $(\calE_{\psi,M}^{\rm{ap}},\calE_{n,M}^{\rm{ap}})$ satisfy ~\eqref{eq:FG} for $M\geq 8K+12s+48$, the solution to ~\eqref{eq:UV} can be expressed as
\begin{align}
    (U,V)=\calT(uv,\De(|u|^2))-\calT(\calE_{\psi,M}^{\rm{ap}},\calE_{n,M}^{\rm{ap}}).
\end{align}
We denote by $\calT_{\rm{nl}}(u,v)=\calT(uv,\De(|u|^2))$, $\calT_0=\calT(\calE_{\psi,M}^{\rm{ap}},\calE_{n,M}^{\rm{ap}})$, $\ti \calT(u,v)=\calT_{\rm{nl}}(u,v)-\calT_0$ for convenience.
 
\subsection{Nonlinear estimates}\label{subsec:nl}

In this section we establish the estimates on $\calT_{\rm{nl}}(u,v)$.

\begin{lem}\label{lem:nl}
    Under the assumption on $s,K,M,T_1$ as in Lemma~\ref{lem:UV}, there holds
    \begin{equation}\label{eq:nl-estimate}
        \begin{aligned}
             &\|\calT_{\rm{nl}}(u_1,v_1)-\calT_{\rm{nl}}(u_2,v_2)\|_{\calH_K^s}\\
             \lesssim_{s,K,M}& (T-T_1)^{K-s-2}\left(\|(u_1,v_1)\|_{\calH_K^s} +\|(u_2,v_2)\|_{\calH_K^s} \right)\|(u_1-u_2,v_1-v_2)\|_{\calH_K^s}.
        \end{aligned}
    \end{equation}
\end{lem}
\begin{proof}[Proof of Lemma~\ref{lem:nl}]
    There holds
    \begin{equation}
        \begin{aligned}
            \|u_1v_1-u_2v_2\|_{H^{s+1}}\lesssim & \|(u_1-u_2)v_1\|_{H^{s+1}}+\|u_2(v_1-v_2)\|_{H^{s+1}}\\
            \lesssim & \|u_1-u_2\|_{H^{s+1}}\|v_1\|_{H^{s+1}}+\|u_2\|_{H^{s+1}}\|v_1-v_2\|_{H^{s+1}}\\
            \lesssim & (T-t)^{2K}\left(\|(u_1,v_1)\|_{\calH_K^s} +\|(u_2,v_2)\|_{\calH_K^s} \right)\|(u_1-u_2,v_1-v_2)\|_{\calH_K^s},\\
            \|\p_t(u_1v_1-u_2v_2)\|_{H^s}\lesssim & (T-t)^{2K}\left(\|(u_1,v_1)\|_{\calH_K^s} +\|(u_2,v_2)\|_{\calH_K^s} \right)\|(u_1-u_2,v_1-v_2)\|_{\calH_K^s}, 
        \end{aligned}
    \end{equation}
    and
    \begin{equation}
        \begin{aligned}
            \| \De_r(|u_1|^2)-\De_r(|u_2|^2)\|_{H^{s }}\lesssim & \|u_1-u_2\|_{H^s}\|u_1\|_{H^{s+2}}+\|u_2\|_{H^s}\|u_1-u_2\|_{H^{s+2}}\\
            &+\|u_1-u_2\|_{H^{s+1}}\left(\|u_1\|_{H^{s+1}}+\|u_2\|_{H^{s+1}}\right)    \\
            \lesssim & (T-t)^{2K}\left(\|(u_1,v_1)\|_{\calH_K^s} +\|(u_2,v_2)\|_{\calH_K^s} \right)\|(u_1-u_2,v_1-v_2)\|_{\calH_K^s}
        \end{aligned}
    \end{equation}
    Then we take $C_0=\left(\|(u_1,v_1)\|_{\calH_K^s}+\|(u_2,v_2)\|_{\calH_K^s}\right)\|(u_1-u_2,v_1-v_2)\|_{\calH_K^s}$ in ~\eqref{eq:FG} and use ~\eqref{eq:UV-estimate} to obtain ~\eqref{eq:nl-estimate}.
\end{proof}

\subsection{Construction of the solution with finite regularity}\label{subsec:finite-regularity}

In this section we construct the solution of finite regularity. 

\begin{prop}
    \label{prop:fixed-point}
    Under the assumption of $s,K,M,T_1$ in Lemma~\ref{lem:UV},  there exists a solution $(u,v)\in \calH_{K}^s([T_1,T))$ to 
    \begin{equation}\label{eq:uv-1}
        \begin{aligned}
            \left\{
         \begin{array}{rl}
              i\p_t u+\De_r u-n_M^{\rm{ap}}u-\psi_{M}^{\rm{ap}}v &=uv-\calE_{\psi,M}^{\rm{ap}},  \\
              \p_{tt}v-\De_r v-2\Re\De(\overline{\psi_{M}^{\rm{ap}}}u) &=\De(|u|^2)-\calE_{n,M}^{\rm{ap}}.
         \end{array}
    \right.
        \end{aligned}
    \end{equation}
    Moreover, such a solution is unique in the closed unit ball of $\calH_K^s$ and satisfies $\|(u,v)\|_{\calH_{K}^s([T_1,T))}<1$. 
\end{prop}
\begin{proof}[Proof of Proposition ~\ref{prop:fixed-point}]
    Combining ~\eqref{eq:source-estimate} and ~\eqref{eq:UV-estimate}, there exists a constant $C_{s,K,M}$ depending only on $s,K,M$ such that
    \begin{align}\label{eq:T-estimate}
        \|\ti \calT(u,v)\|_{\calH_K^s}\leq C_{s,K,M}(T-T_1)^{K-s-2}\left( \|(u ,v )\|_{\calH_K^s}^2 +(T-T_1)^{\frac{M}{4}-6-2s-2K}\right).
    \end{align}
    Fix $s,K,M$ such that $K=K(s)\gg 1$ and $M=M(s,K)\gg K$ satisfying all the requirements in Lemma ~\ref{lem:UV}. Then for $T_1<T$ sufficiently close to $T$, ~\eqref{eq:nl-estimate} and ~\eqref{eq:T-estimate} yields that $\ti \calT$ is contractive on the closed unit ball of $\calH_K^s$, which yields the existence and uniqueness of the solution in the closed unit ball. For $T_1$ sufficiently close to $T$, ~\eqref{eq:T-estimate} yields $\|(u,v)\|_{\calH_{K}^s([T_1,T))}<1$.
\end{proof}

\subsection{Construction of smooth solution}\label{subsec:smooth}

We prove that the solutions to ~\eqref{eq:Za} we constructed in Proposition~\ref{prop:fixed-point} are actually the same.

\begin{prop}\label{prop:smooth}
    There exists $T_0<T$ such that there exists a solution  to ~\eqref{eq:Za} defined on $[T_0,T)$ such that for any even integer $s\geq 2$ and integer $K\geq K_s\gg 1$  for some $K_s$ depending on $s$, there exists an integer $M$ and a time $T_1$ satisfying $T_0\leq T_1<T$ such that 
    \begin{align}
        (\psi,n)=(\psi_{M}^{\rm{ap}},n_M^{\rm{ap}})+(u,v), \quad  \forall \ t\in[T_1,T),
    \end{align}
    where $(\psi_{M}^{\rm{ap}},n_M^{\rm{ap}})$ is defined in ~\eqref{eq:ap-def} and $(u,v)$ is given in Proposition~\ref{prop:fixed-point}.
\end{prop}

\begin{proof}[Proof of Proposition~\ref{prop:smooth}]
    Fix an even integer $s_0\gg 1$, Lemma~\ref{lem:UV} and Proposition~\ref{prop:fixed-point} yields the existence of integers $K_0$, $M_0\geq 4K_0+16s_0+48$ and a time $T_0<T$ such that  there exists a solution $(\psi_0,n_0)$ to ~\eqref{eq:Za} taking the form
    \begin{align}
        (\psi_0,n_0)=(\psi_{M_0}^{\rm{ap}},n_{M_0}^{\rm{ap}})+(u_0,v_0), \quad  \forall \ t\in[T_0,T)
    \end{align}
    with $(u_0,v_0)\in\calH_{K_0}^{s_0}([T_0,T))$. 
    
    For any even integers $s_1\geq s_0$ and $K_1\geq K_0+K_{s_1}$ where $K_{s_1}$ is given in Lemma~\ref{lem:UV}, fix a $M_1\geq M_0+M_{K_1,s_1}$ sufficiently large and $T_1\in (T_0,T)$ sufficiently close to $T$ so that Proposition ~\ref{prop:fixed-point} gives a solution $(\psi_1,n_1)=(\psi_{M_1}^{\rm{ap}},n_{M_1}^{\rm{ap}})+(u_1,v_1)$ with $(u_1,v_1)\in\calH_{K_1}^{s_1}([T_1,T))$.

    We prove that there exists $T_2\in[T_1,T)$ such that $(\psi_0,n_0)=(\psi_1,n_1)$ for $T_2\leq t<T$. Write
    \begin{align}
        (\psi_1,n_1)=&(\psi_{M_0}^{\rm{ap}},n_{M_0}^{\rm{ap}})+(u_1,v_1)+(\psi_{M_1}^{\rm{ap}}-\psi_{M_0}^{\rm{ap}},n_{M_1}^{\rm{ap}}-n_{M_0}^{\rm{ap}}),
    \end{align}
    then $ (u_1,v_1)+(\psi_{M_1}^{\rm{ap}}-\psi_{M_0}^{\rm{ap}},n_{M_1}^{\rm{ap}}-n_{M_0}^{\rm{ap}})$ is also the solution to ~\eqref{eq:uv-1} with $M$ replaced by $M_0$.  Lemma ~\ref{lem:ap-diff} yields that
    \begin{equation}
        \begin{aligned}
            \|(\psi_{M_0}^{\rm{ap}}-\psi_{M_1}^{\rm{ap}},n_{M_0}^{\rm{ap}}-n_{M_1}^{\rm{ap}})\|_{\calH_{K_0}^{s_0}([T_2,T))}\lesssim_{M_1}  (T-T_2)^{\frac{M_0}{4}-K_0-2s_0-4},
        \end{aligned}
    \end{equation}
    which combined with $\|(u_1,v_1)\|_{\calH_{K_0}^{s_0}([T_2,T))}\leq \|(u_1,v_1)\|_{\calH_{K_1}^{s_1}([T_1,T))}<1$ implies that $ (u_1,v_1)+(\psi_{M_1}^{\rm{ap}}-\psi_{M_0}^{\rm{ap}},n_{M_1}^{\rm{ap}}-n_{M_0}^{\rm{ap}})$ lies in the closed unit ball of $\calH_{K_0}^{s_0}([T_2,T)) $ for $T_2<T$ sufficiently close to $T$. Then using the same fixed point argument as in the proof of Proposition ~\ref{prop:fixed-point} we obtain that $(\psi_0,n_0)=(\psi_1,n_1)$ for $T_2\leq t<T$. Then taking $(\psi,n)=(\psi_0,n_0)$ yields the desired result.

\end{proof}

\begin{proof}[Proof of Theorem ~\ref{thm:main}]
    We consider $(\psi,n)$ constructed in Proposition ~\ref{prop:smooth}. After translation of time we may assume $T_0=0$. Lemma~\ref{lem:Cauchy} yields that $\psi,n\in \mathop{\cap} \limits_{s=0}^{\infty}H^s$, which implies ~\eqref{eq:main-0}.  The estimates ~\eqref{eq:main-rate}, ~\eqref{eq:main-core} and  ~\eqref{eq:main-remote} are immediate results of Proposition~\ref{prop:core}, Lemma~\ref{lem:iteration-solve}, ~\eqref{eq:core-app}, Lemma~\ref{lem:N-infity}, ~\eqref{eq:ap-def}, Proposition~\ref{prop:fixed-point} and Proposition~\ref{prop:smooth}.
\end{proof}

\bibliographystyle{plainurl}
\bibliography{Zakharov}

\end{document}